\documentclass[fleqn,11pt,letterpaper]{article}
\usepackage{amsfonts}
\usepackage{amsmath}
\usepackage[final]{showkeys}
\usepackage{amsmath,color}
\usepackage{latexsym,amscd,graphicx,dsfont,mathrsfs,amssymb}
\usepackage{appendix}
\usepackage{bm}

\newtheorem{theorem}{Theorem}
\newtheorem{lemma}{Lemma}

\newtheorem{definition}{Definition}

\begin{document}

\setcounter{section}{0}
\newcommand{\bez} {B\'{e}zier }

\renewcommand{\theequation}{\thesection.\arabic{equation}}
\renewcommand{\thefigure}{\arabic{figure}}
\renewcommand{\thetheorem}{\thesection.\arabic{theorem}}
\renewcommand{\thelemma}{\thesection.\arabic{lemma}}
\renewcommand{\thecorollary}{\thesection.\arabic{corollary}}
\renewcommand{\thedefinition}{\thesection.\arabic{definition}}
\renewcommand{\thetable}{\arabic{table}}
\renewcommand{\theremark}{\thesection.\arabic{remark}}
\renewcommand{\theexample}{\thesection.\arabic{example}}
\newcommand{\IR}{I\!\!R}

\title{Isogeometric Analysis for Surface PDEs with Extended Loop Subdivision}

\author{ {\sf Qing Pan${}^{a},$\quad Timon Rabczuk${}^{b},$\quad Gang Xu${}^{c},$\quad Chong Chen${}^{d,}$\thanks{Corresponding author. E-mail address:
chench@lsec.cc.ac.cn}}
\\[1mm]\textit{${}^{a}$ \!\!\small Key Laboratory of High Performance Computing and Stochastic Information Processing}
\\[-1mm] \textit{\small  (Ministry of Education of China), Hunan Normal University, Changsha, China}
\\[-1mm] \textit{${}^{b}$ \!\!\small Institute of Structural Mechanics, Bauhaus Universit${\ddot{a}}$t-Weimar, Weimar, Germany}
\\[-1mm] \textit{${}^{c}$ \!\!\small Department of Computer Science, Hangzhou Dianzi University, Hangzhou, China}
\\[-1mm] \textit{${}^{d}$ \!\!\small LSEC, ICMSEC, Academy of Mathematics and Systems Science, Chinese Academy of }
\\[-1mm] \textit{\small Sciences, Beijing, China}
}

\date{}
\maketitle

\begin{abstract}
We investigate the isogeometric analysis for surface PDEs based on the extended Loop subdivision approach.
The basis functions consisting of quartic box-splines corresponding to each subdivided control mesh are utilized to
represent the geometry exactly, and construct the solution space for dependent variables
as well, which is consistent with the concept of isogeometric analysis.
The subdivision process is equivalent to the $h$-refinement of NURBS-based isogeometric
analysis. The performance of the proposed method is evaluated by solving various surface PDEs, such as
surface Laplace-Beltrami harmonic/biharmonic/triharmonic equations, which are
defined on different limit surfaces of the extended Loop subdivision for different initial control meshes.
Numerical experiments demonstrate that the proposed method has desirable performance in terms of the
accuracy, convergence and computational cost for solving the above surface PDEs defined on both
open and closed surfaces. The proposed approach is proved to be second-order accuracy in the sense of $L^2$-norm by
theoretical and/or numerical results, which is also outperformed over the standard linear finite
element by several numerical comparisons.

\end{abstract}

\setcounter{theorem}{0} \setcounter{equation}{0}
\setcounter{figure}{0}
\section{Introduction}
\label{intro}
The isogeometric analysis (IGA) was introduced by Hughes et al. \cite{Hughes2006,Hughes2010} to
replace the traditional finite elements by Non-Uniform Rational B-Splines (NURBS) \cite{Piegl97,Sederberg08} or
T-splines \cite{Sederberg03,Sederberg04,Li06,Bavilevs10}. The concept of IGA shows great potential in developing
the seamless integration between computer-aided design (CAD) and finite element method (FEM), which possesses higher
accuracy than traditional FEM. The $h$-refinement, $p$-refinement and even $k$-refinement can be easily
implemented by knot insertion and/or order elevation to improve the simulation accuracy without changing
the geometry, which vanish the necessity of subsequent communication with CAD system. From the analytical
point of view, hierarchical bases are preferred. Several constructions of such spaces have been proposed,
the most common of which are (Truncated) Hierarchical B-Splines (THB) \cite{Bert2}, Locally Refined
(LR)-Splines \cite{Dokken13,Dokken14} and Polynomial/Rational Splines over Hierarchical
T-Meshes (PHT/RHT splines) \cite{NguyenJ11,NguyenH11,Nguyen15,Nguyen17,Deng08}. A review
including an efficient IGA code can be found in \cite{Timon15}. A computation framework is reused in IGA on a set of three-dimensional models with similar semantic features in \cite{XuGang17}.

On the other hand, surface PDEs have been attained extensive attentions in a variety of subjects, such as
structural mechanics, biology, fluid dynamics, and image processing, etc.
\cite{Ding09,Cho96,Cantwell14,Patelli17,Fedkiw03,Lai11,Ratz06,Dziuk08,Reuter09}. A
comprehensive review about this topic can be founded in \cite{Bartezzaghi17}.
Although the FEMs were developed to solve such kind of problems (see, for instance,
\cite{Dziuk88,Dziuk07,Olshanskii14,Mekchay11}), which cannot possess the ability to represent the geometry exactly \cite{Hughes2010}.
To overcome this deficiency, the NURBS-based IGA was proposed to solve the high-order surface PDEs in \cite{Bartezzaghi17}.
In this article, we will further study the subdivision-based IGA, and use this method to solve surface PDEs.

Subdivision is a powerful technique in surface modeling, which provides a simple and efficient method to generate smooth
surfaces with arbitrary topology structures. Moreover, it is capable of recovering sharp
features of surfaces with creases and corners. Subdivision surfaces
play a key role in computer graphics and numerical analysis. A class of piecewise-smooth surface representations
was introduced based on subdivision to reconstruct smooth surface from scattered
data \cite{Hoppeetal}. Thin-shell finite element analysis was used to describe both the geometry and associated displacement
fields \cite{Cirak}. The limit representation of Loop subdivision for triangular meshes was combined
with the diffusion model to arrive at a discretized version of the diffusion problem \cite{BajajXu03}.
Mixed finite element methods based on subdivision technology were used to
construct high-order smooth surfaces with specified boundary conditions \cite{panxu20116}. Truncated
hierarchical Catmull-Clark subdivision was developed recently to support local refinement and generalize
truncated hierarchical B-splines to arbitrary topology \cite{Jessica15}.

As a compatible technique of NURBS, subdivision surfaces are capable of the refinability of B-spline techniques.
There recently have been a few works on the application of the subdivision-based IGA approach. Volumetric IGA based
on Catmull-Clark solids was investigated in \cite{Umlauf2010}. For complex physical domain,
Powell-Sabin splines were used as IGA tools for advection-diffusion-reaction problems \cite{SMS12}. The bivariate
splines in the rational Bernstein-B\'ezier form over the triangulation were applied in \cite{NX14}. A reproducing
kernel triangular B-spline-based FEM was proposed to solve PDEs \cite{JZX13}. Collocated
isogeometric boundary element methods and unstructured analysis-suitable T-spline surfaces were coupled for
linear elastostatic problems in \cite{SELBHS13}. A new generalized surface and IGA elements having the vertices
of the irregular quad mesh through complementing bicubic splines and biquartic splines near irregularities in the mesh
layout was presented in \cite{KNP16}. A framework for geometric design and IGA on unstructured quadrilateral
meshes was proposed in \cite{TSH17}. A new type of Hermite bases for bicubic spline defined over a rectangular
mesh with arbitrary topology was investigated in the framework of IGA \cite{WMGN17}. Second-order PDEs on
surfaces with the IGA was considered in \cite{Dede15}.

\paragraph{Contributions} In the present paper, we introduce an alternative numerical method, namely,
the isogeometric analysis based on the extended Loop subdivision approach (IGA-Loop), for solving the PDEs on subdivision surfaces.
The finite elements consisting of quartic box-splines corresponding to each
subdivided control mesh are utilized to represent the geometry of interest, and construct the
solution space for dependent variables as well. The subdivision process is equivalent to
the $h$-refinement of NURBS-based isogeometric analysis. The basis functions induced by
the extended Loop subdivision possess the ability to represent the geometry exactly which
is consistent with the concept of isogeometric analysis. The performance of the proposed
method is evaluated by solving various surface PDEs, such as
surface Laplace-Beltrami harmonic/biharmonic/triharmonic equations, which are
defined on various limit surfaces of the extended Loop subdivision for different initial control meshes.
Numerical experiments demonstrate that the proposed method has desirable performance in terms of the
accuracy, convergence and computational cost for solving the above classical surface PDEs defined on both open and closed surfaces.
The proposed approach is proved to be second-order accuracy in the sense of $L^2$-norm by
theoretical and numerical results for the surface Laplace-Beltrami harmonic equation, which is outperformed over the standard linear finite
element by several numerical comparisons. Furthermore, the $L^2$-norm error with the second-order convergence rate can be
observed for the surface Laplace-Beltrami biharmonic/triharmonic equations.

\paragraph{Outline} Let us outline the content of the paper. The required mathematical preliminaries are introduced in Section \ref{math_pre}.
The isogeometric analysis based on the extended Loop subdivision is investigated in Section \ref{loopintro}.
We apply the proposed method to solve second/fourth/sixth-order surface PDEs in Section \ref{SurfacePDE}. The Section \ref{numer_ex}
presents several numerical experiments with comparison to the linear finite element method.
Finally, the Section \ref{conclu} concludes the paper.

\setcounter{theorem}{0} \setcounter{equation}{0}
\setcounter{figure}{0}
\section{Mathematical Preliminaries}
\label{math_pre}

Here we introduce some mathematical preliminaries which are requisite for treating the PDEs on surface.

\subsection{Surface Parameterization}

Let ${\cal S} = \{{\bm x}(u,v) \in \mathbb{R}^{3}: [u,v]^{{\rm T}} \in \Omega
\subset \mathbb{R}^{2} \}$ be a piece of regular smooth parametric
surface. The tangent plane and unit normal vector at ${\bm x}$
of ${\cal S}$ are defined as
\[
T_{\bm x} {\cal S} = {\rm span} \{ {\bm x}_u, {\bm x}_v\} \quad {\rm
and} \quad {\bm n} = \frac{{\bm x}_u \times {\bm x}_v}{\|
{\bm x}_u \times {\bm x}_v \| },
\]
respectively, where ${\bm x}_u$ and ${\bm x}_v$ are the coordinate tangent
vectors. The {\it first fundamental form} of $\cal S$ is
\[ I = \langle {\rm
d}{\bm x}, {\rm d}{\bm x} \rangle = g_{11} {\rm d} u {\rm d} v +
2 g_{12} {\rm d} u {\rm d} v + g_{22} {\rm d} u {\rm d} v,
\]
where the coefficients $g_{\alpha \beta} = \langle {\bm
x}_{u^{\alpha}}, {\bm x}_{u^{\beta}} \rangle$, $\alpha, \beta = 1,2$
with $(u^1, u^2) = (u, v)$.  The {\it second fundamental form} of $\cal S$ is
\[
I\!I = -\langle {\rm d}{\bm x}, {\rm d}{\bm n} \rangle = b_{11} {\rm
d}u {\rm d}u + 2 b_{12} {\rm d}u {\rm d}v + b_{22} {\rm
d}v {\rm d}v
\]
with the coefficients $b_{\alpha \beta} = \langle {\bm
x}_{u^{\alpha} u^{\beta}}, {\bm n} \rangle = - \langle {\bm
x}_{u^{\alpha}}, {\bm n}_{u^{\beta}} \rangle$.
Set
\[
g = \det[g_{\alpha \beta}], \quad [g^{\alpha \beta}] =
[g_{\alpha \beta}]^{-1} \quad {\rm
and} \quad  b = \det[b_{\alpha \beta}] \quad\text{for}\quad \alpha, \beta = 1,2.
\]

\paragraph{Area} The positive value
\begin{equation}
\label{area}
\int_{\cal S} {\rm d} A = \iint_{\Omega} |{\bm x}_u \wedge {\bm
x}_v| {\rm d} u {\rm d} v
\end{equation}
is called the {\it area} of surface $\cal S$.

\paragraph{Mean Curvature} First we introduce Weingarten transformation ${\cal W}: T_{\bm x} {\cal S}
\rightarrow T_{\bm x} {\cal S}$ which is a linear mapping, and satisfies
\[
{\cal W}({\bm x}_u) = -{\bm n}_u \quad {\rm and} \quad {\cal
W}({\bm x}_v) = -{\bm n}_v.
\]
The linear mapping can be represented by a $2 \times 2$ matrix $S =
[b_{\alpha \beta}][g^{\alpha \beta}].$ The eigenvalues $k_1$ and
$k_2$ of $S$ are the principal curvatures of $\cal S$. Their
arithmetic average is the {\it mean curvature}:
\[
H := \frac{k_1+k_2}{2} =\frac{{\rm
tr} (S)}{2} = \frac{b_{11}g_{22} - 2b_{12}g_{12} +
b_{22}g_{11}}{2g}.
\]

\subsection{Surface Differential Operators}

Let $C^1({\cal S})$ denote the set of function $f:{\cal S}
\rightarrow \mathbb{R}$ with continuous derivatives at any ${\bm x}
\in {\cal S}$. Analogously, we can define the sets $C^2({\cal S})$, $C^3({\cal S}),$ $\ldots$.
For the sufficiently differentiable $f({\bm
x}) = f({\bm x}(u, v))$ defined on surface $\cal S$, we denote
\[
f_u = \frac{\partial f({\bm x}(u,v))}{\partial u} \quad {\rm and} \quad f_v = \frac{\partial f({\bm x}(u,v))}{\partial
v}.
\]
For simplicity, we set $f_1 := f_u$, $f_2 := f_v$, $f_{11} := f_{uu}$, $f_{12} := f_{uv}$ and $f_{22} := f_{vv}$.

\paragraph{Tangential Gradient Operator} Let $f \in C^1({\cal S})$,  then the {\it tangential gradient
operator} $\nabla_{\cal S}$ acting on $f$  at ${\bm x} \in {\cal S}$
is defined as
\begin{equation}
\label{tgo1} \nabla_{\cal S} f({\bm x}) =  [{\bm x}_u, {\bm
x}_v][g^{\alpha \beta}][f_u,f_v]^{\rm T}  \in \mathbb{R}^3.
\end{equation}
By the definition of $g^{\alpha \beta}$ we have
\begin{equation}
\label{tgo2} \nabla f({\bm x}) = [f_u,f_v]^{\rm T} = [{\bm x}_u, {\bm
x}_v]^{\rm T} \nabla_{\cal S} f.
\end{equation}
Obviously, $\nabla_{\cal S} f \in T_{\bm x} {\cal S}$ and hence
$\langle \nabla_{\cal S} f({\bm x}), {\bm n}\rangle = 0$.

\paragraph{Divergence Operator} Let ${\bm v}$ be a $C^1$ smooth vector field on surface $\cal
S$. Then the {\it divergence operator} ${\rm div}_{\cal S}$ acting on ${\bm v}$ is
defined as
\[
{\rm div}_{\cal S} ({\bm v}) =
\frac{1}{\sqrt{g}}\left[\frac{\partial}{\partial u},
\frac{\partial}{\partial v} \right] \left[\sqrt{g} [g^{\alpha
\beta}][{\bm x}_u, {\bm x}_v]^{\rm T}{\bm v} \right].
\]

\paragraph{Laplace-Beltrami Operator} Let $f\in C^2({\cal S})$, then $\nabla_{\cal S} f$ is a smooth
vector field on surface $\cal S$. The {\it
Laplace-Beltrami operator} (LBO) $\Delta_s$ acting on $f$ is defined
as
\begin{eqnarray*}
\Delta_{\cal S} f = {\rm div}_{\cal S} (\nabla_{\cal S} f) .
\end{eqnarray*}
With the definitions of $\Delta_{\cal S} $ and ${\rm div}_{\cal S}$,
we derive
\[
\Delta_{\cal S} f =
\frac{1}{g}(g_{22}f_{11}+g_{11}f_{22}-2g_{12}f_{12}).
\]

\subsection{Global Properties of Surface}

Here we introduce the concept of the global surface, and derive some global properties,
for example {\it Green formula}. There will be quite useful for the numerical
analysis of surface PDEs.

\begin{definition}
\label{glob1}
The ${\cal S} = \bigcup\limits_{\alpha \in \Lambda} {\cal
S}_{\alpha}$ is called a global orientable surface in $\mathbb{R}^3$ if it satisfies:
\begin{itemize}
\item[(1)] Each ${\cal S}_{\alpha} = \{{\bm x}_{\alpha}(u, v)\in \mathbb{R}^3: [u, v]^{\rm T} \in {\Omega}_{\alpha}
\subset \mathbb{R}^2\}$ is a regular surface in $\mathbb{R}^3$.
\item[(2)] For any $\alpha, \beta \in \Lambda$, if ${\cal S}_{\alpha} \cap {\cal S}_{\beta} \neq
\varnothing$, ${\cal S}_{\alpha} \cap {\cal S}_{\beta}$ is still a
piece of surface and the composite mapping
\[
{\bm x}^{-1}_{\beta}\circ {\bm x}_{\alpha}: \Omega_{\alpha \beta} =
{\bm x}^{-1}_{\alpha}({\Omega}_{\alpha} \cap
{\Omega}_{\beta})\subset {\Omega}_{\alpha} \longrightarrow
\Omega_{\beta \alpha} = {\bm x}^{-1}_{\beta}({\Omega}_{\alpha} \cap
{\Omega}_{\beta})\subset {\Omega}_{\beta}
\]
is differentiable.
\item[(3)] For every ${\cal S}_{\alpha}$, there is an orientation ${\bm n}_{\alpha}$.
Moreover ${\bm n}_{\alpha} = {\bm n}_{\beta}$
for ${\bm x} \in {\cal S}_{\alpha} \cap {\cal S}_{\beta}$.
\end{itemize}
\end{definition}

\begin{theorem}[Green formula for LBO]
Let ${\cal S}$ be an orientable surface, \textcolor{black} {and}
$\Omega$ be a subregion of ${\cal S}$ with a piecewise smooth
boundary $\partial{\Omega}$. Let ${\bm n}_c$ be the outward unit normal along the boundary $\partial{\Omega}$. Then for
a given smooth vector field ${\bm v} \in C^1(\Omega)$, we have
\begin{equation}
\label{greenf} \int_{\Omega}\langle {\bm v},\nabla_s f\rangle + f\,
{\rm div}_{\cal S}({\bm v}){\rm d}A = \int_{\partial \Omega} f \langle {\bm
v},{\bm n}_c \rangle {\rm d} s.
\end{equation}
\end{theorem}

\subsection{Sobolev Spaces on Surface}

The FEM is applied to solve PDEs on surface, therefore the theory of Sobolev spaces on surface is requisite.
Assume that ${\cal S}$ is a sufficiently smooth surface. For a given constant $k$ and a function $f \in C^{\infty} ({\cal S})$,
denote $\nabla^{k} f$ the $k$-th order covariant derivative of function $f$, with the convention $\nabla^{0}f = f$. Let
\[
C_{k}({\cal S}) = \left\{ f \in C^{\infty}({\cal S}): \int_{\cal S} |\nabla^{j} f|^{2} {\rm d}A \leq \infty~\text{for}~j = 0,\ldots, k \right\}.
\]
We have the following definition of Sobolev space $H^{k}({\cal S})$.
\begin{definition}
\label{sursobo}
Let ${\cal S}$ be a compact surface with at least $k$-th order smoothness.
Sobolev space $H^{k}({\cal S})$ is the completion of $C_{k}({\cal S})$ in the sense of norm
\begin{equation}
\label{ssobo}
\|f\|_{H^{k}({\cal S})} := \left( \sum\limits^{k}_{j=0} \int_{\cal S} |\nabla^{j} f|^{2} {\rm d}A \right)^{\frac{1}{2}}.
\end{equation}
\end{definition}
For the compact surface $\cal S$, we have
\[
C_{k}({\cal S}) = C^{\infty}({\cal S})\subset C^{k}({\cal S}) \subset H^{k}({\cal S}).
\]

\setcounter{theorem}{0} \setcounter{equation}{0}
\setcounter{figure}{0}
\section{Isogeometric Analysis Based on Extended Loop Subdivision}
\label{loopintro}

Isogeometric analysis is a recently developed computational approach that the
solution space for dependent variables is represented in terms of the same functions
which represent the geometry. The isogeometric procedures can be developed based on
NURBS, A-patches or subdivision surfaces \cite{Hughes2010}. In this work, we investigate
isogeometric analysis for surface PDEs based on subdivision surfaces.

A subdivision surface is a method of representing a smooth surface, which can be generated
through a progressively iterative subdivision process
starting from an initial control mesh by a prescribed subdivision scheme. The subdivision schemes used to
generate smooth surfaces can be separated into two categories: interpolatory and
approximatory. In the first category, the vertex positions of the initial mesh
are fixed, and only the positions of new added vertices need to be computed
in each subdivision step (see \cite{Kobbelt97, Zorin96}). However, for the second category,
both the old and new vertex positions are required to update during each refinement (see \cite{Catmull78, Loop1}).
In general, the surfaces generated by approximating schemes have better quality than those
generated by interpolating ones. In this article, the extended Loop subdivision is considered,
which falls into the second category.

\subsection{Extended Loop Subdivision Scheme}
\label{subsec_scheme}

Let $\Omega^0_h$ be the input closed triangular mesh, which denotes as the initial control mesh of the conventional Loop subdivision.
In each refinement step, each triangle is subdivided into four sub-triangles in which all the new
vertices are generated by the weighted average of the old vertices, and all the old vertex positions on the refined mesh
$\Omega^{k+1}_h$ are computed by the weighted average of the vertex positions on the mesh $\Omega^k_h$. Let ${\bm p}^{k}_0$ be a vertex with one-ring
neighbors ${\bm p}_j^k$ ($j=1, \dots, n$) on $\Omega^k_h$, where $n$ is the valence of ${\bm p}^{k}_0$. The positions of the new
generated vertices ${\bm p}^{k+1}_i$ on $\Omega^{k+1}_h$ corresponding to the edges of the previous mesh are computed by
\begin{equation}
\label{edge}
{\bm p}^{k+1}_i = \frac{1}{8}(3{\bm p}^{k}_0 + 3{\bm p}^{k}_i + {\bm p}^{k}_{i-1} + {\bm p}^{k}_{i+1}), \quad i = 1, \ldots, n,
\end{equation}
where the subscript index $i$ is regarded as modulo $n$. The old vertex gets new position according to the rule
\begin{equation}
\label{ver}
{\bm p}^{k+1}_0 = (1-n\alpha) {\bm p}^{k}_0 + \alpha \sum^{n}_{j=1} {\bm p}_j^k,
\end{equation}
where
\begin{equation}\label{alphavalue}
\alpha = \frac{1}{n}\left[\frac{5}{8} - \left(\frac{3}{8} + \frac{1}{4}\cos\frac{2\pi}{n}\right)^2\right]
\end{equation}
proposed by Loop in \cite{Loop1} is named as \textit{Loop's weight} \cite{Xu12}.

The conventional Loop subdivision scheme is suitable only for subdividing closed triangular control meshes
without boundaries. Actually, for a lot of geometric modeling problems, the surfaces are usually constructed
in a piecewise manner with fixed boundaries. In such a case, Loop subdivision scheme cannot be implemented
near the boundaries of the control mesh. To overcome this deficiency, the extended Loop subdivision scheme
was proposed by Biermann et al. in \cite{Biermann2000}, which can treat triangular control meshes with
boundaries. Note that the subdivision rule for boundaries is the same as that of cubic B-spline, and the control
vertices on boundaries are treated as the control vertices of cubic B-spline curves with spaced knots.
The extended Loop subdivision scheme is described as follows.

\subsubsection{Vertex Refinement}

The vertices can be separated into three types: corner vertex, boundary vertex and interior vertex. For different type of vertex,
the refined strategy is given as:
 \begin{itemize}
\item[(1)] Corner vertex: The corner vertices are to be interpolated, meaning which are fixed.
\item[(2)] Boundary vertex: Let ${\bm x}_i$ be a boundary vertex, and let ${\bm
x}_l$ and ${\bm x}_r$ be its two neighbor vertices on the boundary,
then ${\bm x}_i$ is updated by $\frac 18 {\bm x}_l+ \frac 34 {\bm
x}_i + \frac 18 {\bm x}_r$.
\item[(3)] Interior vertex: Use the conventional Loop scheme as (\ref{edge}) and (\ref{ver}).
 \end{itemize}

\subsubsection{Edge Refinement}

Overall, we divide the edges into boundary, sub-boundary and interior edges. Boundary edges lie on the boundaries, which are the features
of control mesh in general. Sub-boundary edges are not boundary edges but adjacent to
the boundary vertices. The rest are the interior edges. For different type of edge, the edge refinement is performed as:
 \begin{itemize}
\item[(1)] Boundary edge: Let $[{\bm x}_l {\bm x}_r]$ be a boundary edge.
The newly added vertex on this edge is the average of ${\bm x}_l$
and ${\bm x}_r$, i.e., $\frac{1}{2} ({\bm x}_l + {\bm x}_r$).
\item[(2)] Sub-boundary edge: Let $[{\bm x}_i {\bm x}_j]$ be a sub-boundary
edge with ${\bm x}_i$ being a boundary vertex, ${\bm x}_l$ and ${\bm
x}_r$ being the two wing neighbor vertices of $[{\bm x}_i {\bm
x}_j]$, then the newly added vertex on this edge is defined as
$(\frac 34 - \gamma) {\bm x}_i + \gamma {\bm x}_j + \frac{1}{8}{\bm
x}_l + \frac{1}{8}{\bm x}_r$ where $\gamma = \frac 12 - \frac 14
\cos \theta_k$ with $\theta_k = \frac{\pi}{k}$ for a boundary point,
$\theta_k =\frac{\alpha}{k}$ for  a convex corner point, and
$\theta_k = \ \frac{2\pi- \alpha}{k}$ for a concave corner point.
\item[(3)] Interior edges: Use the conventional Loop scheme as (\ref{edge}) and (\ref{ver}).
\end{itemize}

\subsection{Limit Form of the Extended Loop Subdivision}
\label{subsec_limitform}

The explicit representation for the limit position of each vertex can be derived, which is stated as the following lemma.

\begin{lemma}[\cite{Loop1}]
\label{limitposi}
Let ${\bm p}^0_0$ be a control vertex of valence $n$ on the mesh $\Omega^0_h$, and ${\bm p}^0_i$, $i = 1, \ldots, n$,
be its one-ring neighbor control vertices. The vertices sequence $\{{\bm p}^k_0\}$ converges to a unique point
\begin{equation}
\label{limit} {\bm p}^{\infty}_{0} = (1-nl){\bm p}^0_0 + l \sum\limits^{n}_{j=1}
{\bm p}^0_j, \quad l = \frac{1}{n+ 3/(8\alpha)},
\end{equation}
as the subdivision step $k\rightarrow \infty$, where the $\alpha$ is defined by (\ref{alphavalue}).
\end{lemma}

\begin{figure}[ht!]
\centerline {
\begin{tabular}{cc}
\includegraphics[scale=0.18]{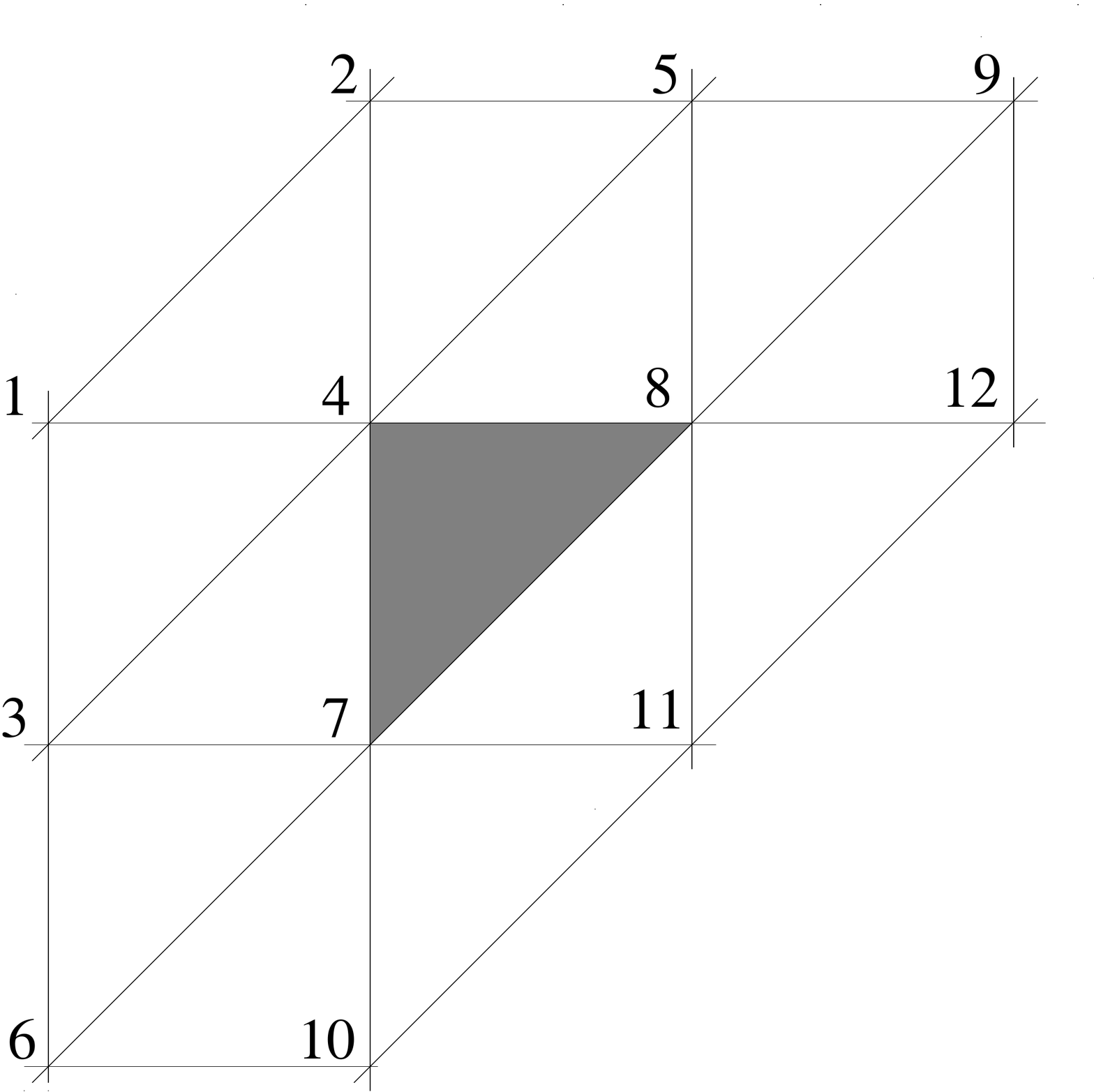} &
\includegraphics[scale=0.28]{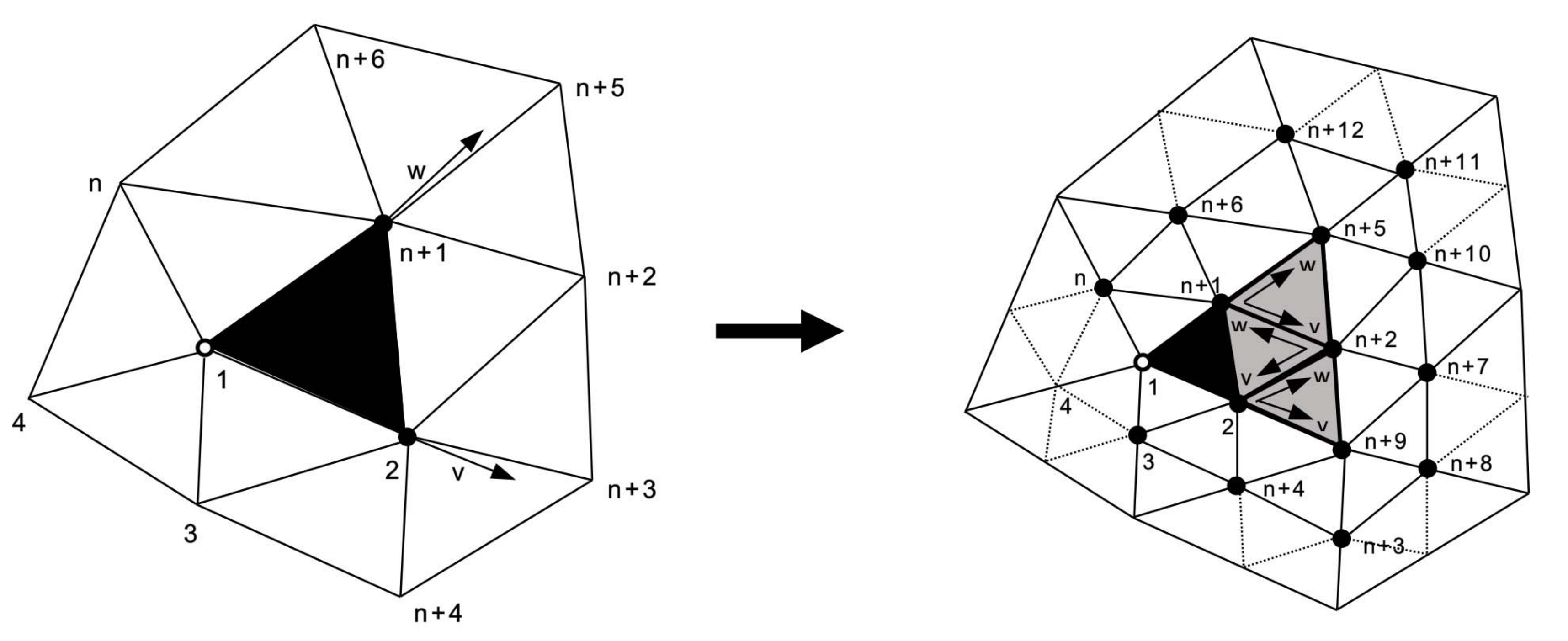} \\
($a$) &  ($b$)
\end{tabular}
}  \caption{{\small ($a$) A single regular triangular patch defined by 12 control vertices. ($b$) Subdividing an irregular patch once generates
three sub-patches and enough control vertices for evaluating three
shaded sub-patches. Here the vertex labelled `1' is extraordinary of
valence 5. }} \label{subsub}
\end{figure}

The limit surface of the extended Loop subdivision is $C^2$ everywhere except at the irregular
vertices (i.e., having valence other than six) where is $C^1$ (see \cite{Schweitzer96}).
By Lemma (\ref{limitposi}), we can evaluate the position of the limit surface at any finite subdivision level and
at any vertex by weighted averaging the vertex and its neighbors. However, it is not feasible
to compute any point on the limit surface. Fortunately, there exists the explicit expression for the
regular Loop subdivision surface patch. If all vertices of the triangle have valence six
and none of its two-ring neighbor vertices is a boundary vertex, the resulting surface patch
is called {\it regular}. The regular patch can be exactly described by
a quartic box-spline, more precisely, which is formulated by the linear combination of 12 basis functions:
\begin{equation}
\label{loopevalua}
{\bm p}(\xi,\eta) =  \sum_{i=1}^{12} B_i(\xi,\eta) {\bm p}^k_i,
\end{equation}
where $(1-\xi-\eta, \xi, \eta)$ are the barycentric coordinates of the shaded unit reference triangle and ${\bm p}^k_i$ are the
indexed vertices of the control mesh in Figure \ref{subsub} ($a$). The analytic expressions of the basis functions $B_i$ are given as those in \cite{Xu12}.
Each triangle of the control mesh can be regarded as a parametric domain, which corresponds to one
triangular patch of the limit surface.

If a triangle is {\it irregular}, i.e., at least one of
its vertices has a valence other than six or one of its two-ring neighbor
vertices is a boundary vertex, the resulting surface patch cannot be
represented by a quartic box-spline. For evaluating irregular patches, we use the
fast scheme proposed by Stam in \cite{Stam983}, in which the mesh needs to be
subdivided repeatedly until the parameter values of interest are interior
to a regular patch (See Figure \ref{subsub} ($b$)). Hence,
the basis functions given in (\ref{loopevalua}) can be used to represent the geometry exactly.

\subsection{Finite Element Space for IGA-Loop}
\label{fefs}

Let ${\cal S}$ be a considered smooth surface, which is generated by the limit of the extended Loop subdivision scheme.
And we denote $\Omega_h^k$ the $k$-th refined control mesh of the discretized representation ${\cal S}_h$, where ${\bm p}^k_{i}$ denotes
its $i$-th control vertex and $h$ is the length of maximal edge.
Each control vertex ${\bm p}^k_{i}$ on $\Omega_h^k$ can be evaluated recursively, which is described
in Section \ref{subsec_scheme}. Each triangle patch of $\Omega_h^k$
can be locally represented by a linear combination of explicit box-spline tiles as (\ref{loopevalua}).
The boundaries of ${\cal S}_h$ are represented as the cubic
B-spline curves which are preserved as the subdivision proceeds. The subdivision step is equivalent to
the $h$-refinement of the NURBS-based isogeometric analysis.

Let us define the basis functions of the finite element
space for IGA-Loop. For each control vertex ${\bm p}^k_i$ on the control mesh $\Omega^k_{h}$ of the surface ${\cal S}_h$,
including boundary control vertices, we associate it with a basis function $\phi_i$, where $\phi_i$ is defined as the limit
of the extended Loop subdivision scheme applying to the zero control values everywhere except at ${\bm p}^k_i$ where it is one.
Note that the basis functions have a compact support within two-ring neighbors. The mesh $\Omega^k_{h}$, formulated as piecewise
patches, is served as the parametric domain of the basis functions. Actually the set of quartic box-splines corresponding to
each subdivided computable control patch are utilized to represent the geometry of interest and the solution space for dependent variables.

The control mesh $\Omega^k_{h}$, as a  piecewise linear surface, is
served as the definition domain of the basis function $\phi_i$. The
mapping from $\Omega^k_{h}$ to $\phi_i$ is defined by a dual
subdivision process. More precisely, when the extended Loop
subdivision scheme is applied to the control function values
recursively, the linear subdivision scheme (each triangle is
partitioned into four equal-sized sub-triangles) is applied to the
control mesh correspondingly. The limit of the former is $\phi_i$
and that of later is $\Omega^k_{h}$ itself.

The basis functions share some properties with the well-known
B-spline basis. These properties are important in the proposed
method. Now let us describe them as follows.

\begin{itemize}

\item[(1)]
{\it Positivity}. The weights of the extended subdivision rules are
positive. Hence the basis function $\phi_i$ is nonnegative
everywhere and positive around ${\bm p}^k_i$.

\item[(2)]
{\it Locality}. It is known that the limit value at a control vertex
is a linear combination of the one-ring neighbor values. Hence, the
limit value is zero at a control vertex if the control values on the
one-ring neighbor control vertices are zero. Therefore, the support
of the basis function is within the two-ring neighborhood.
\item[(3)]
{\it Partition of Unity}. Since all the subdivision rules have the
properties that the weights are summed to one.
Therefore, if we choose all the control values as one. The control values after one
subdivision step are still one. This implies that $\sum \phi_i = 1$. This property is called partition of unity.

\item[(4)]
{\it Interpolatory Properties at the Boundary}. The extended
subdivision rules on the boundary do not involve the interior
control vertices. Hence the basis functions for the interior control
vertices are zero at the boundary.

\item[(5)] {\it Tangential Property}. Let $ {\bm p}^k_i$ be a control vertex, with
non of its one-ring neighbor control vertices \textcolor{black} {are}
boundary control vertices. Then $\nabla_s \phi_i$ vanishes on the
boundary. This fact can be observed by considering the
eigen-decomposition of the control vertices. Let ${\bm x}^{k} \in
\mathbb{R}^{(n+1)\times 3}$ be a vector consisting of one-ring
neighbor control vertices of ${\bm p}_i^{k}$ at the subdivision
level $k$, $S\in \mathbb{R}^{(n+1)\times(n+1)}$ be the local
subdivision matrix that converts ${\bm x}^{k}$ to ${\bm x}^{k+1}$,
i.e.,
\[
{\bm x}^{k+1} = S {\bm x}^{k} = S^k {\bm x}^{1}, \ \ \ k = 1,
2, \ldots.
\]
Here $n$ stands for the valence of ${\bm p}_i^{k}$. Suppose ${\bm
x}^{1}$ is decomposed into
\[
{\bm x}^{1} = {\bm e}_0{\bm a}^{\rm T}_0 + {\bm e}_1{\bm a}^{\rm T}_1 + {\bm
e}_2{\bm a}^{\rm T}_2 + \cdots + {\bm e}_n{\bm a}^{\rm T}_n, \ \ {\bm a}_j \in
\mathbb{R}^3,
\]
where ${\bm e}_0$, ${\bm e}_1$, $\ldots$, ${\bm e}_n$ are the
eigenvectors of $S$. Here we assume that these eigenvectors are
arranged in the order of non-increasing eigenvalues $\lambda_j$.
Then
\[
{\bm x}^{k+1} = \lambda_0^k {\bm e}_0{\bm a}^{\rm T}_0 + \lambda_1^k
{\bm e}_1{\bm a}^{\rm T}_1 + \lambda_2^k {\bm e}_2{\bm a}^{\rm T}_2 + \cdots +
\lambda_n^n{\bm e}_n{\bm a}^{\rm T}_n,
\]
where $\lambda_0 = 1$, $\lambda_1 = \lambda_2 < 1$. It is well-known
that the limit position at the center is ${\bm a}_0$. The tangent
direction at this vertex are ${\bm a}_1$ and ${\bm a}_2$, and ${\bm
a}_j$ is given by ${\bm a}_j^{\rm T} = \tilde{\bm e}_j^{\rm T}{\bm x}^{1}$.
$\tilde{\bm e}_j$ are the left eigenvectors of $S$ with
normalized condition $\tilde{\bm e}_j^{\rm T} {\bm e}_j = 1$.
The analysis above is valid for control function values.
 The fact that $\nabla_{\cal S} \phi_i$ vanishes on the
boundary implies that the tangent vector of the subdivision surface
on the boundary determined by the boundary and sub-boundary control
vertices. This is similar to B\'ezier and B-spline surfaces.

\item[(6)]
{\it Linear independency}. The functions $\phi_i$, $i = 0, \ldots, m$, are defined
by the extended Loop subdivision scheme.  This fact can be derived from the
unique solvability of the following interpolation problem:

\begin{lemma}[\cite{Xuerror05}]
\label{inter00}  Let ${\bm p}^{\infty}_{i}$ be the limit position of the
control vertex ${\bm p}^{k}_{i}$ with its 1-ring neighbor vertices ${\bm
p}^k_j, j = 1,2,\ldots, n_i$, $u({\bm p}^{\infty}_i)$ be the $i$-th
interpolation function value, and ${\rm u}({\bm p}^k_i)$ be the $i$-th control
function value, where $i = 1,\ldots, m$. The system
\[
(1 - n_i l_i){\rm u}({\bm p}^k_i) + l_i \sum\limits^{n_i}_{j=1} {\rm u}({\bm
p}^k_{j}) = u({\bm p}^{\infty}_i),\ \ i =1,\ldots, m,
\]
is always solvable uniquely. Here the $l_{i}$ is defined as (\ref{limit}).
\end{lemma}

\item[(7)]
{\it Interpolant Error}. With the basis of Lemma \ref{inter00}, we derive the following
interpolant error estimation:
\begin{lemma}[\cite{PanXu15}]
\label{interpolant} Assume that $\Omega_h^k$ is the $k$-th refined control mesh for the smooth surface ${\cal S}$  generated
by the extended Loop subdivision, where $h$ stands for the length of the maximal edge. For $u \in H^2({\cal S})$,
there exists an interpolation function $u_{I} \in H^{2}({\cal S}_h)$ such that
\begin{equation}
\label{interpo} \| u-u_{I}\|_{H^{l}({\cal S})} \leq ch^{2-l} \|u\|_{H^{2}({\cal S})},\quad l = 0, \ 1.
\end{equation}
\end{lemma}

\end{itemize}

\setcounter{theorem}{0} \setcounter{equation}{0}
\setcounter{figure}{0}
\section{Applications to Surface PDEs}
\label{SurfacePDE}

We consider the approximations for the Laplace-Beltrami harmonic, biharmonic and
triharmonic problems on different surfaces by the proposed IGA-Loop strategy. The numerical examples we provide in Section \ref{numer_ex} include both open and closed
surfaces, more concretely, which are the Laplace-Beltrami harmonic problem on a quarter of cylinder and one-eighth of sphere, the Laplace-Beltrami biharmonic problem on a cylinder, and the Laplace-Beltrami triharmonic problem on a unit sphere. The numerical error will be analyzed which shows the consistent convergence rate with (\ref{surffem}). For comparison, we also consider the approximations of these three problems by standard linear element discretization. Hence we adopt the mixed formulations of the high-order Laplace-Beltrami biharmonic and triharmonic problems.

\subsection{Surface PDEs}

\paragraph{Laplace-Beltrami Harmonic Problem} The Laplace-Beltrami harmonic equation with zero boundary
condition can be written as
\begin{equation}
\label{secondorder} \left\{
\begin{array}{l}
- \Delta_{\cal S} u  = f, \quad \text{in}~{\cal S},
\\[0.5em]
u= 0, \quad \text{on}~\partial {\cal S},
\end{array}
\right.
\end{equation}
where the boundary $\partial {\cal S}$ of surface $\cal S$ is Lipschitz continuous.
Here $f: {\cal S}  \rightarrow \mathbb{R}$ is a given
sufficiently regular function, and $\Delta_{\cal S}$ is the Laplace-Beltrami operator.
The weak form of (\ref{secondorder}) is given as
follows:
\begin{equation}
\label{weakcont1}
\left\{
\begin{array}{l}
{\rm Find}\  u \in H^1_{0}({\cal S}) \ {\rm such\  that} \\
[0.5em]
\displaystyle \int_{\cal S} \nabla_{\cal S} u \cdot \nabla_{\cal S}
\varphi \ {\rm d} A = \int_{\cal S} f \varphi \ {\rm d}A,  \quad \forall \varphi \in
H^1_0({\cal S}).
\end{array}
\right.
\end{equation}

\paragraph{Laplace-Beltrami Biharmonic Problem} The Laplace-Beltrami biharmonic equation with homogeneous boundary conditions reads as
\begin{equation}
\label{fourthorder} \left\{
\begin{array}{l}
\Delta_{\cal S}^2 u  = f, \quad \text{in}~{\cal S},
\\[0.5em]
u= \frac{\partial u}{\partial {\bm n}} = 0, \quad \text{on}~\partial {\cal S},
\end{array}
\right.
\end{equation}
where ${\bm n}$ is the outward directed unit vector normal to the boundary $\partial {\cal S}$. Let $\Delta_{\cal S} u = -v$, then the mixed weak formulation of (\ref{fourthorder}) reads as
\begin{equation}
\label{weakcont2}
\left\{
\begin{array}{l}
{\rm Find}\  (u,v) \in H^1_{0}({\cal S}) \times  H^{1} ({\cal S})\ {\rm such\  that}
\\ [0.8em]
\displaystyle \int_{\cal S} \nabla_{\cal S} u \cdot \nabla_{\cal S} \varphi \ {\rm d} A  - \int_{\cal S} v \varphi \   {\rm d} A  =0,   \quad \forall \varphi \in H^1_{0}({\cal S}),
\\ [0.8em]
\displaystyle \int_{\cal S}  \nabla_{\cal S} v \cdot \nabla_{\cal S} \psi \ {\rm d} A  = \int_{\cal S} f \psi \ {\rm d}A,  \quad \quad  \quad \forall \psi \in
H^1({\cal S}),
\end{array}
\right.
\end{equation}
where $v: {\cal S} \rightarrow  \mathbb{R}$  are the auxiliary unknowns.

\paragraph{Laplace-Beltrami Triharmonic Problem} The Laplace-Beltrami triharmonic equation with homogeneous boundary conditions is formulated as
\begin{equation}
\label{sixthorder} \left\{
\begin{array}{l}
- \Delta_{\cal S}^3 u = f, \quad \text{in}~{\cal S},
\\[0.5em]
u = \frac{\partial u}{\partial {\bm n}} = \Delta_{\cal S} u = 0, \quad \text{on}~\partial {\cal S}.
\end{array}
\right.
\end{equation}
Let $\Delta_{\cal S} u = -v$ and $\Delta_{\cal S} v = -w$, then the mixed formulation of (\ref{sixthorder}) reads as
\begin{equation}
\label{weakcont3}
\left\{
\begin{array}{l}
{\rm Find}\  (u,v,w) \in H^1_{0}({\cal S}) \times  H^{1}_0 ({\cal S}) \times  H^{1} ({\cal S}) \ {\rm such\  that}
\\ [0.8em]
\displaystyle \int_{\cal S} \nabla_{\cal S} u \cdot \nabla_{\cal S} \varphi \ {\rm d} A  - \int_{\cal S} v \varphi \   {\rm d} A  =0,   \quad \forall \varphi \in H^1_{0}({\cal S}),
\\ [0.8em]
\displaystyle \int_{\cal S} \nabla_{\cal S} v \cdot \nabla_{\cal S} \psi \ {\rm d} A  - \int_{\cal S} w \psi \   {\rm d} A  =0,   \quad \forall \psi \in H^1_{0}({\cal S}),
\\ [0.8em]
\displaystyle \int_{\cal S}  \nabla_{\cal S} w \cdot \nabla_{\cal S} \zeta \ {\rm d} A  = \int_{\cal S} f \zeta \ {\rm d}A,  \quad  \quad \quad \forall \zeta \in H^1({\cal S}),
\end{array}
\right.
\end{equation}
where $v, w: {\cal S} \rightarrow  \mathbb{R} $ are the auxiliary unknowns.

\subsection{IGA-Loop Discrertization of Surface PDEs}
\label{fed}

The finite element space defined by the limit form of the extended Loop subdivision is employed to
represent the geometry and also represent the solution space for the discretization of the
surface PDEs. The input triangular mesh serves as the initial control mesh of the extended Loop
subdivision. Such function is $C^2$-continuous everywhere except at the extraordinary
vertices with $C^1$-continuity. Using globally $C^1$-continuous basis functions induced by the extended Loop subdivision
yields our IGA finite dimensional space that is the subspace of the Sobolev space $H^2({\cal S})$.
For each control vertex ${\bm x}_{i}$ of the limit surface ${\cal S}$,
we associate it with a basis function $\phi_{i}$, where $\phi_{i}$ is defined
by the limit of the extended Loop subdivision for the zero control values everywhere except at ${\bm x}_{i}$ , where it is one.
Hence the support of $\phi_{i}$ is a piecewise function, and covers the two-ring neighborhood of the vertex ${\bm x}_{i}$.

The function can be locally parameterized on the unit triangle
defined by $T = \{(\xi, \eta)\in \mathbb{R}^2: \xi \ge 0, \eta \ge 0, \xi
+ \eta \le 1 \}$ where $(1-\xi-\eta, \xi, \eta)$ are the barycentric coordinates of
the unit triangle. Using this parameterization, the discretized
representation of the surface domain is $ {\cal S}_{h} = \bigcup_{i= 1}^k {\mathcal T}_{i}, \ \mathring{\mathcal T}_{i} \cap
\mathring{ \mathcal T}_{j} = \varnothing \ {\rm for} \ i \ne j$, where $\mathring{\mathcal T}_{i}$ is the interior of
triangular patch $\mathcal{T}_{i}$. Each triangular patch can
be parameterized locally as
\[
{\bm e}_{i}: \ \  T \rightarrow \mathcal{T}_{i};  \ \ (\xi, \eta) \  \mapsto \ {\bm e}_{i}(\xi, \eta).
\]
Denote ${\bm e}_{j}, j = 1,\ldots, m_{i}$ be the two-ring neighborhood patches around the control vertex ${\bm x}_{i}$. Then if ${\bm e}_{j}$ is regular, the explicit box-spline expression as in (\ref{loopevalua}) exists for $\phi_{i}$ on ${\bm e}_{j}$. If ${\bm e}_{j}$ is irregular, local subdivision, as described in Section 3.2, is needed around ${\bm e}_{j}$ until the parameter values of interest are interior to a regular patch. The parameterization has no overlap. Each control vertex ${\bm x}_{i} \in {\cal S}$ has its unique parameter coordinates.
Using the set of the basis functions $\{\phi_{i}\}$, the limit surface ${\cal S}$ of the extended Loop subdivision is expressed as
\begin{equation}
\label{patch} {\cal S} = {\bm x}(\xi,\eta) = {\bm x}( x(\xi,\eta),y(\xi,\eta),z(\xi,\eta)) = \sum\limits^{n}_{i=1} \phi _i(\xi,\eta) {\bm
x}_i.
\end{equation}
Each control triangular surface patch of the surface ${\cal S}$ is defined locally by only a few related basis functions, since the
supports of the basis are compact. The surface boundaries are represented as the cubic B-spline curves
which are preserved as the subdivision proceeds. With the parameterization, the differential operators on the surface as described in Section 2 can be computed directly, and the computation of the function integration on the surface is replaced by
\[
\int_{\cal S} u \ {\rm d} A := \sum_{i} \int_{T} u({\bm e}_{i}(\xi, \eta)) \sqrt{g}\ {\rm d} \xi {\rm d} \eta.
\]

\subsubsection{Numerical Discretization}

In what follows, we discretize the weak formulations of the Laplace-Beltrami based problems. Let ${\bm x}_j, j = 1,\ldots, n_{0}$ be the interior control vertices of the surface  $\cal S$, and ${\bm x}_j, j = n_{0}+1,\ldots, n$ be its boundary control vertices. Recalling (\ref{weakcont1}), (\ref{weakcont2}) and (\ref{weakcont3}), by means of the basis functions $\phi_{j} \in H^2({\cal S}_h)$ introduced by the limit form of the extended Loop subdivision, we have the following discrete description for the function $u^h$  to be determined
\[
u^h = \sum^{n_{0}}_{j=1} \phi_j u^{h}_{j}  + \sum^{n}_{j=n_{0}+1} \phi_j u^{h}_{j},\ \ \ \
\]
where $u^h_j \ (j = n_{0}+1, \ldots, n)$ is the known boundary conditions. The other two auxiliary functions $v^h$ and $w^h$ are represented as
\[
v^{h} = \sum^{n}_{j=1} \phi_j v^{h}_j \quad  {\rm and} \quad w^{h} = \sum^{n}_{j=1} \phi_j w^{h}_j.
\]
For simplicity, denote
\[
\begin{array}{l}
{\bm U}_{K} = [u^{h}_{1}, \ldots, u^{h}_{K}]^{\rm T} \in \mathbb{R}^{K}
\end{array}
\]
the solution vector, and the other two auxiliary solution vectors
\[
{\bm V}_{K} = [v^{h}_{1}, \ldots, v^{h}_{K}]^{\rm T} \in \mathbb{R}^{K}\quad  {\rm and} \quad
{\bm W}_{K} = [w^{h}_{1}, \ldots, w^{h}_{K}]^{\rm T} \in \mathbb{R}^{K}.
\]
The coefficient matrices and the right-hand side terms are represented respectively as
\[
\begin{array}{lll}
{\bm M}_{K \times L} = [m_{ij}]^{K,L}_{ij=1}, \quad {\bm K}_{K \times L} = [k_{ij}]^{K,L}_{ij=1} \quad   {\rm and} \quad {\bm B}_{K} = [b_i]^{K}_{i=1},
\end{array}
\]
whose elements are
\[
\begin{array}{lll}
m_{ij} =  \displaystyle -\int_{\cal S} \phi_{i} \phi_{j}  \ {\rm d} A, \quad
k_{ij} = \displaystyle  \int_{\cal S} [ (\nabla_{\cal S} \phi)^{\rm T}_{i}  \nabla_{\cal S} \phi_{j} ] \ {\rm d} A \quad
{\rm and} \quad b_{i} = \displaystyle \int_{\cal S} f \phi_{i}  \ {\rm d} A.
\end{array}
\]

Take the test functions $\varphi = \phi_{i}, i = 1,\ldots,n_{0}$, for the Laplace-Beltrami harmonic problem (\ref{weakcont1}) on a quarter of cylinder and one-eighth of sphere in Test Suite 1. Considering the zero boundary condition, we obtain the linear system as
\begin{equation}
\label{second}
{\bm K}_{n_{0} \times n_{0}} {\bm U}_{n_{0}}  = {\bm B}_{n_{0}}.
\end{equation}

Choose the two sets of the test functions $\varphi = \phi_{i}, i = 1,\ldots,n_{0}$, and $\psi = \phi_{i}, i = 1,\ldots, n$, for the Laplace-Beltrami biharmonic problem (\ref{weakcont2}) on a cylinder in Test Suite 2. Considering the homogeneous boundary condition, we achieve the following linear system
\begin{equation}
\label{forth}
\left[
\begin{array}{cc}
{\bm K}_{n_{0} \times n_{0}} &  {\bm M}_{n_{0} \times n}  \\
{\bm 0} & {\bm K}_{n \times n}  \\
\end{array}
 \right] \left[
\begin{array}{c}
{\bm U}_{n_{0}} \\  {\bm V}_{n}
\end{array}
\right] = \left[
\begin{array}{c}
{\bm 0}  \\ {\bm B}_{n}
\end{array}
\right].
\end{equation}

Taking the three sets of the test functions $\varphi = \psi = \zeta = \phi_{i}, i = 1,\ldots,n$, for the Laplace-Beltrami triharmonic problem (\ref{weakcont3}) on a unit sphere in Test Suite 3, the linear system reads as
\begin{equation}
\label{sixth}
\left[
\begin{array}{ccc}
{\bm K}_{n \times n} &  {\bm M}_{n \times n} & {\bm 0} \\
{\bm 0} & {\bm K}_{n \times n} &  {\bm M}_{n \times n} \\
{\bm 0} & {\bm 0} & {\bm K}_{n \times n}
\end{array}
 \right] \left[
\begin{array}{c}
{\bm U}_{n} \\  {\bm V}_{n} \\ {\bm W}_{n}
\end{array}
\right] =  \left[
\begin{array}{c}
{\bm 0} \\  {\bm 0} \\ {\bm B}_{n}
\end{array}
\right].
\end{equation}

\subsubsection{Precompute the Basis Functions}

The framework of our IGA-Loop is the same standard process as the classical FEM.
The related basis functions and their derivatives need to be precomputed for each control patch of ${\cal S}_{h}$
before solving the linear systems. However, those computations are not intuitive in comparison with standard
linear elements since the required two-ring neighbors around each patch have arbitrary topological structure.
And additional geometric data are reflected in the subdivision schemes around boundaries.
We classify the control mesh into interior patches, sub-boundary patches, and boundary patches.
The patches containing boundary vertices are named as boundary patches. The patches adjacent to boundary
patches are called sub-boundary patches. And all of the other patches are called interior pathes.

We use the following algorithms to treat those patches. For interior patches, the Stam's algorithm is applied
to this case (see \cite{Stam983}). The sub-boundary patches can be subdivided into four interior
sub-patches by once, then the algorithm for interior case can be used. In addition, the boundary patches can be
subdivided repeatedly till their sub-patches belong to the sub-boundary case, then the above methods can be
applied to evaluate them.

Note that the Stam's fast evaluation scheme is always suitable for interior patches with only one extraordinary
vertex. Therefore, it is necessary to first subdivide once each patch of the initial
mesh. The evaluation of basis functions over their support elements uses general Gaussian integration, which just
needs a few subdivision steps to bring Gaussian quadrature knots into a box-spline patch. The integrations for
computing the matrix elements are computed by a 12-point Gaussian quadrature rule. That is,
each triangle is subdivided into four sub-triangles and a 3-point Gaussian quadrature rule is employed
on each of the sub-triangles. The 3-point Gaussian quadrature rule has error bound $O(h^3)$, where $h$ is the maximal
edge length. We adopt the adaptive numerical method developed in our former work \cite{PanXu15}.

\subsubsection{Error Estimate}

The interpolant error (\ref{interpo}) obtained in our former work can be applied to the case of second-order surface PDEs with
zero boundary condition (\ref{secondorder}). By means of the general analysis method of the classical finite elements,
we can obtain the $H^1$-norm estimation of the exact solution $u$ and the approximate solution $u^{h}$ resulting from IGA-Loop, and the
$H^0$-norm estimation by  means of the Aubin-Nitsche duality argument. The error estimation is stated as the following theorem.

\begin{theorem}
\label{fem} Consider the Laplace-Beltrami harmonic problem defined on the limit form ${\cal S}$ of the extended Loop subdivision, endowed with the zero boundary condition. Let $u \in H^2({\cal S})$ be the exact solution of the problem,  and $u^h$ be the approximate solution obtained with IGA-Loop strategy, the following error estimates hold:
\begin{equation}
\label{surffem} \| u-u^{h}\|_ {H^{l}({\cal S})} \leq ch^{2-l} \|u\|_{H^{2}({\cal S})},\quad l = 0,\  1,
\end{equation}
where  $c$  is a constant independent of $h$.
\end{theorem}
Furthermore, we also analyze the numerical errors for the cases of the fourth-order and the sixth-order surface PDEs, even defined on closed
surface like the sphere, which show the consistent convergence rates.

\setcounter{theorem}{0} \setcounter{equation}{0}
\section{Numerical Experiments}
\label{numer_ex}

In this section, we present several numerical experiments to show the performance of
our IGA-Loop method by solving the three surface PDEs on different surfaces.
It should be indicated that the numerical computations
are conducted on the limit surface generated by the extended Loop subdivision,
therefore the integration evaluation of the Gauss-Legendre knots is performed on
the triangulation of the limit surface.

We compare the proposed method (IGA-Loop) with the linear finite element method (FEM-Linear) in
terms of accuracy, convergence and computational cost. That is because FEM-Linear is the most often used finite elements in industry, which has the
same number of control vertices and error estimate as those of IGA-Loop.
The following test suites seek to assess the performance of the proposed IGA-Loop against
different surface PDEs.

\subsection{Test Suite 1: Surface Laplace-Beltrami Harmonic Equation}

Here we consider the tests for the Laplace-Beltrami harmonic equation. The first example is
\begin{equation}
\label{cylinder}
\left\{
\begin{array}{l}
-\Delta_{\cal S} u = f, \quad \text{in}~{\cal S},
\\[0.5em]
u = 0, \quad \text{on}~{\partial\cal S},
\end{array}
\right.
\end{equation}
with the exact solution
\[
u = (1-x)(1-y) {\rm sin} (\pi z)
\]
on a quarter of the cylinder ${\cal S}_1= \{(x,y,z): x^2+y^2 =1 \ \& \ x \geq 0 \ \& \  y \geq 0 \ \& \ 0 \leq z \leq 2 \}$
with calculated $f$.

\begin{figure}[ht!]
\centerline {
\begin{tabular}{ccc}
\includegraphics[scale=0.21]{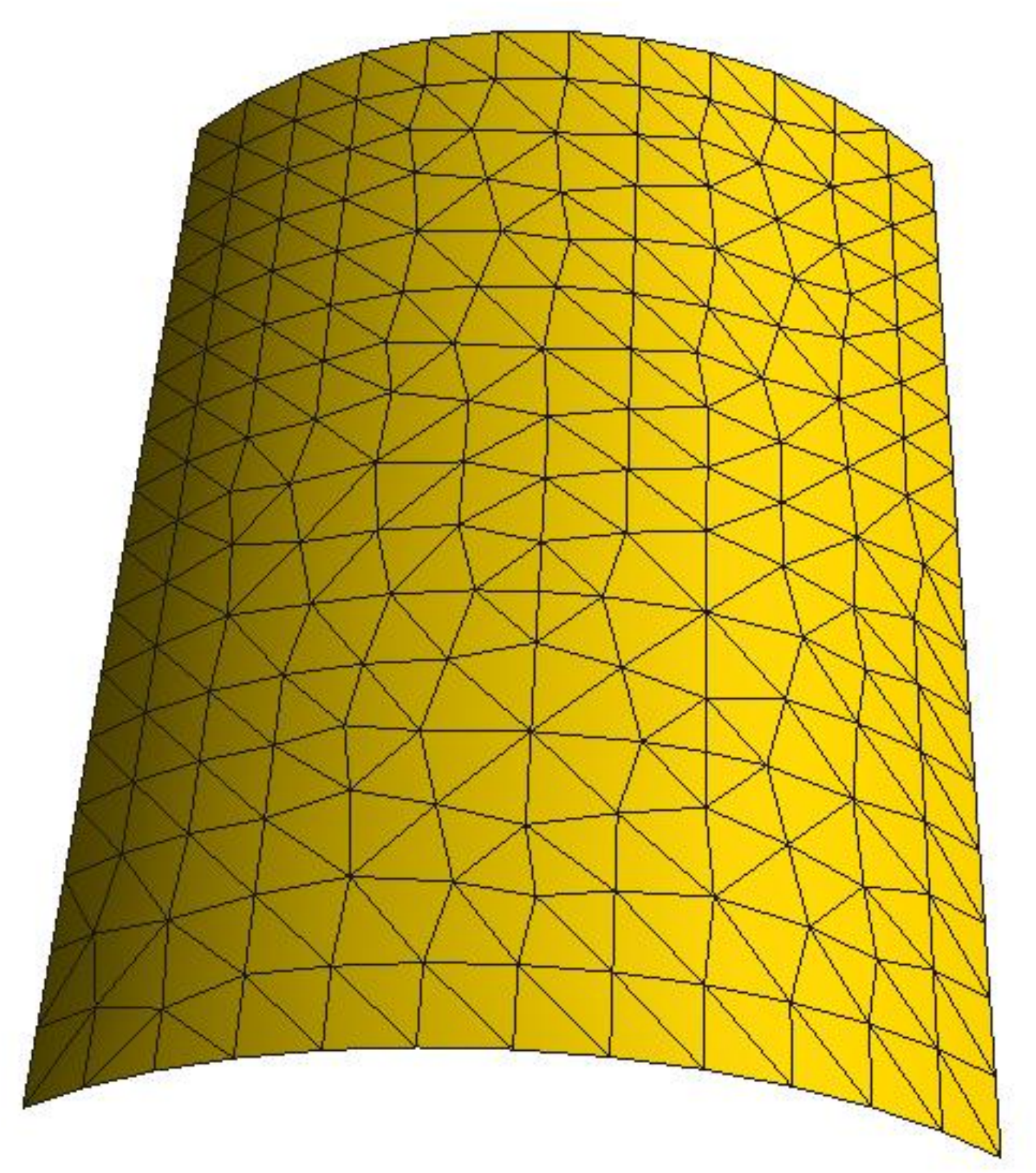}  & \ \ \ \ \ \ \
\includegraphics[scale=0.21]{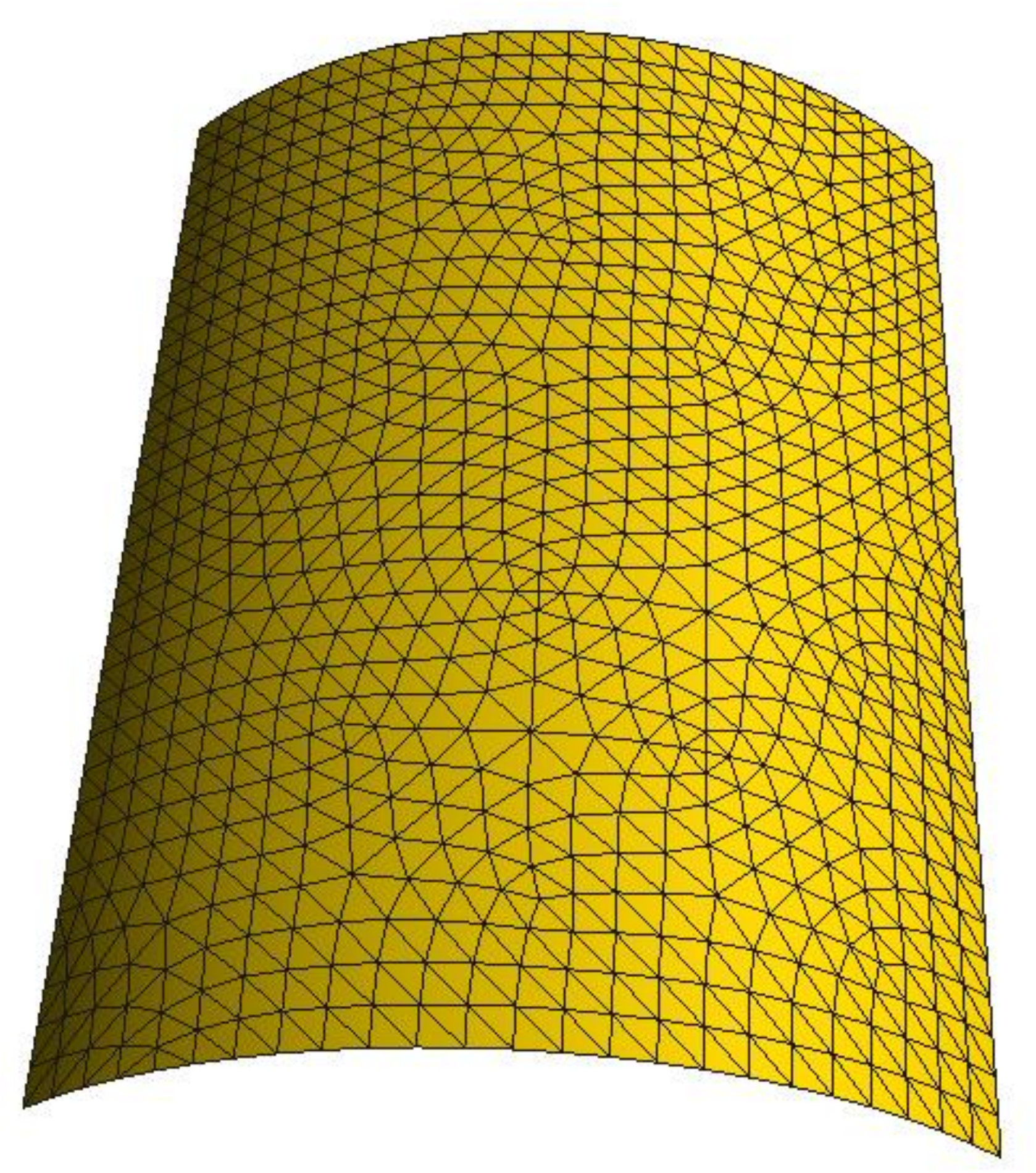}  & \ \ \ \ \ \ \
\includegraphics[scale=0.21]{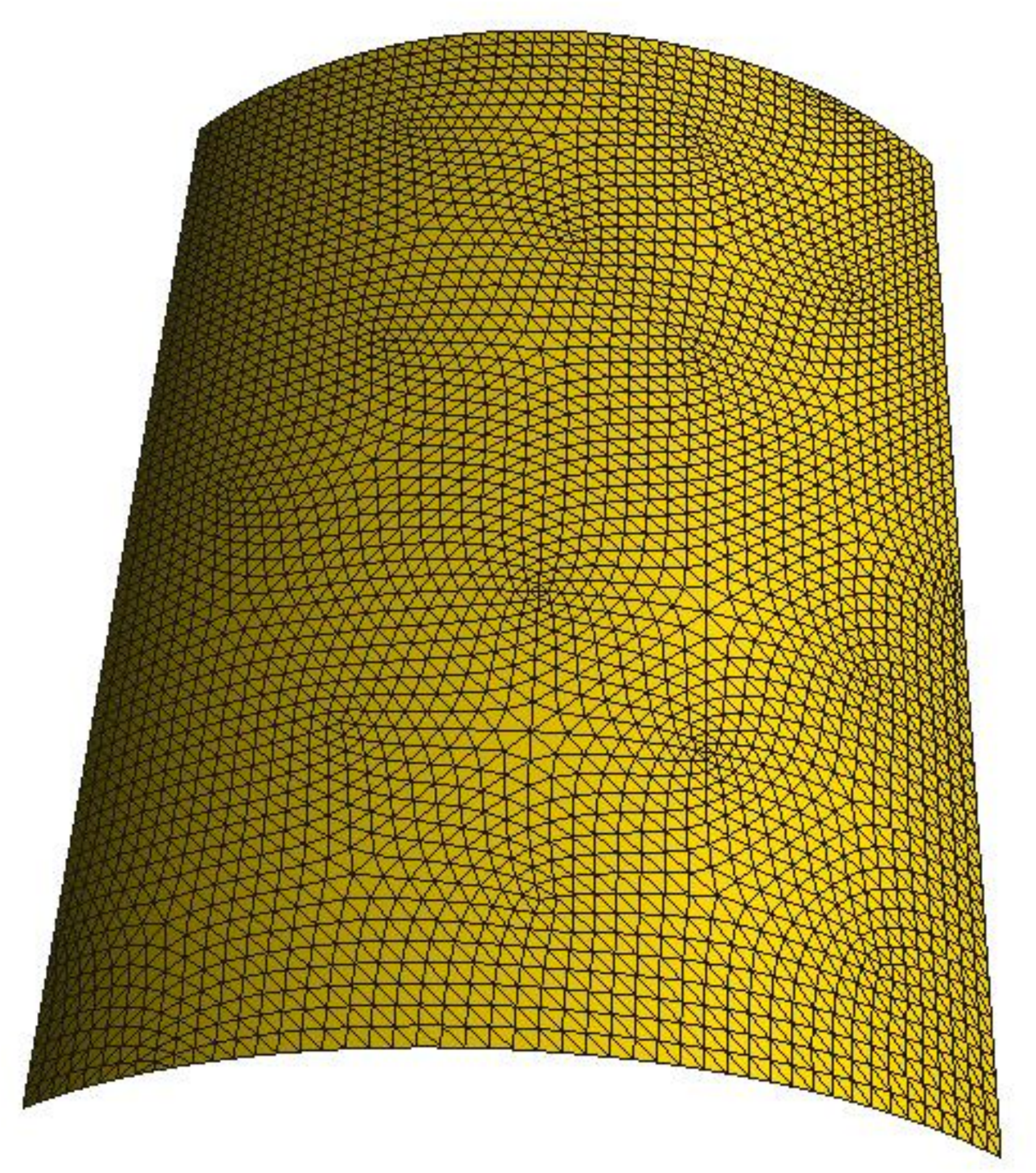}  \\
($a$) &    ($b$)  &   ($c$)  \\
\includegraphics[scale=0.22]{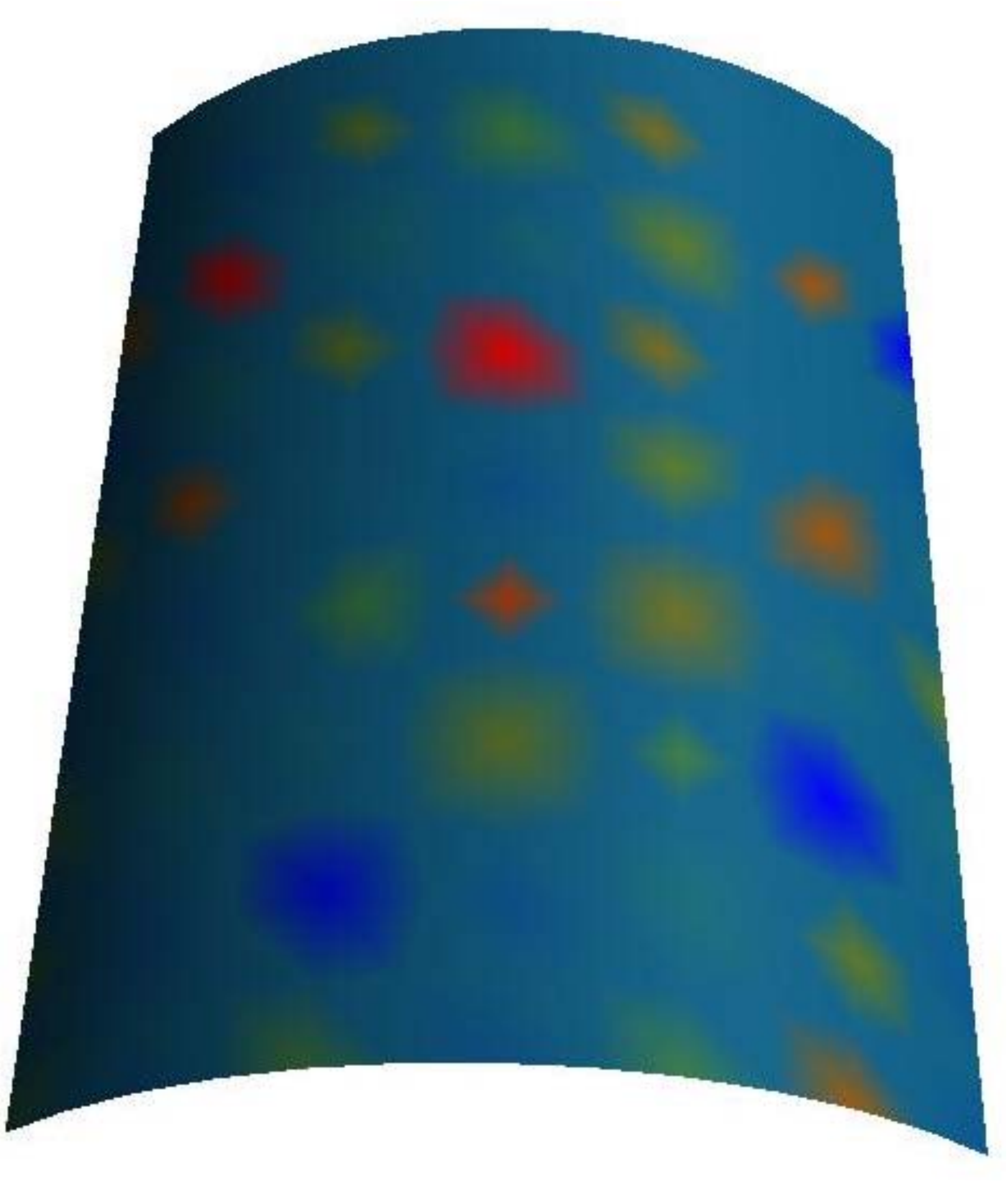}  & \ \ \ \ \ \ \
\includegraphics[scale=0.22]{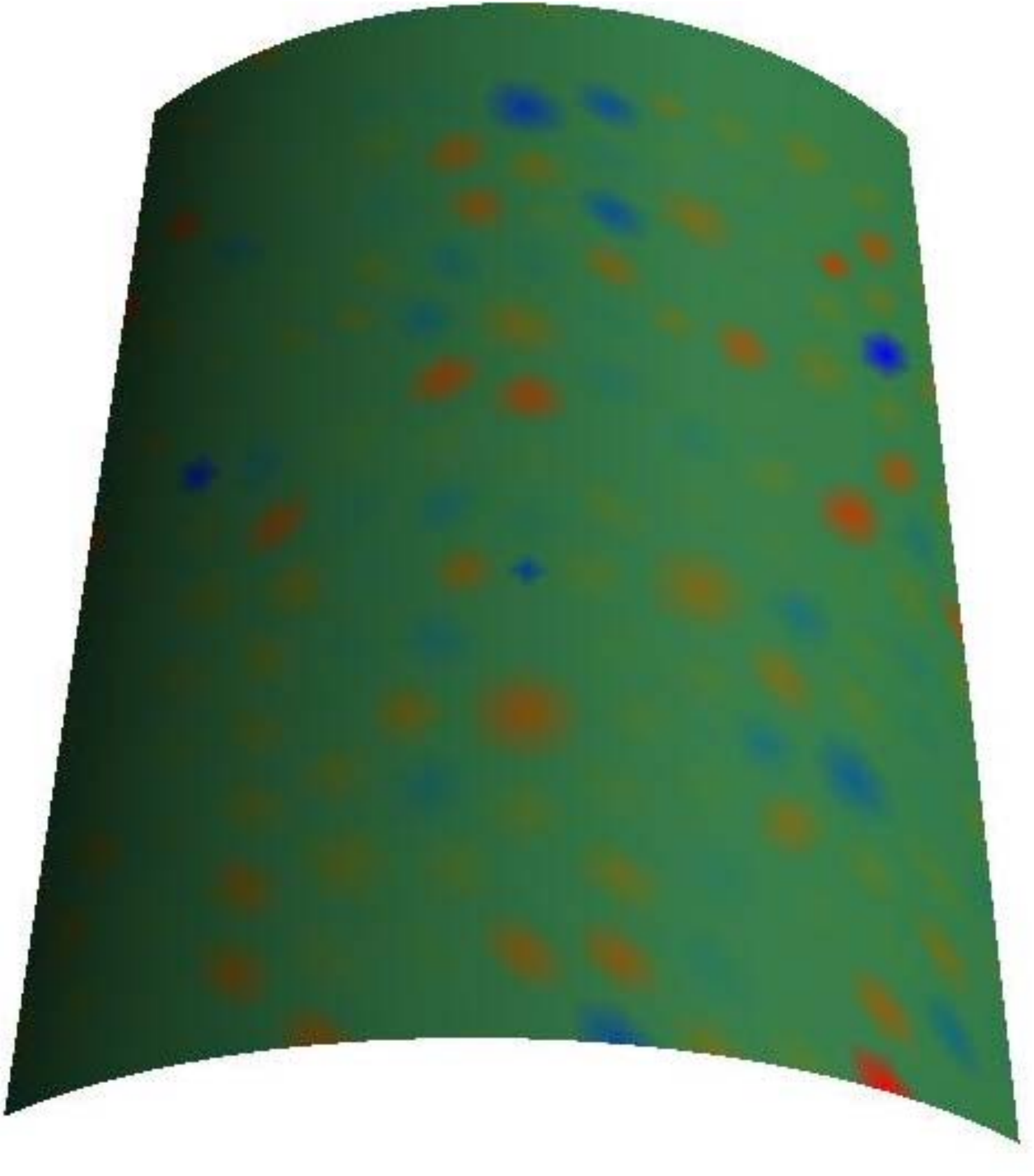}  & \ \ \ \ \ \ \
\includegraphics[scale=0.22]{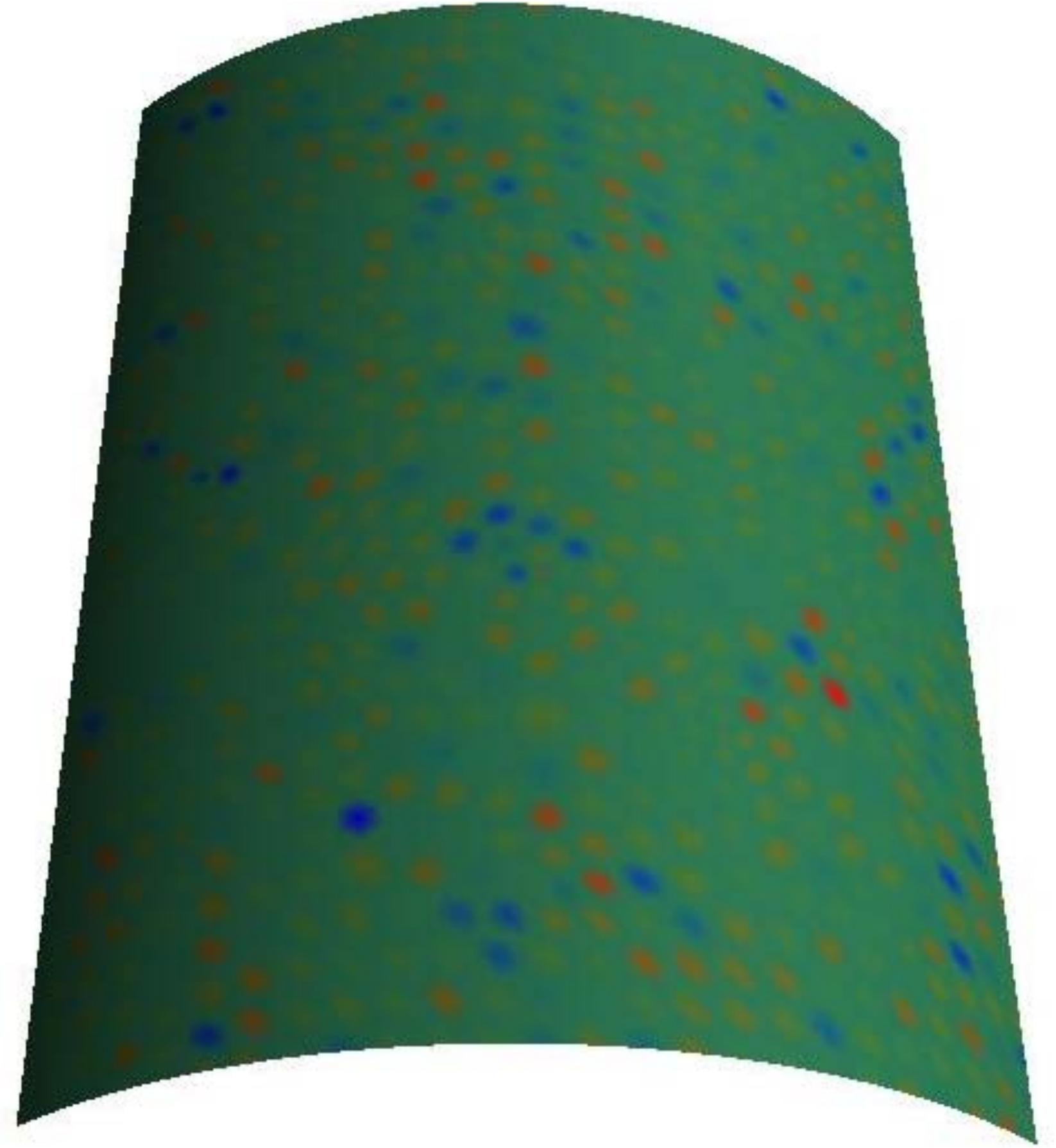}  \\
($a'$) &  ($b'$) &   ($c'$)  \\
\includegraphics[scale=0.23]{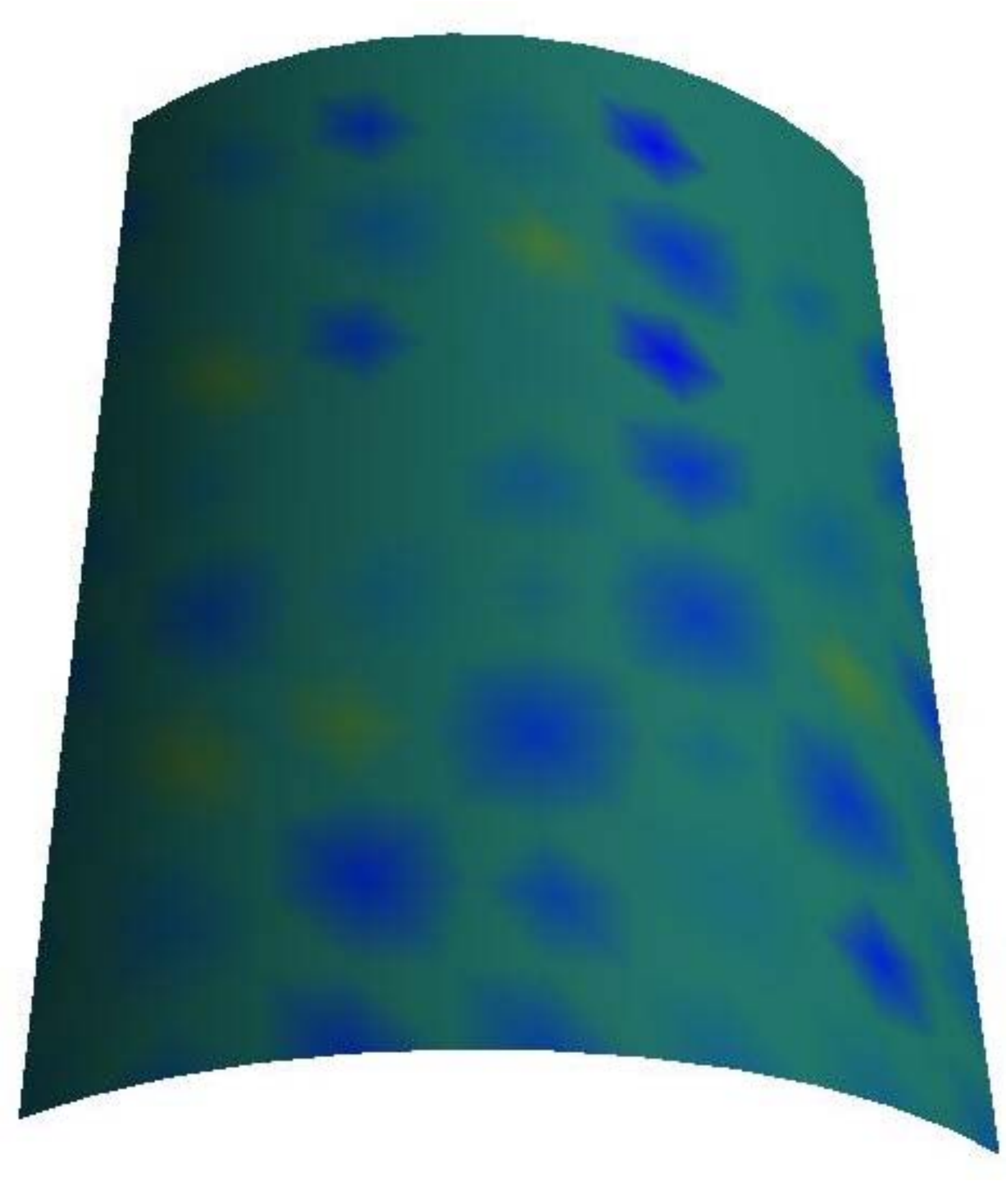}  & \ \ \ \ \ \ \
\includegraphics[scale=0.23]{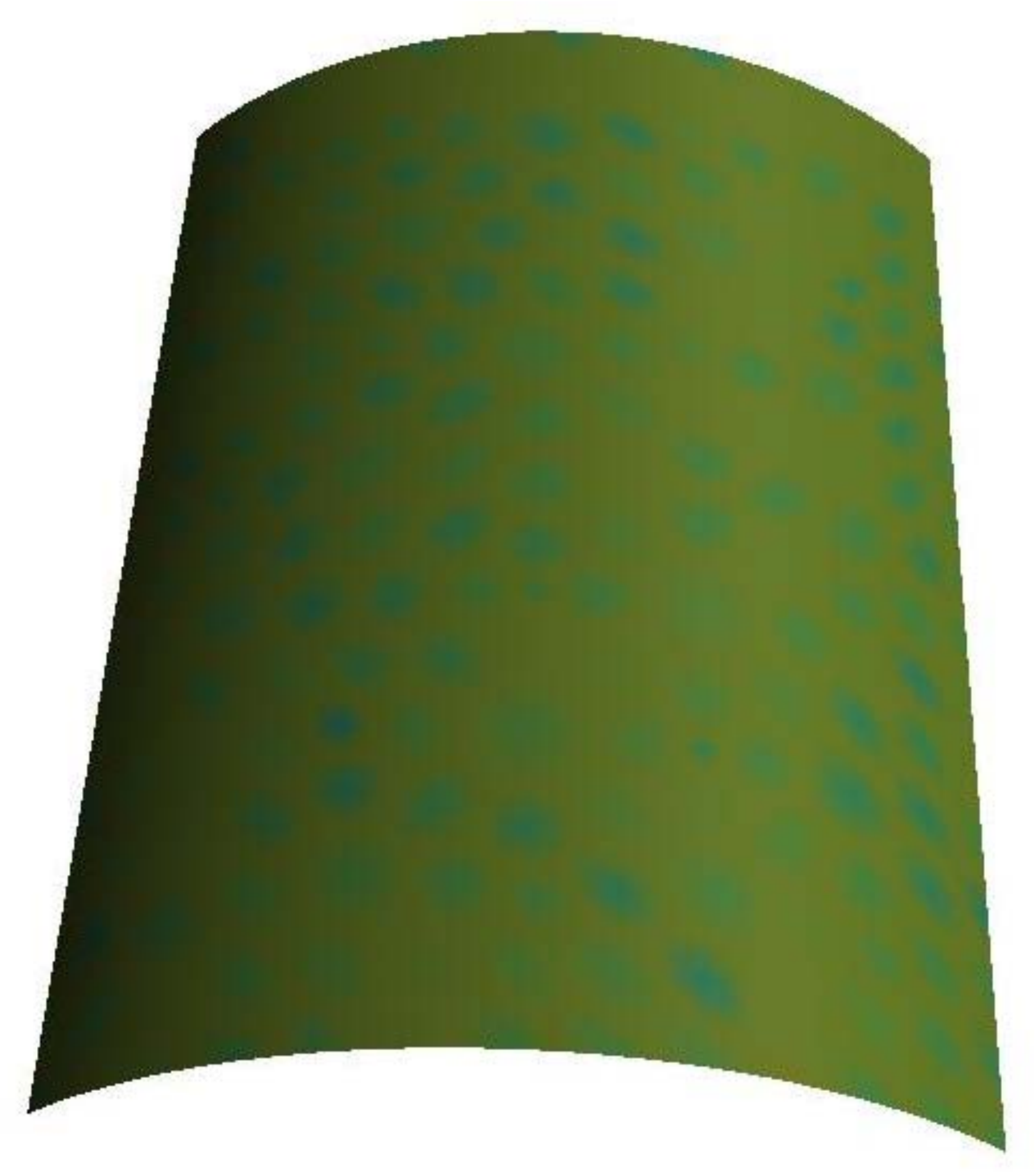}  & \ \ \ \ \ \ \
\includegraphics[scale=0.23]{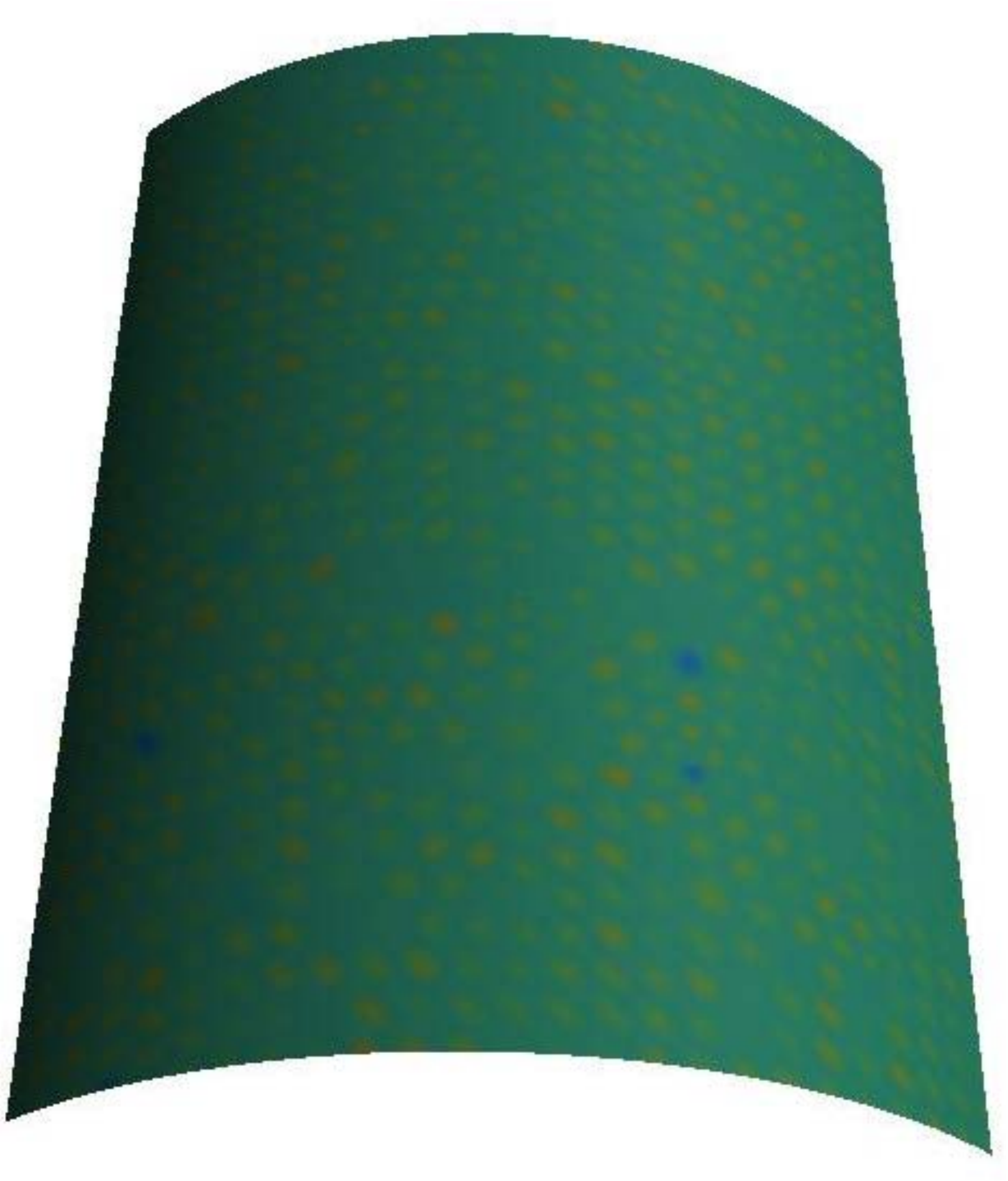}  \\
($a''$) &  ($b''$) & ($c''$)  \\
\end{tabular}
}  \caption{{\small Surface ${\cal S}_1$. Top row: ($a$), ($b$)
and ($c$) show control meshes from coarse to refined, where only once refinement is
implemented from ($a$) to ($b$), and ($b$) to ($c$). Middle row: The corresponding error distribution of
$u-u^h$ results from FEM-Linear shown in ($a'$), ($b'$), and ($c'$). Bottom row: The
corresponding error distribution of the proposed IGA-Loop is respectively shown in ($a''$),
($b''$), and ($c''$).}} \label{quacy}
\end{figure}

\begin{figure}[htb!] \centerline {
\begin{tabular}{c}
\includegraphics[scale=0.7]{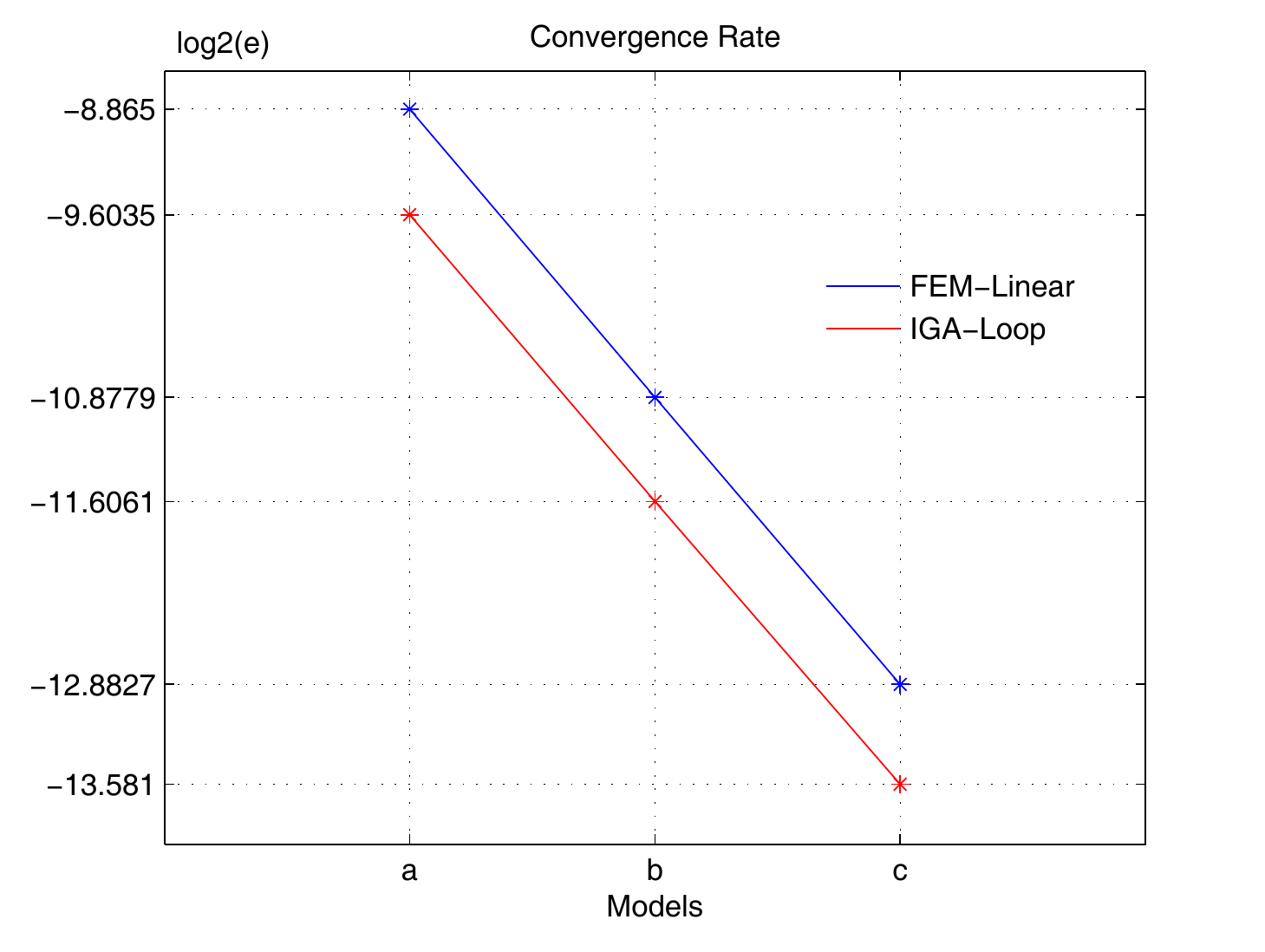}
\end{tabular}
} \caption{{\small Surface ${\cal S}_1$. Comparison of the convergence rate of
 the errors versus the refinement times between FEM-Linear and the proposed
 IGA-Loop. Here the numbers a, b and c on the $x$-axis correspond to the models
 of Figure \ref{quacy} ($a$), ($b$)  and ($c$) respectively, and $e$ on
 the $y$-axis is the $L^2$-norm error.}} \label{errquacy}
\end{figure}

We solve the above problem based on three different control meshes from coarse to refined, as shown
in the top row of Figure \ref{quacy}. The total numbers of vertices/patches are 221/384, 825/1536, and 3185/6144
respectively. They have the same limit surface ${\cal S}_1$, namely, one quarter
of a cylinder. The computational accuracy is indicated by the error distribution $u-u^h$.
As shown in Figure \ref{quacy}, the accuracy of IGA-Loop is higher than that of FEM-Linear. In addition, we also
compare the convergence rate of IGA-Loop against successive refinement. The error is progressively
decreasing as the control mesh becoming finer and finer as shown in Figure \ref{errquacy}. Apart from the visual comparison in
Figure \ref{quacy}, the quantitative comparison of the error between
IGA-Loop and FEM-Linear is performed in Figure \ref{errquacy}. As we can see that IGA-Loop has
higher accuracy.

Consider the second example as the same Laplace-Beltrami harmonic equation as the first one,
but with different exact solution
\[
u= xyz
\]
on a one-eighth of the sphere ${\cal S}_2= \{(x,y,z): x^2+y^2+z^2 =1 \ \& \ x \geq 0 \ \& \  y \geq 0 \ \& \ z \geq 0 \}$
for suitable $f$.

\begin{figure}[ht!]
\centerline {
\begin{tabular}{ccc}
\includegraphics[scale=0.175]{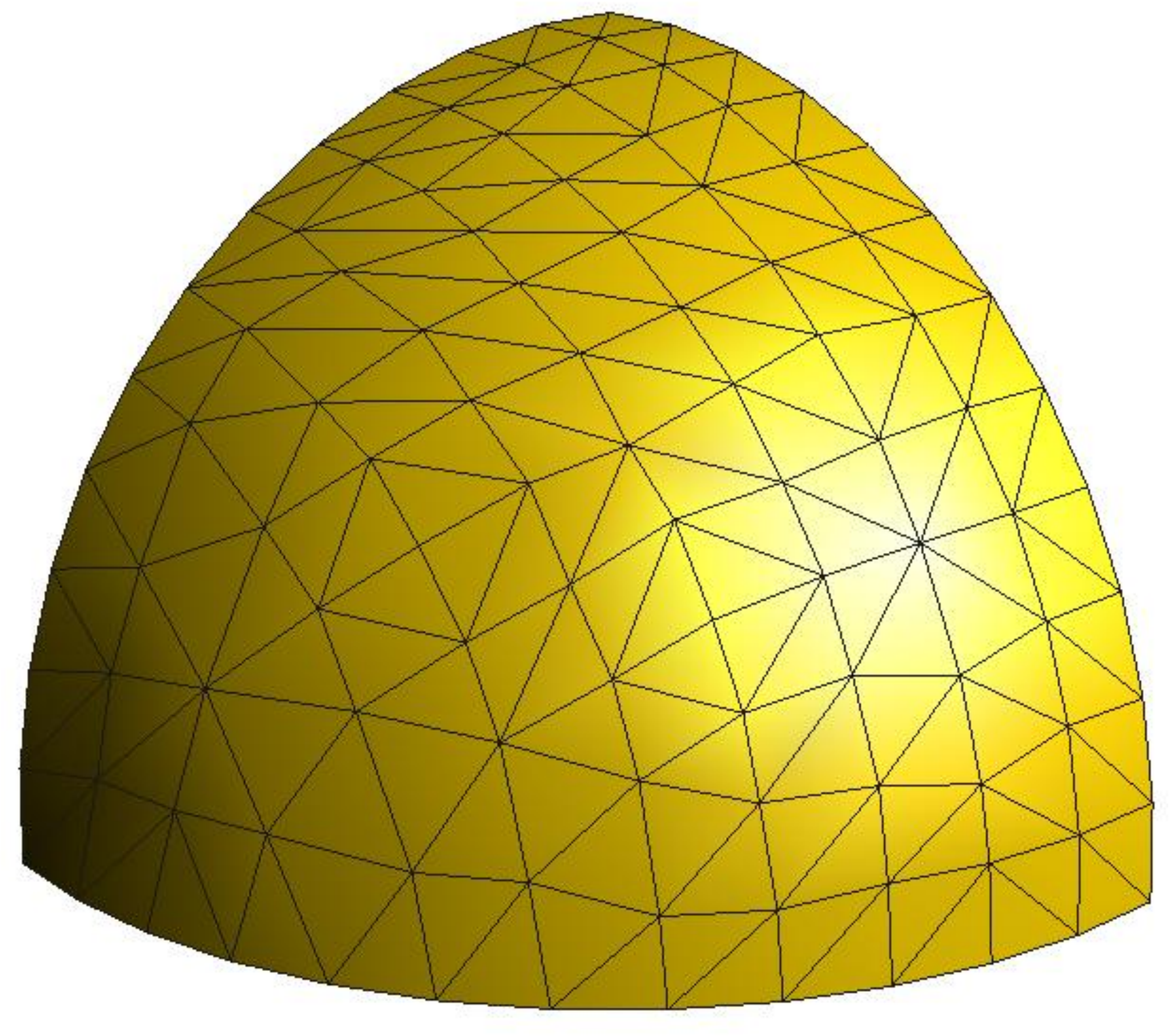}  &
\includegraphics[scale=0.175]{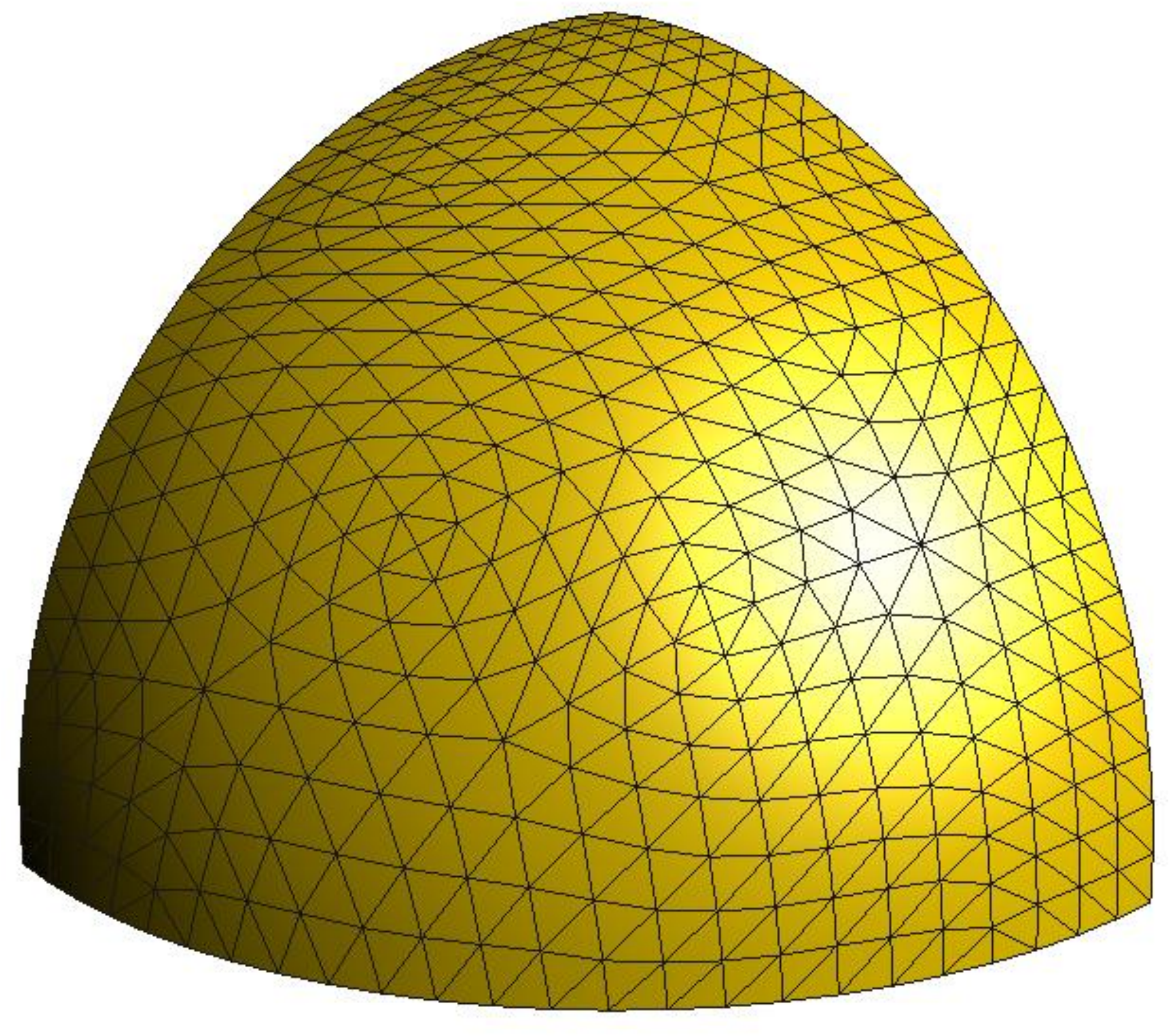}  &
\includegraphics[scale=0.175]{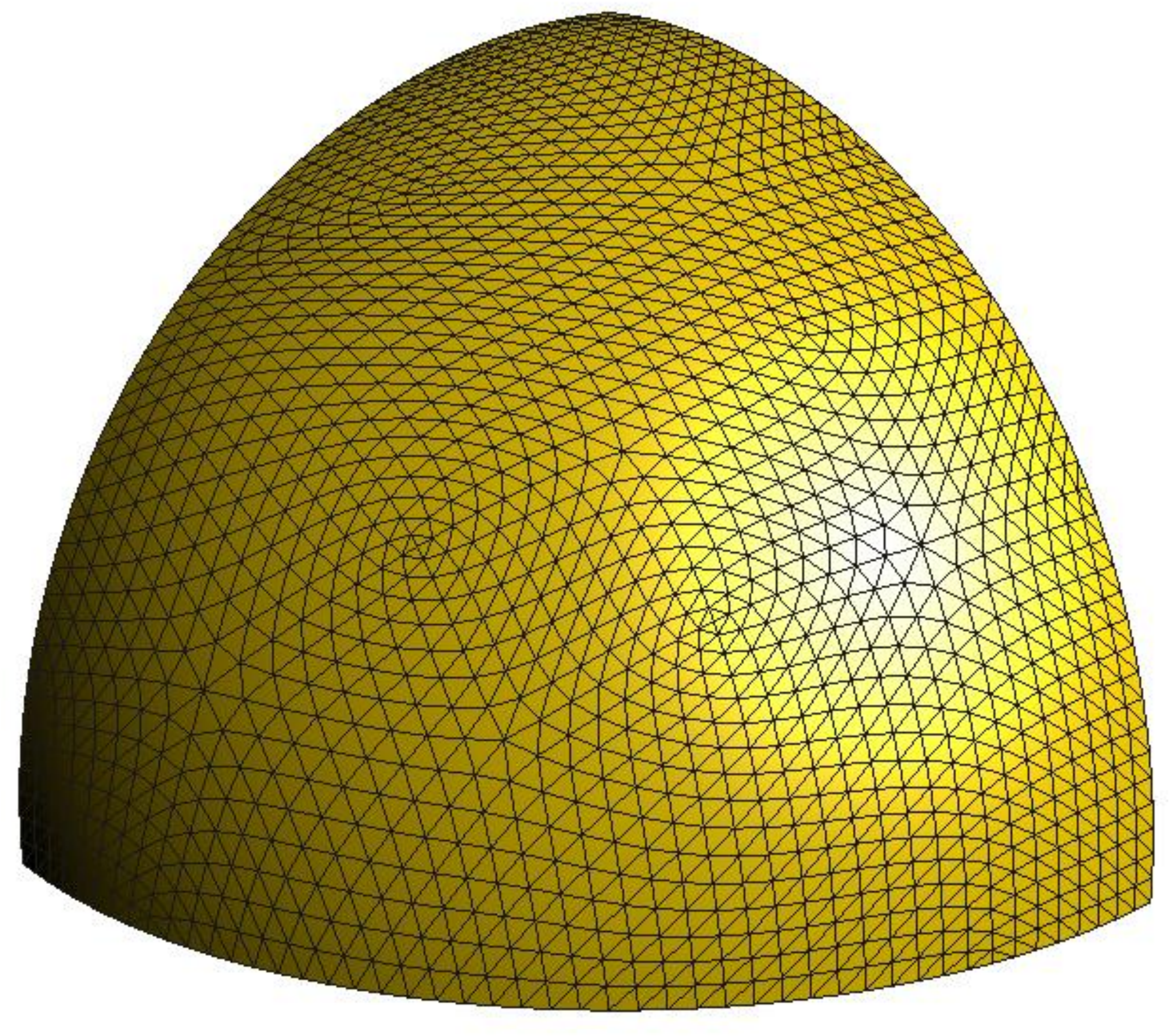}  \\
($a$) &    ($b$)  &   ($c$)  \\
\includegraphics[scale=0.185]{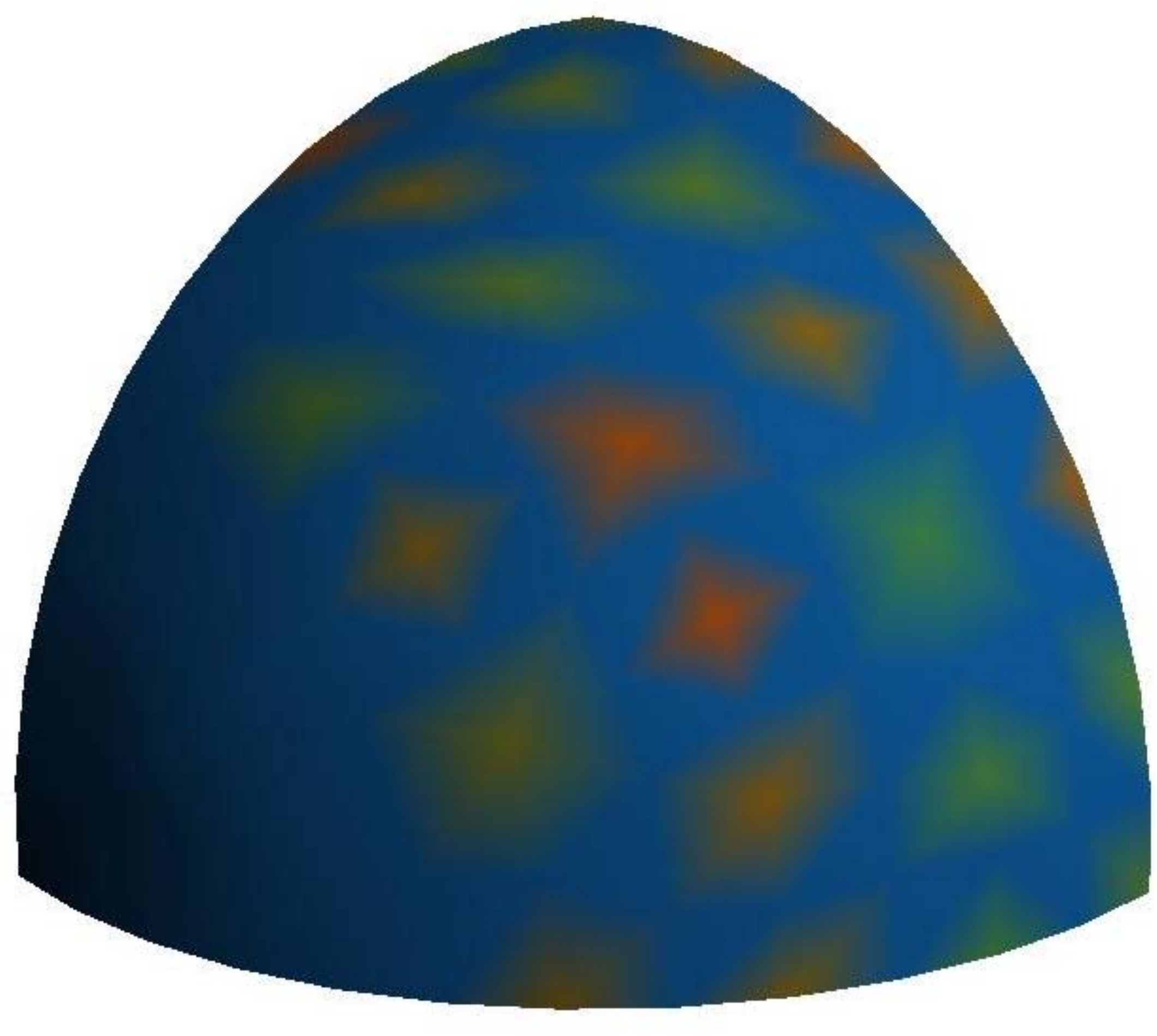}  &
\includegraphics[scale=0.185]{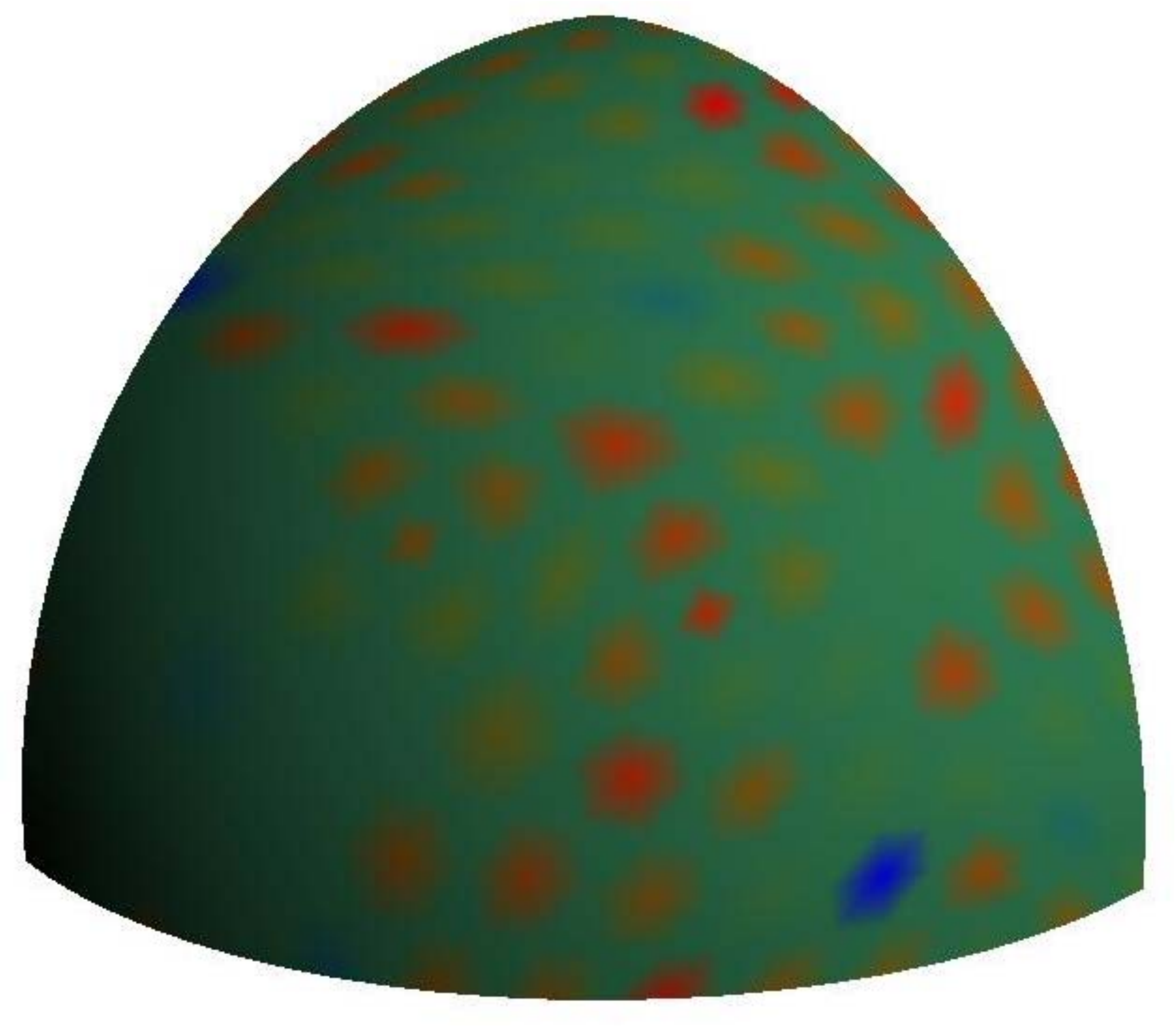}  &
\includegraphics[scale=0.185]{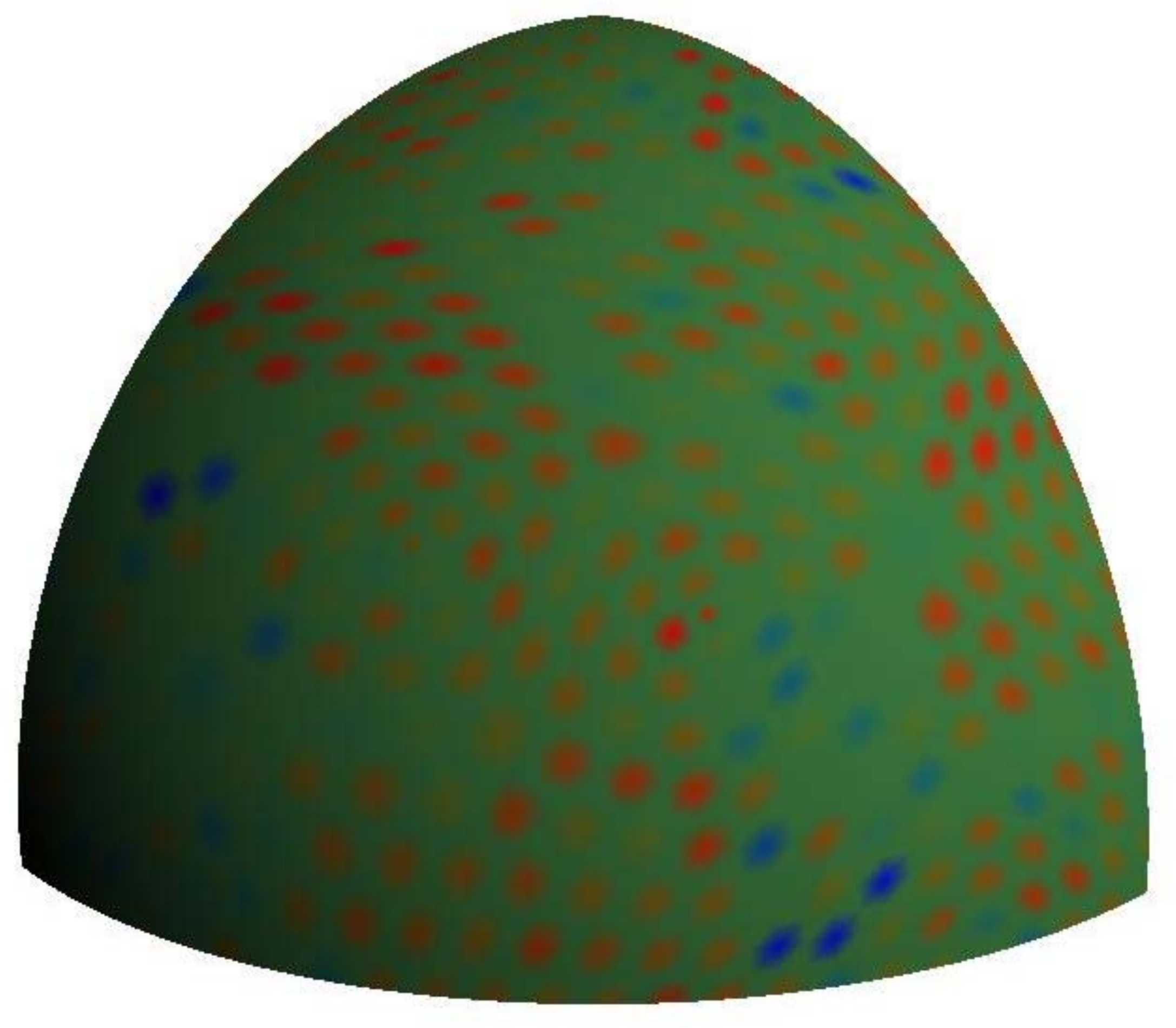}  \\
($a'$) &  ($b'$) &   ($c'$)  \\
\includegraphics[scale=0.2]{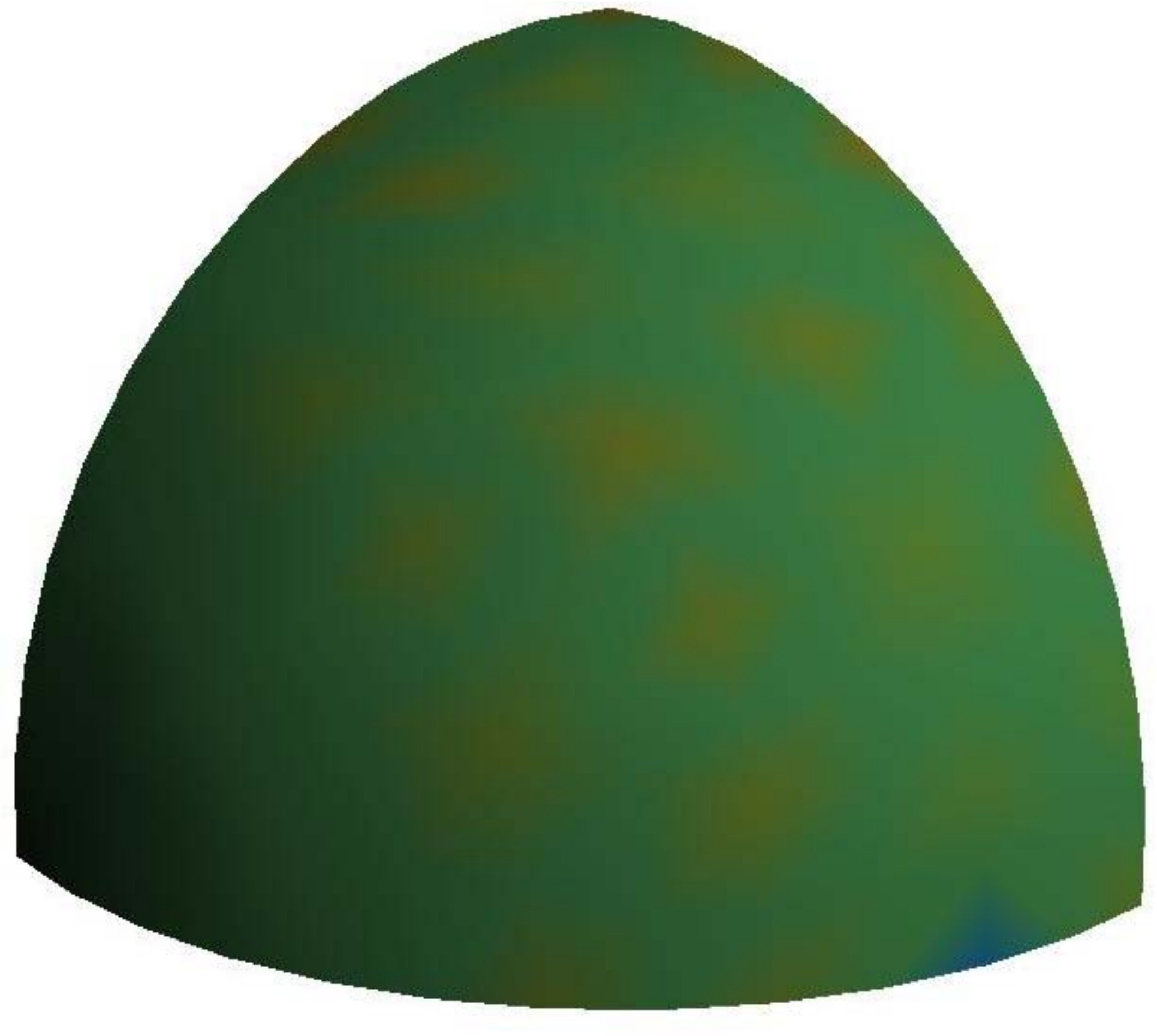}  &
\includegraphics[scale=0.2]{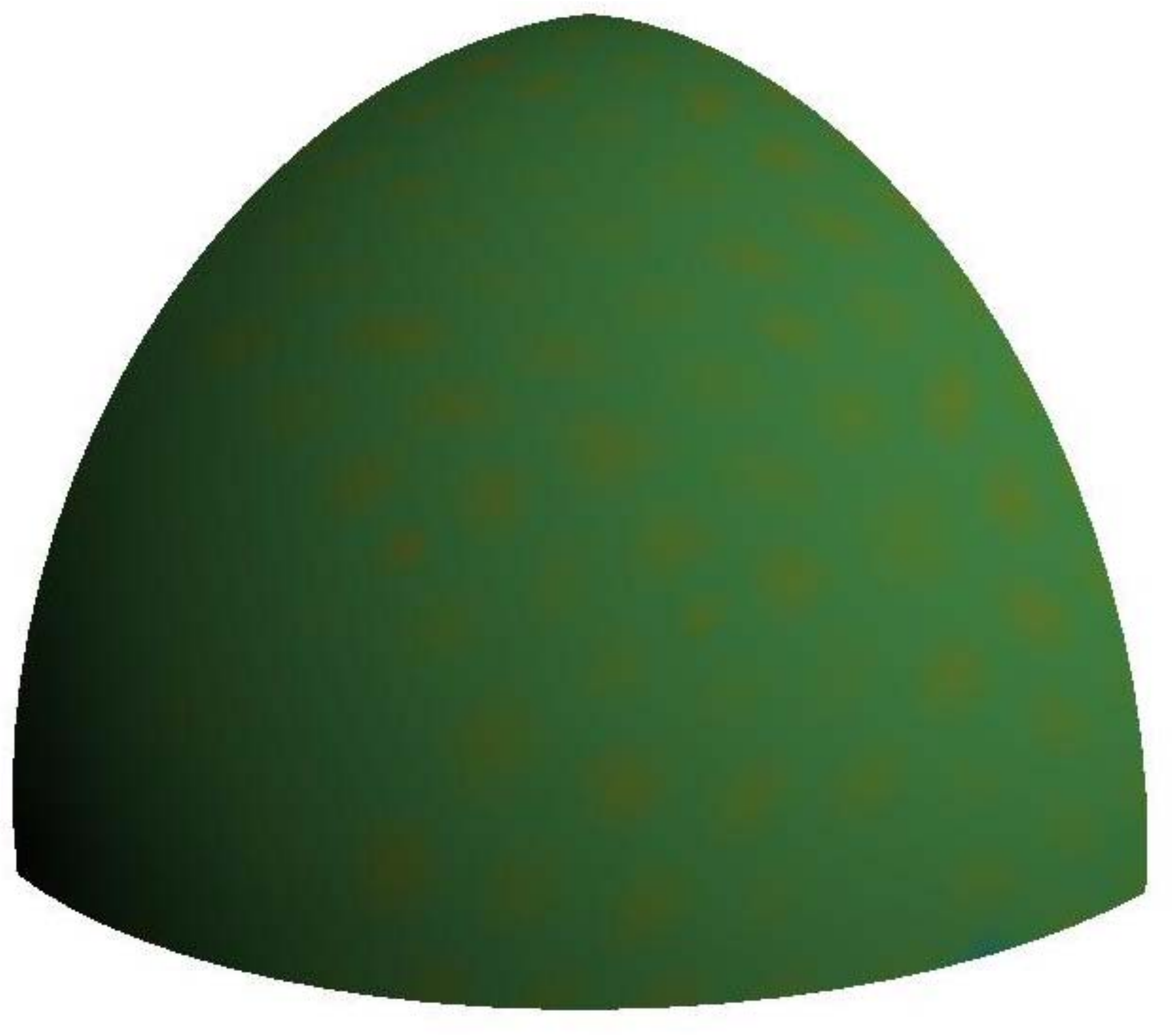}  &
\includegraphics[scale=0.2]{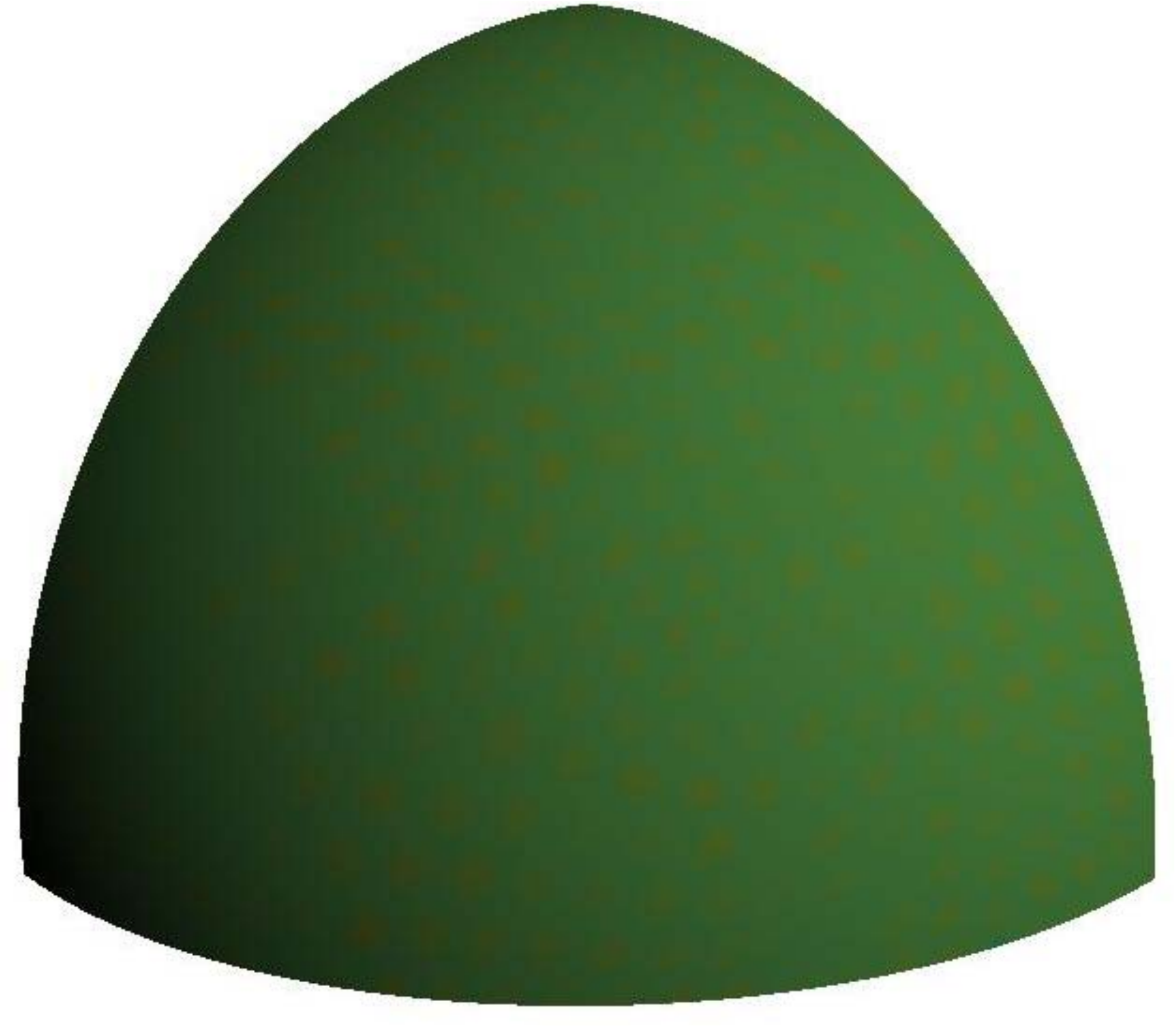}  \\
($a''$) &  ($b''$) & ($c''$)  \\
\end{tabular}
}  \caption{{\small Surface ${\cal S}_2$. Top row: ($a$), ($b$)
and ($c$) are three control meshes where once refinement is
implemented from ($a$) to ($b$), and ($b$) to ($c$). Middle row:
The corresponding distribution of the error
$u-u^h$ results from FEM-Linear shown
in ($a'$), ($b'$), and ($c'$). Bottom row: The error distribution $u-u^h$ results
from IGA-Loop shown in ($a''$), ($b''$), and ($c''$).}} \label{quasph}
\end{figure}

\begin{figure}[htb!] \centerline {
\begin{tabular}{c}
\includegraphics[scale=0.7]{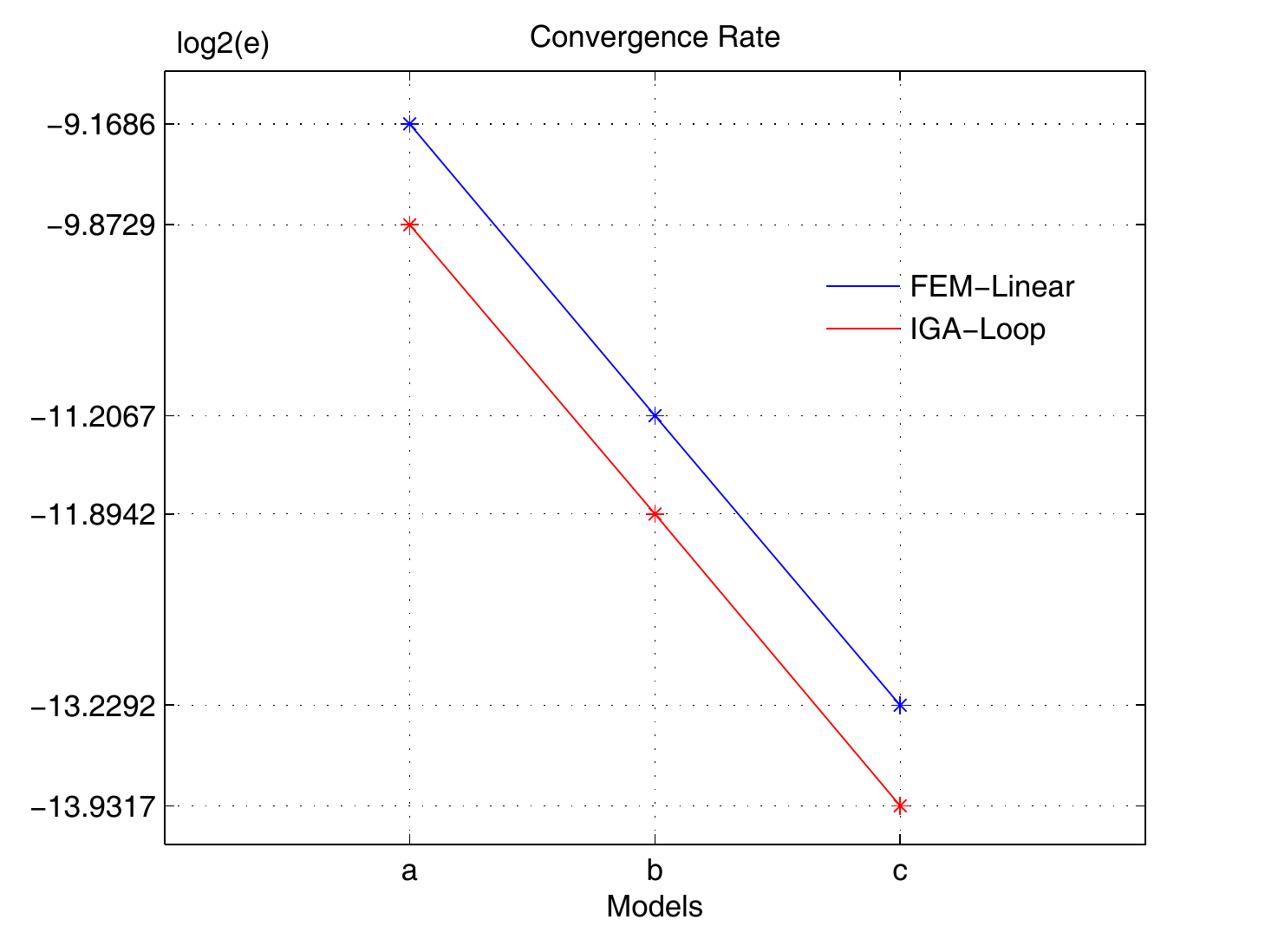}
\end{tabular}
} \caption{{\small Surface ${\cal S}_2$. Comparison of the convergence rate of
 the errors versus the refinement times between FEM-Linear and IGA-Loop. Here the numbers a, b and c on the $x$-axis correspond to the models of Figure \ref{quasph} ($a$), ($b$)  and ($c$) respectively, and $e$ on the $y$-axis is the $L^2$-norm error.}} \label{errquasph}
\end{figure}

As before, we compute the above surface PDE based on
three different control meshes from coarse to refined, as shown
in the top row of Figure \ref{quasph}, which has the same limit
surface ${\cal S}_2$. The total numbers of vertices/patches are 107/176, 389/704, and 1481/2816
respectively. For numerical comparison, FEM-Linear
is used to solve the same problem. As shown in Figure \ref{quasph},
the accuracy of IGA-Loop is higher than that of FEM-Linear by the profile
of error distribution. We further compare the convergence rate of the proposed method against
successive refinement. The error is gradually decreasing as the control mesh
becoming finer and finer as shown in Figure \ref{errquasph}. The quantitative comparison
of the error between the proposed method and FEM-Linear is given in Figure \ref{errquasph}.
As we can see that IGA-Loop has higher accuracy.

In Figure \ref{errquacy} and \ref{errquasph}, the errors $e = \|u-u^h\|_{L^2({\cal S})}$ are computed with the $L^2$-norm of the surface operator,
where both methods have the convergence rate 2. From the comparison of the approximation errors
involved in FEM-Linear and IGA-Loop discretizations, for which the errors versus
the number of subdivision times are plotted, we observe the former is approximately 1.5 times more than the latter. The FEM-Linear
discretization requires the number of degree of freedoms more than the IGA-Loop method, which means IGA-Loop is potentially more efficient.
The explanation of the phenomenon is that the approximation error is introduced by the discretization of the surface geometries using the linear finite elements, instead, exact representation of the surfaces can be achieved by means of IGA-Loop.

We performed the above two numerical tests to demonstrate the mathematical results in Theorem \ref{fem}. For the Laplace-Beltrami biharmonic problem defined on a cylinder and the Laplace-Beltrami triharmonic problem defined on a unit sphere, the error estimate (\ref{surffem}) in norm $L^2({\cal S})$  of Theorem \ref{fem} holds, which will be studied in the following contents.

\subsection{Test Suite 2: Surface Laplace-Beltrami Biharmonic Equation}

In what follows, we solve a fourth-order example of the Laplace-Beltrami biharmonic equation as
\begin{equation}
\label{cylinder_2}
\left\{
\begin{array}{l}
\Delta_{\cal S}^2 u = f, \quad \text{in}~{\cal S},
\\[0.5em]
u = \frac{\partial  u}{\partial {\bm n}} = 0, \quad \text{on}~{\partial\cal S},
\end{array}
\right.
\end{equation}
where the exact solution
\[
u= {\rm sin}^2(2\theta) {\rm sin}^2 (2z)
\]
with $\theta = {\rm atan} (y/x)$ on the
surface ${\cal S}_3= \{(\theta,z): 0 \leq \theta \leq 2\pi \ \& \  0 \leq z  \leq 1 \}$ for suitable $f$.

\begin{figure}[ht!]
\centerline {
\begin{tabular}{ccc}
\includegraphics[scale=0.2]{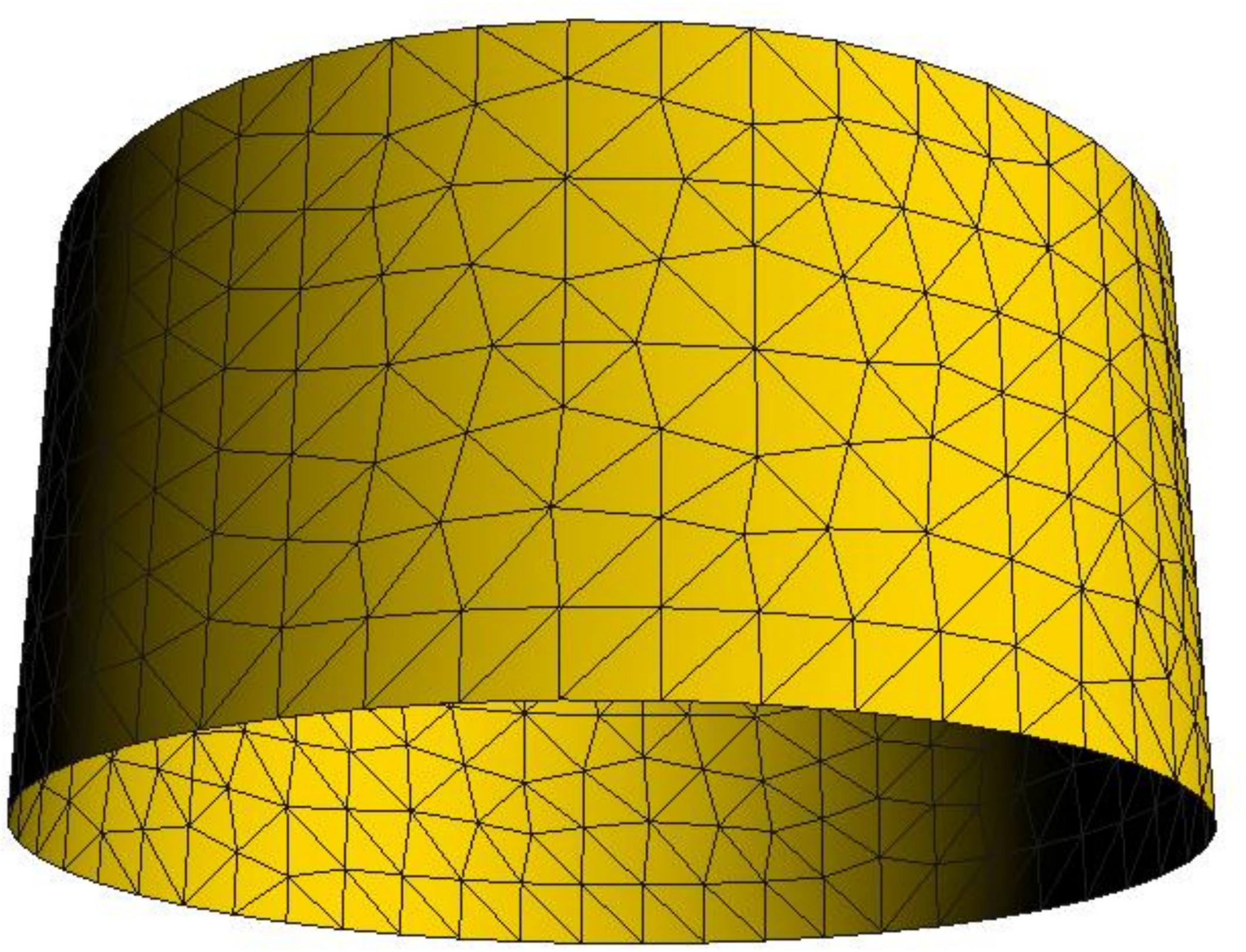}  & \ \ \
\includegraphics[scale=0.2]{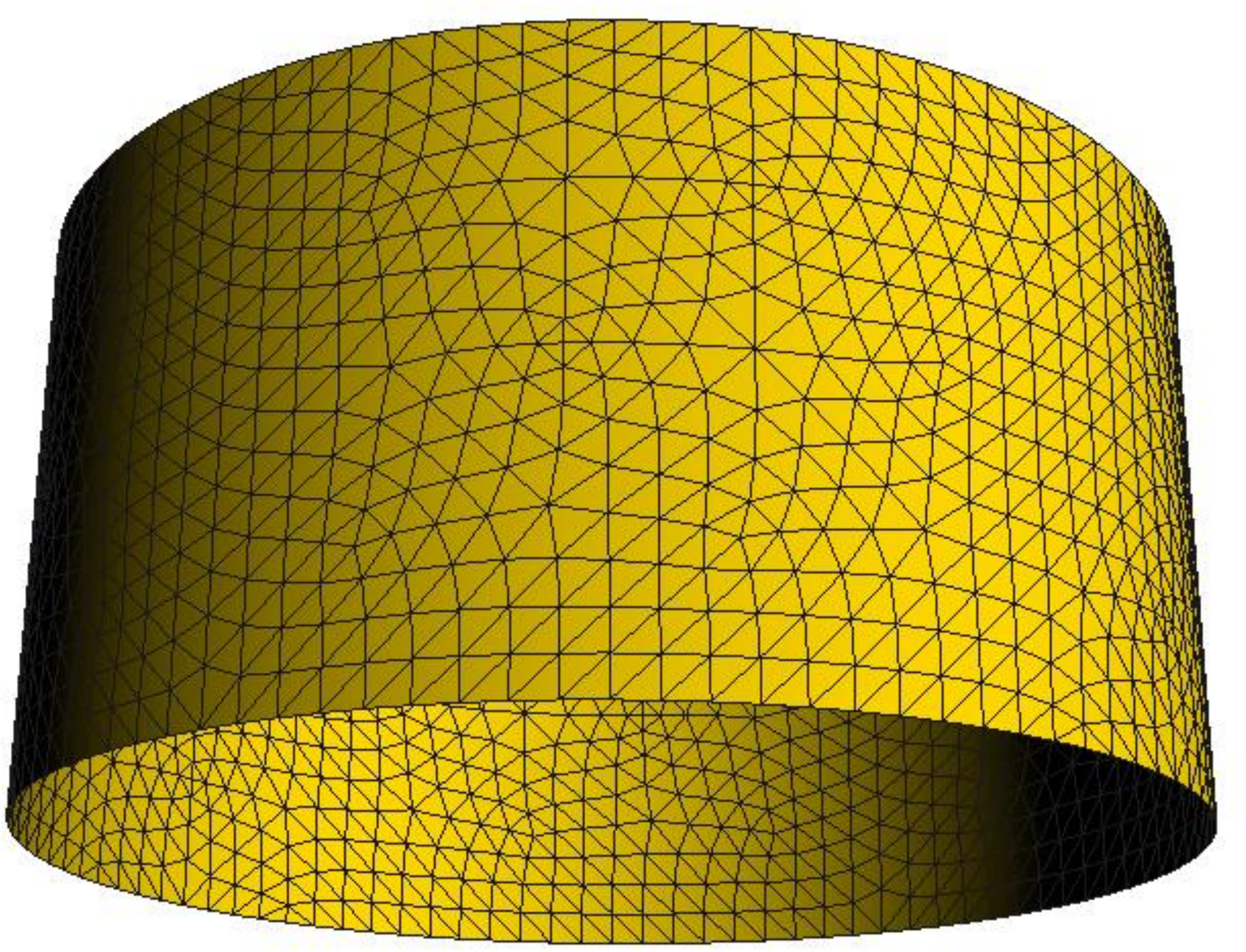}  & \ \ \
\includegraphics[scale=0.2]{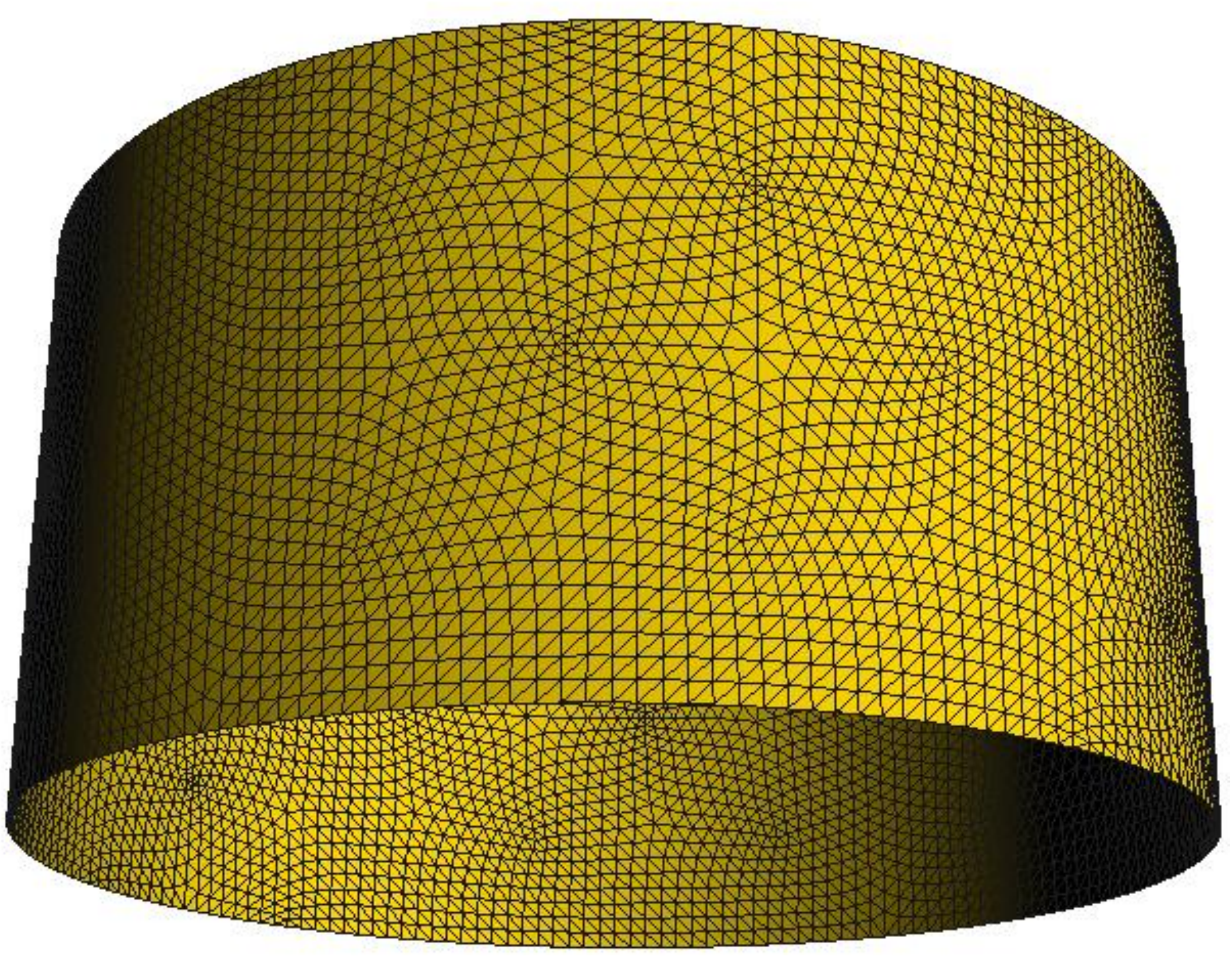}  \\
($a$) &    ($b$)  &   ($c$)  \\
\includegraphics[scale=0.2]{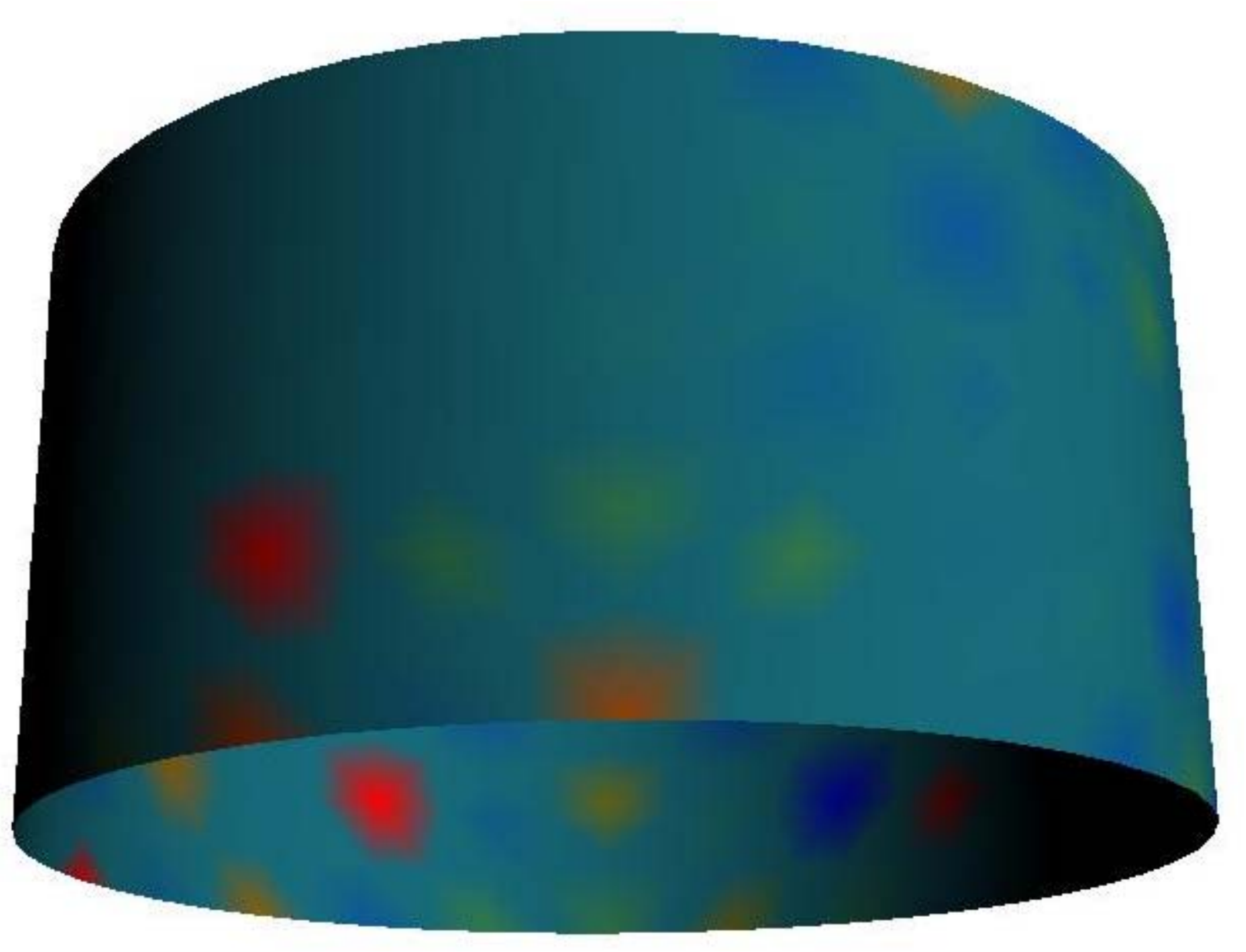}  & \ \ \
\includegraphics[scale=0.21]{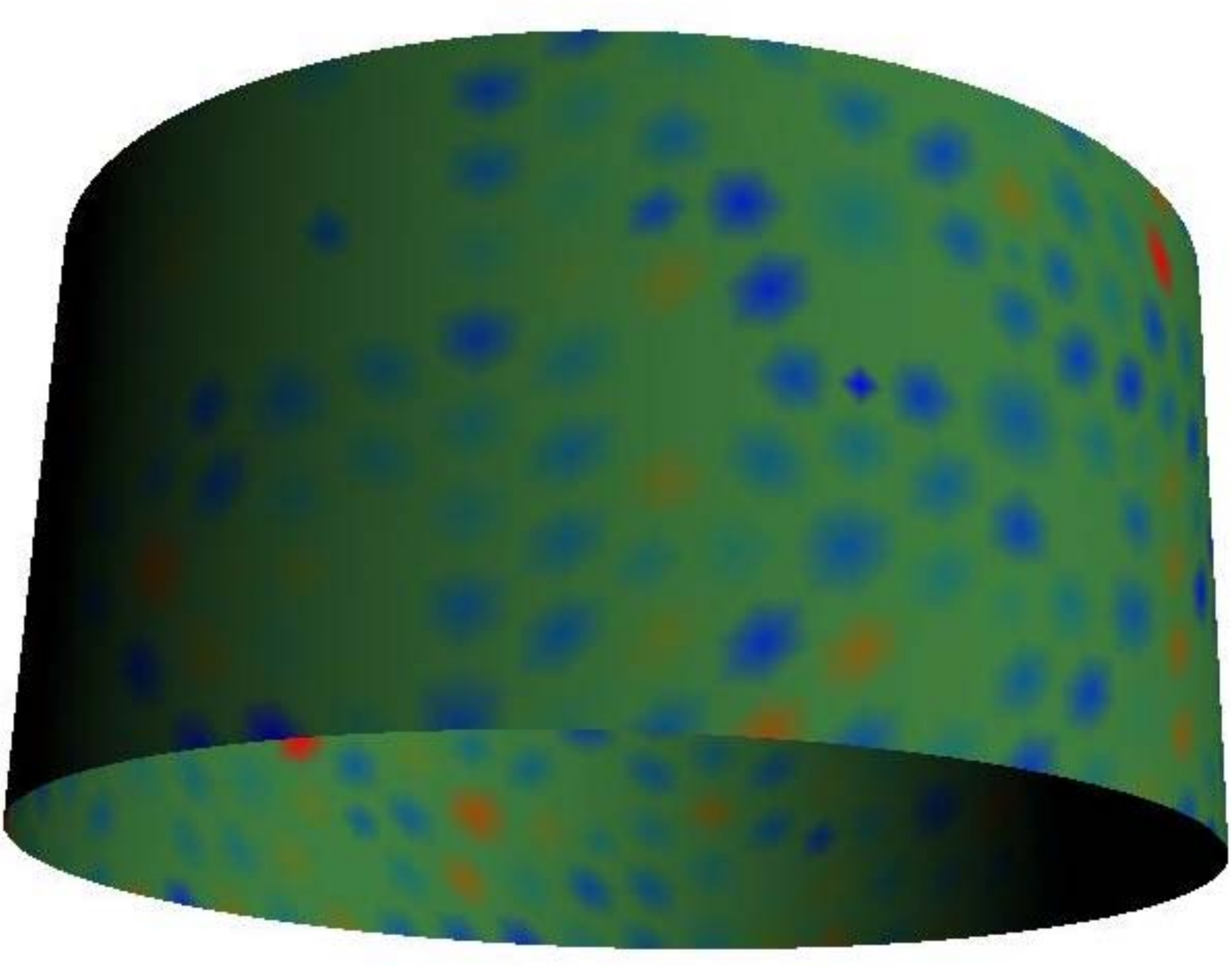}  & \ \ \
\includegraphics[scale=0.21]{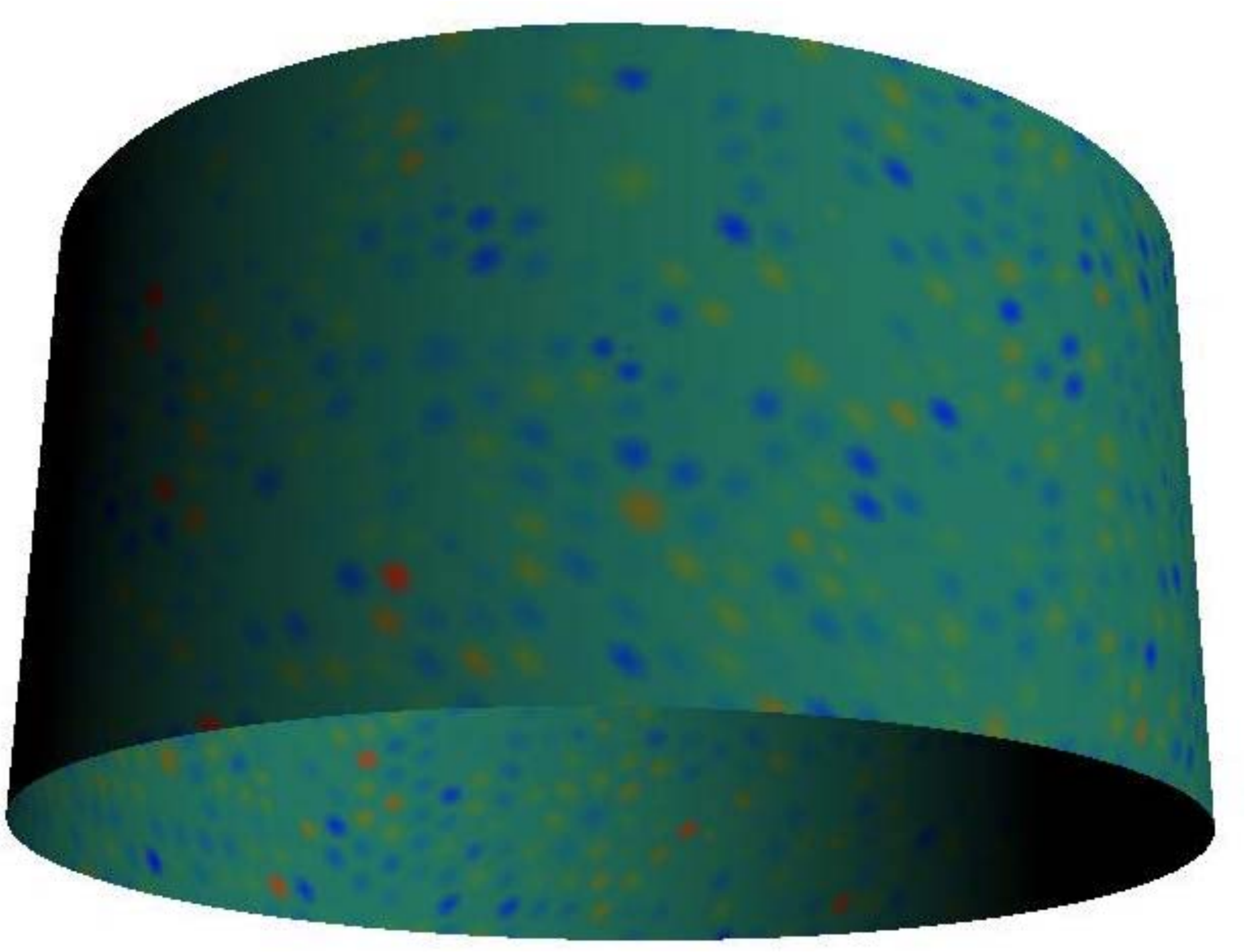}  \\
($a'$) &  ($b'$) &   ($c'$)  \\
\includegraphics[scale=0.2]{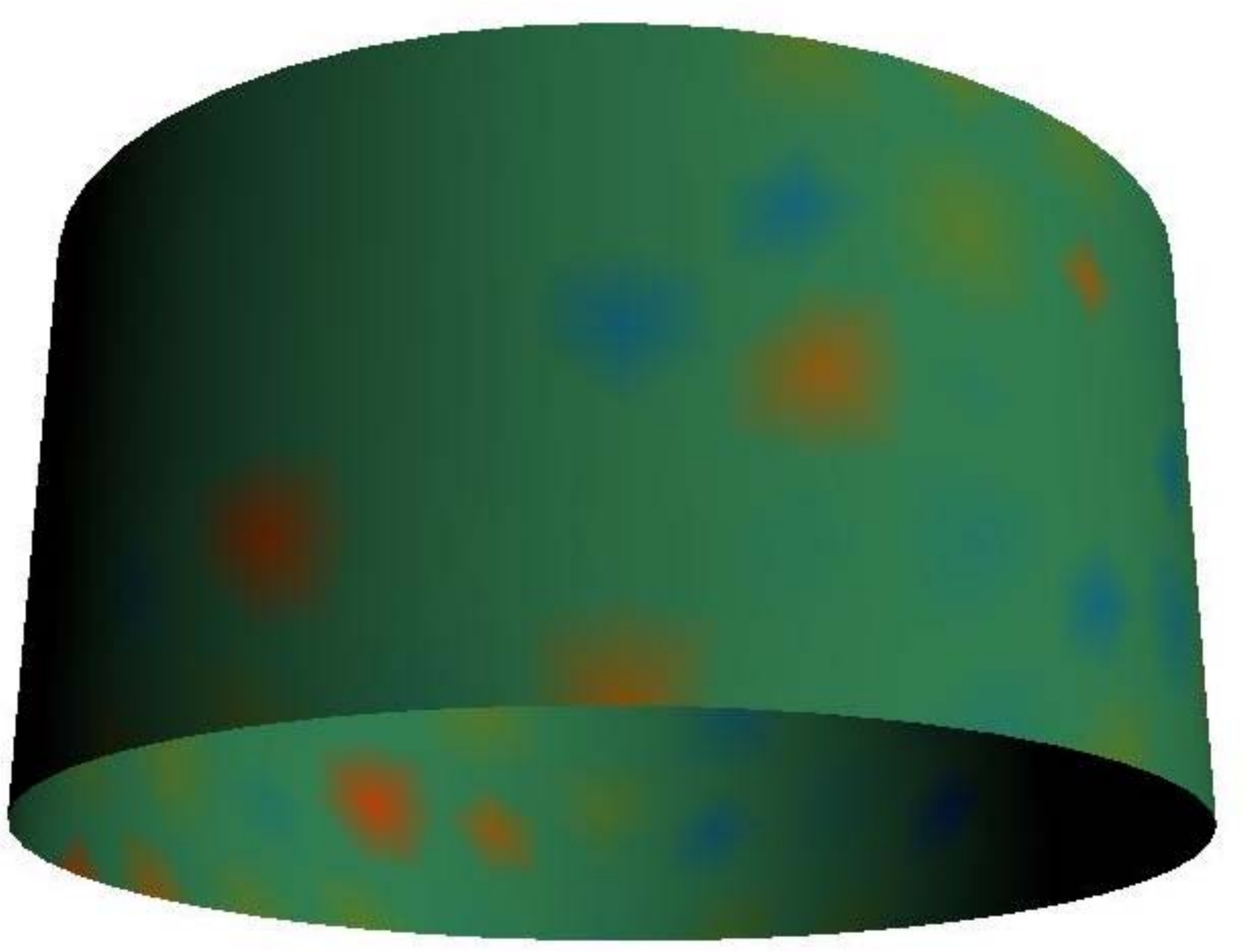}  & \ \ \
\includegraphics[scale=0.2]{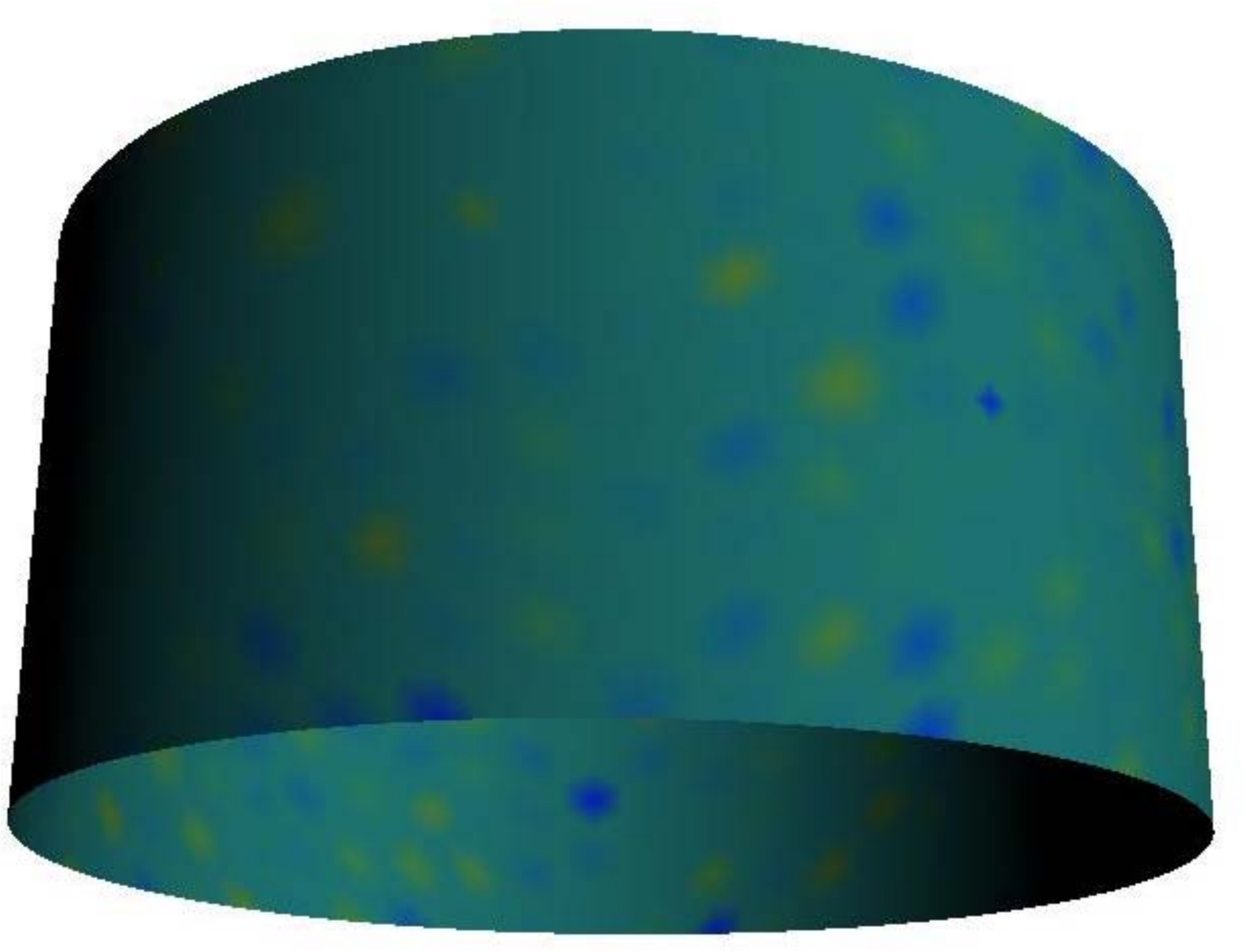}  & \ \ \
\includegraphics[scale=0.2]{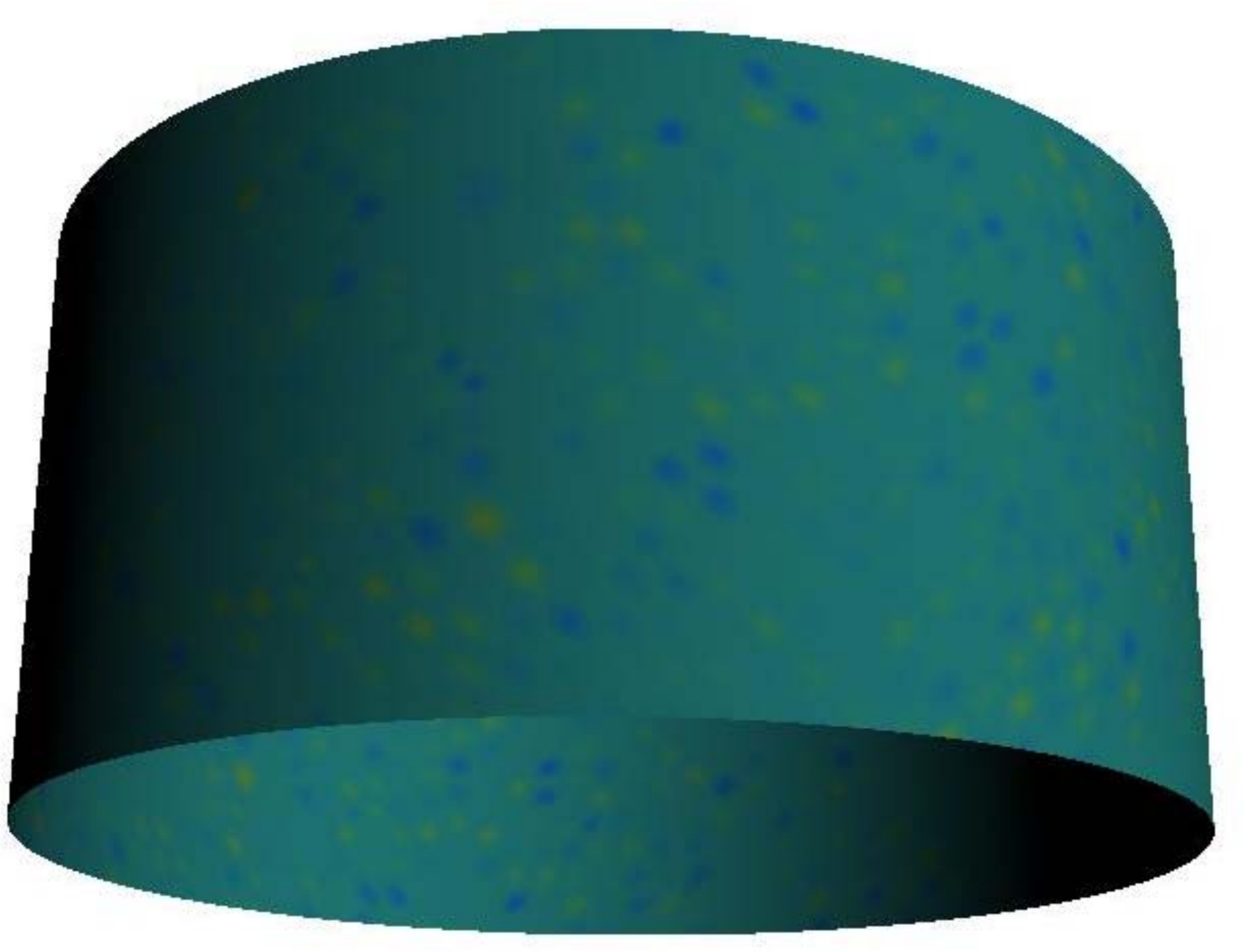}  \\
($a''$) &  ($b''$) & ($c''$)  \\
\end{tabular}
}  \caption{{\small Surface ${\cal S}_3$. ($a$), ($b$)
and ($c$) are three control meshes where one time refinement is
implemented from ($a$) to ($b$), and ($b$) to ($c$). The corresponding distribution of the error
$u-u^h$ resulting from FEM-Linear and our IGA-Loop is respectively shown
in ($a'$), ($b'$), and ($c'$) of the second row, and ($a''$),
($b''$), and ($c''$) of the third row.}} \label{wholecylinder}
\end{figure}

As shown in the top row of Figure \ref{wholecylinder}, three different control meshes from coarse
to refined are given. The total numbers of vertices/patches are 432/768,  1632/3072, and 6336/12288
respectively. They have the same limit surface ${\cal S}_3$ as the extended Loop subdivision proceeds.
We solve the Laplace-Beltrami biharmonic equation based on the proposed IGA-Loop. The computational accuracy is shown by the error distribution $u-u^h$.
As shown in Figure \ref{wholecylinder}, the accuracy of IGA-Loop is higher than that of FEM-Linear.

\begin{figure}[htb!] \centerline {
\begin{tabular}{c}
\includegraphics[scale=0.7]{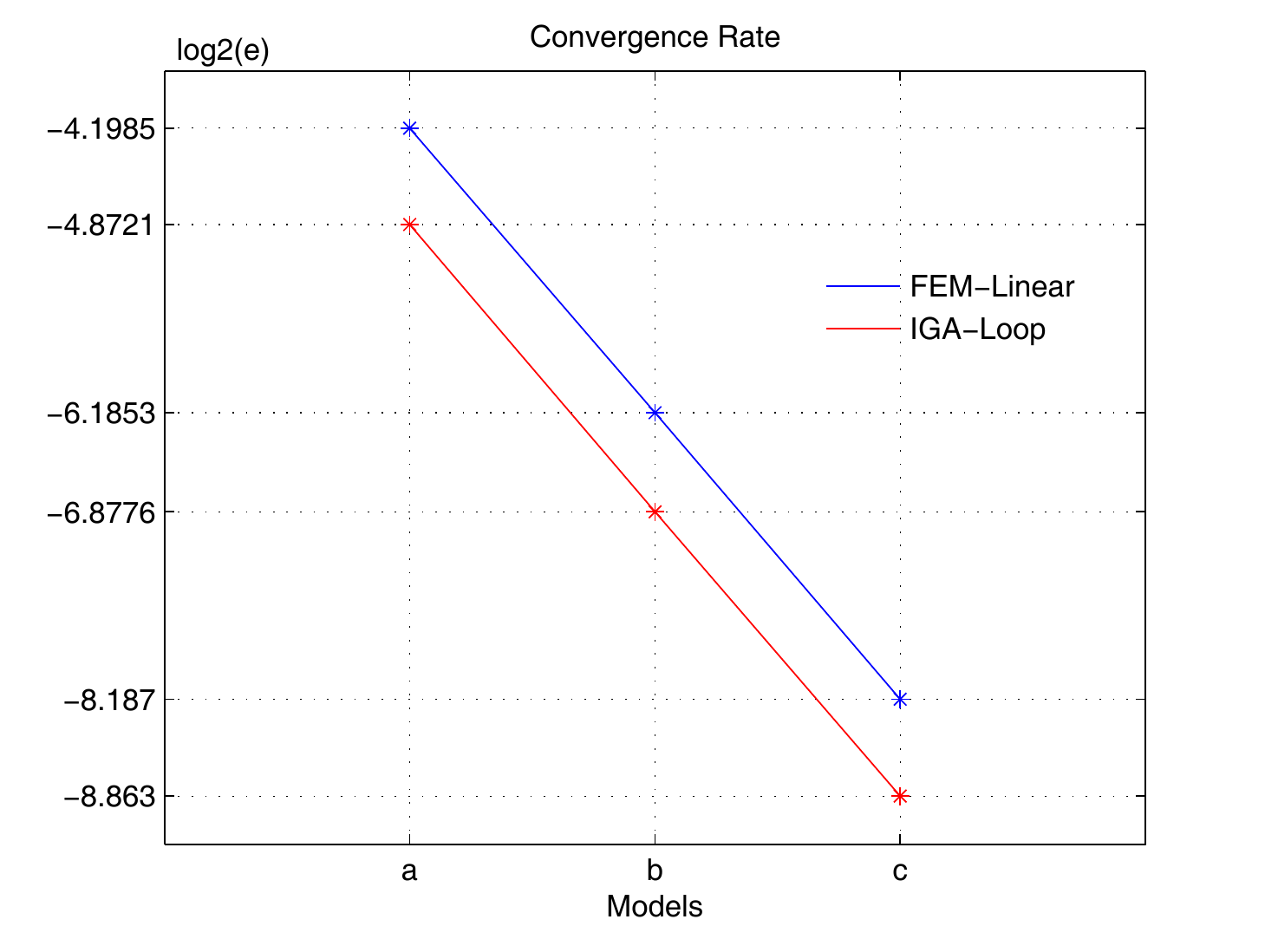}
\end{tabular}
} \caption{{\small Surface ${\cal S}_3$. Comparison of the convergence rate of
 the errors versus the refinement times between FEM-Linear and IGA-Loop. Here the numbers a, b and c on the $x$-axis correspond to the models of Figure \ref{wholecylinder} ($a$), ($b$)  and ($c$) respectively, and $e$ on the $y$-axis is the $L^2$-norm error.}} \label{errcylinder}
\end{figure}

We further compare the convergence rate of IGA-Loop against
successive refinement. The error is gradually decreasing as the control mesh
becoming finer and finer as shown in Figure \ref{errcylinder}. The quantitative comparison
of the error in norm $L^2({\cal S})$  versus the number of subdivision times, obtained with IGA-Loop and FEM-Linear is given in Figure \ref{errcylinder}.
As shown that IGA-Loop has higher accuracy. In Figure \ref{errcylinder}, the errors of norm $L^2({\cal S})$ under three times refinement level is plotted. We observe that the convergence rate is 2 for the error in the sense of norm $L^2({\cal S})$, which is consistent with the error estimate (\ref{surffem}). As shown in Figure \ref{errcylinder}, the approximation errors obtained by the FEM-Linear discretization is approximately 1.6 times more than that of our IGA-Loop discretization. It means that the IGA-Loop approximation only requires a smaller number of degree of freedoms than FEM-Linear to achieve the same accuracy.

\subsection{Test Suite 3: Surface Laplace-Beltrami Triharmonic Equation}

Furthermore, a sixth-order Laplace-Beltrami triharmonic equation is solved as following
\begin{equation}
\label{sphere_1}
- \Delta_{\cal S}^3 u = f,
\end{equation}
where for calculated $f$, the exact solution is
\[
u= xyz
\]
on the closed sphere ${\cal S}_4$ with the radius $r=1$.

\begin{figure}[ht!]
\centerline {
\begin{tabular}{ccc}
\includegraphics[scale=0.2]{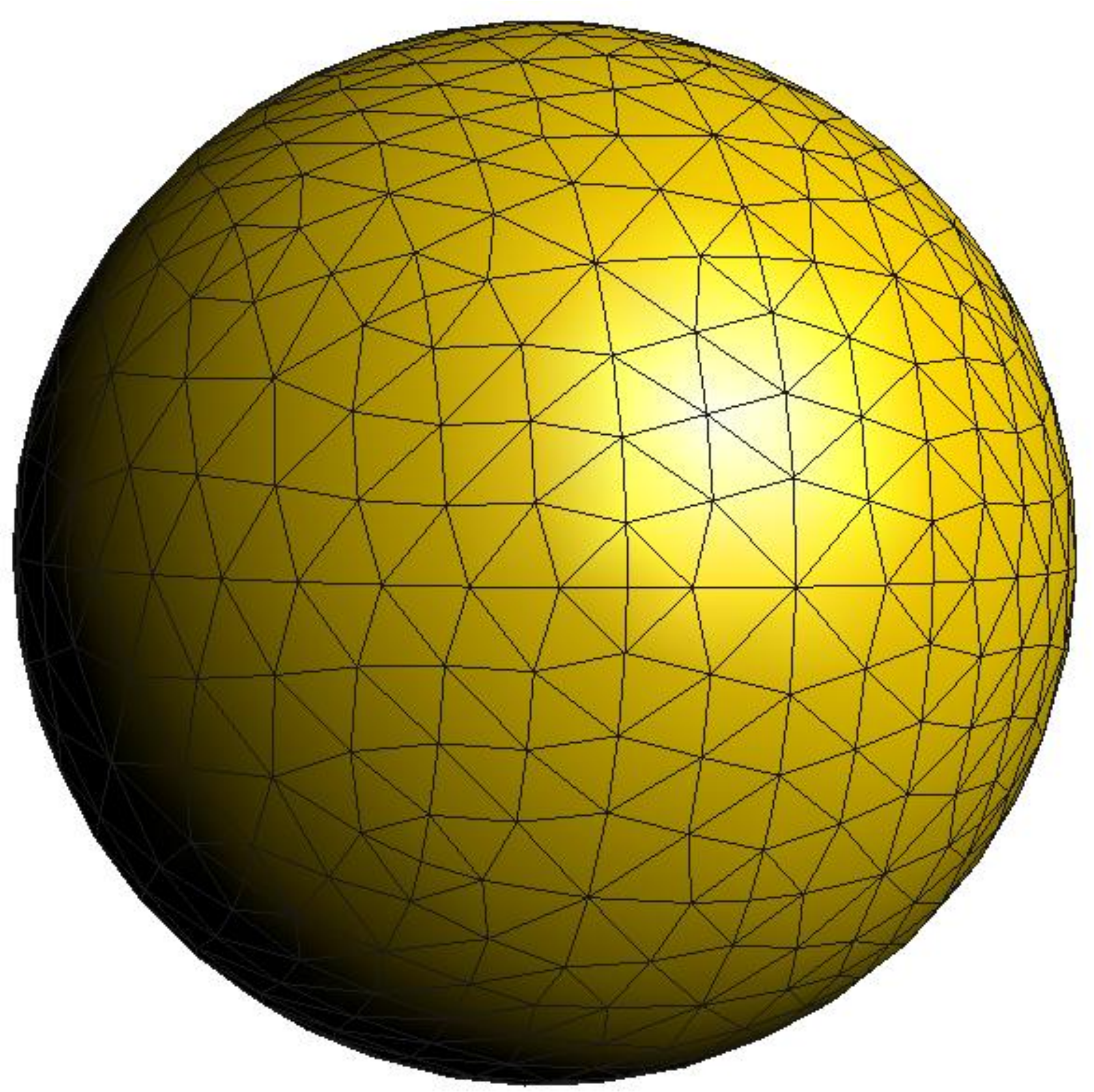}  & \ \ \ \ \
\includegraphics[scale=0.2]{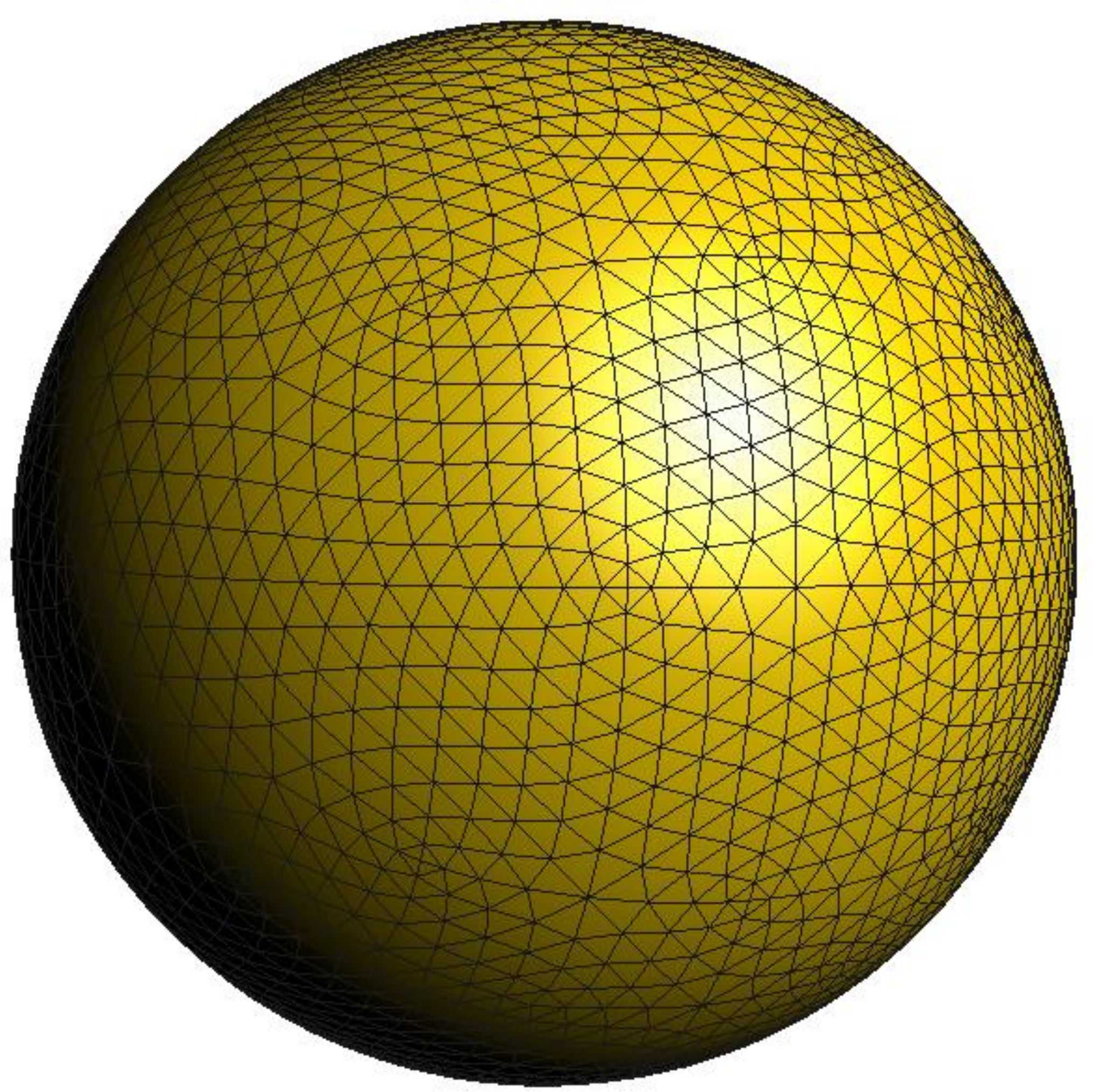}  & \ \ \ \ \
\includegraphics[scale=0.2]{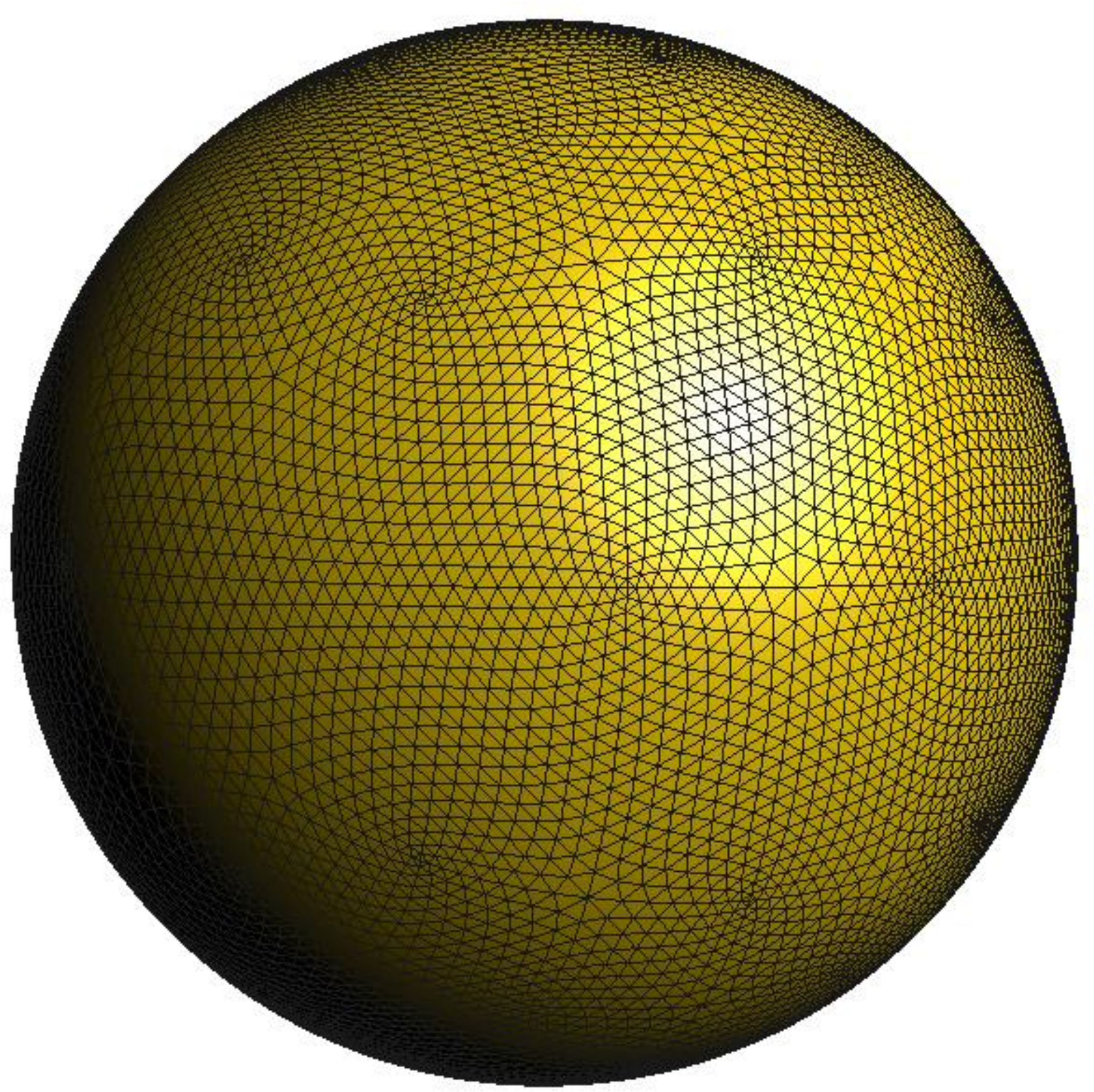}  \\
($a$) &    ($b$)  &   ($c$)  \\
\includegraphics[scale=0.2]{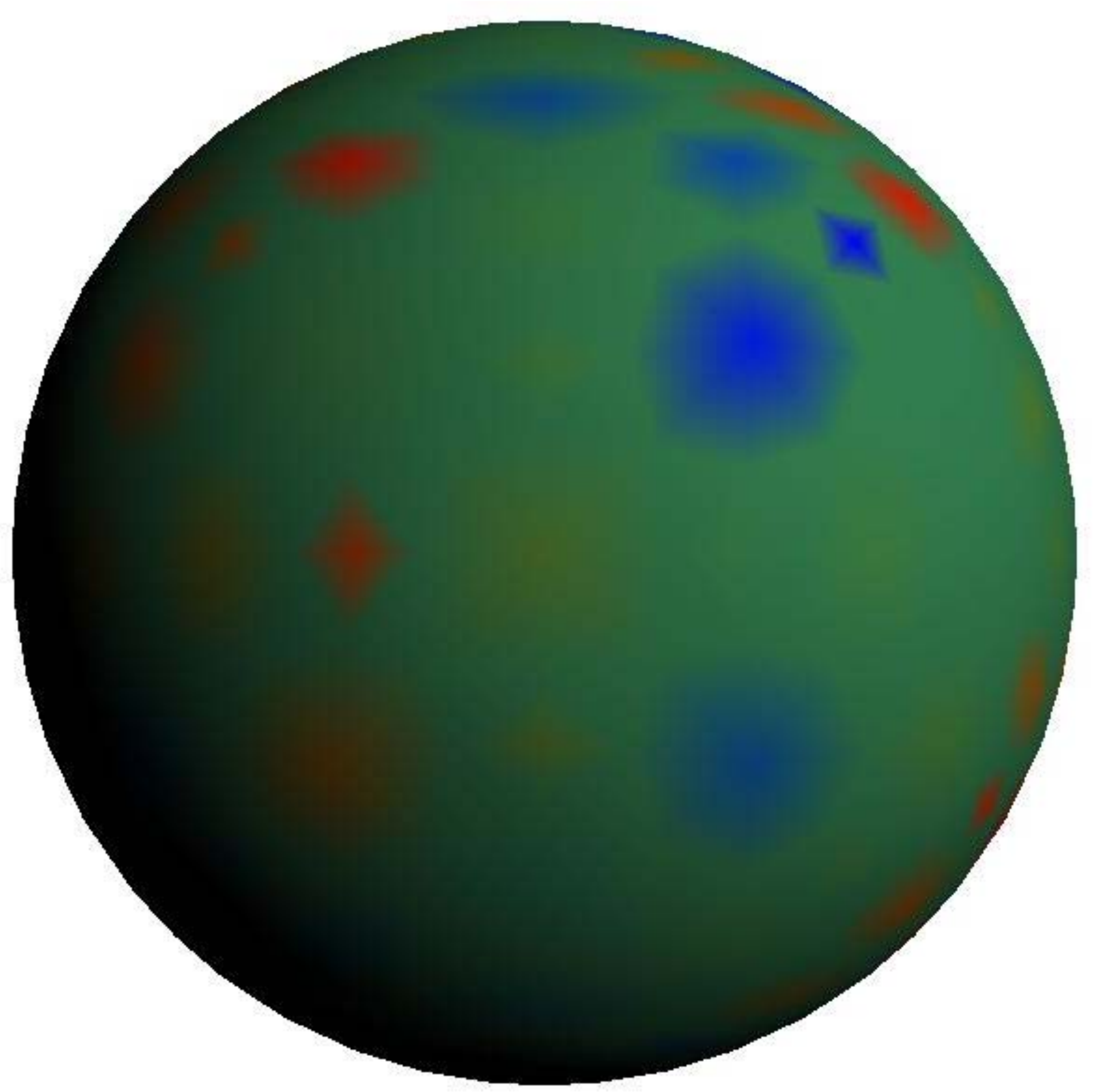}  & \ \ \ \ \
\includegraphics[scale=0.2]{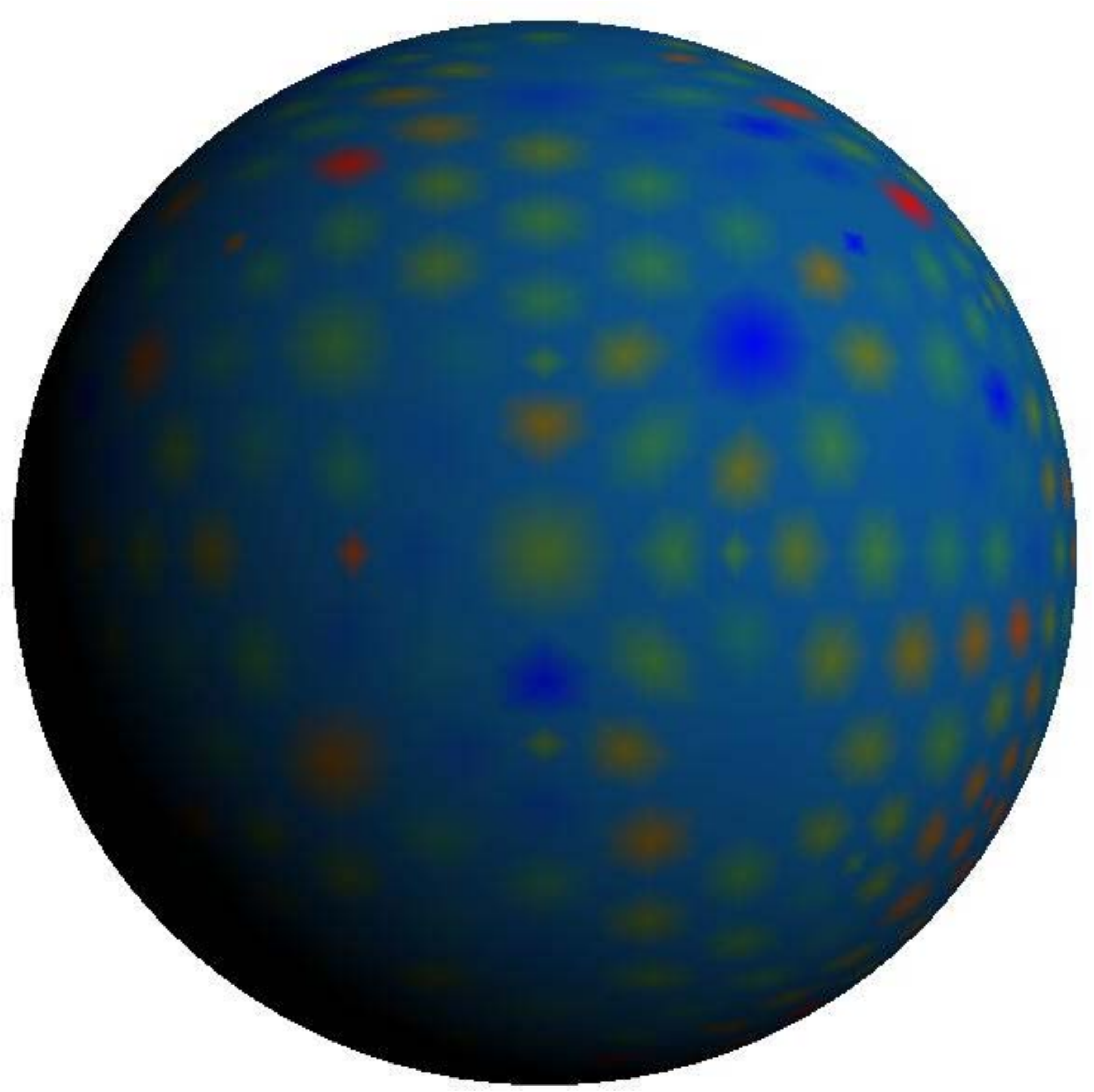}  & \ \ \ \ \
\includegraphics[scale=0.2]{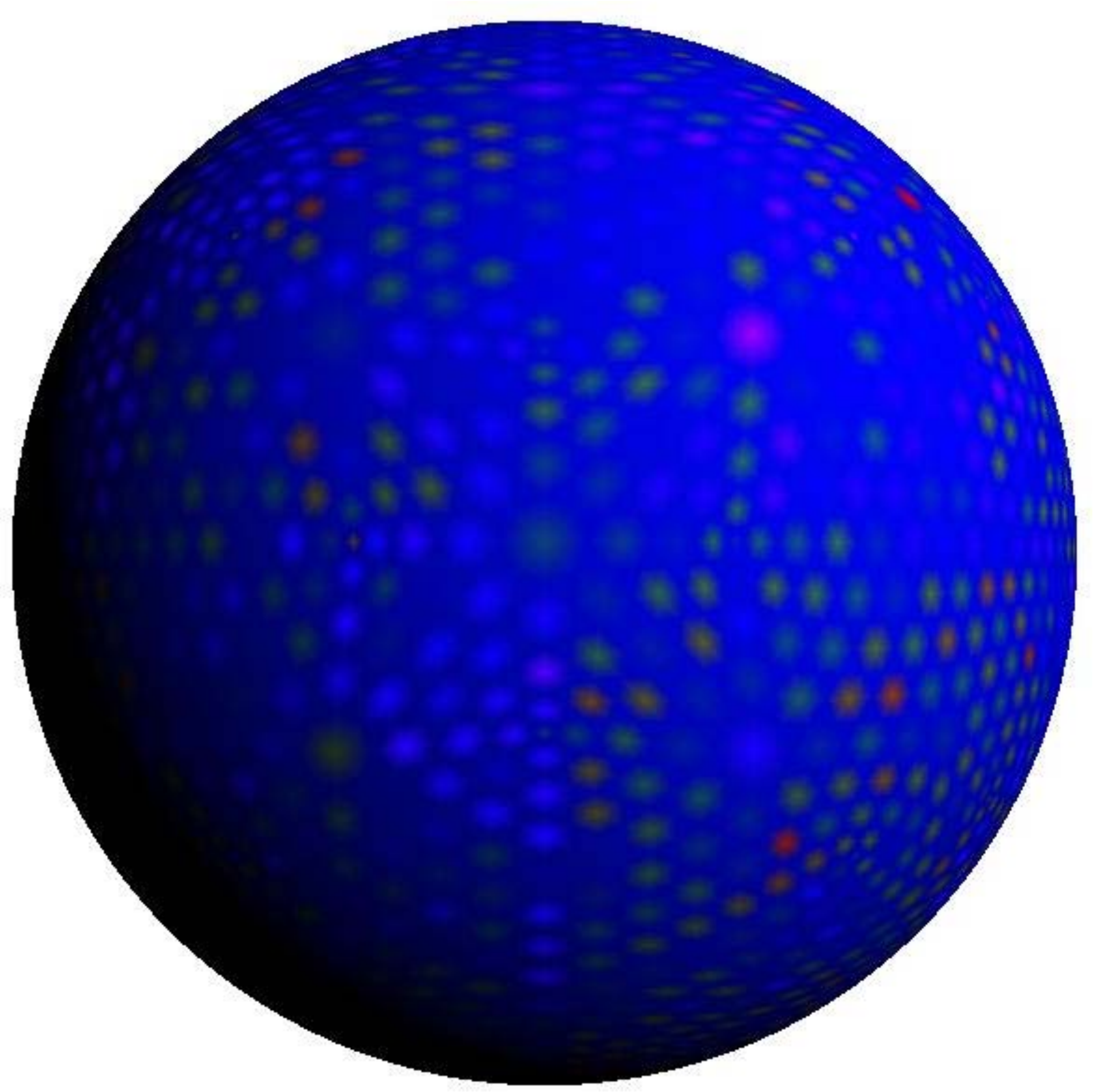}  \\
($a'$) &  ($b'$) &   ($c'$)  \\
\includegraphics[scale=0.2]{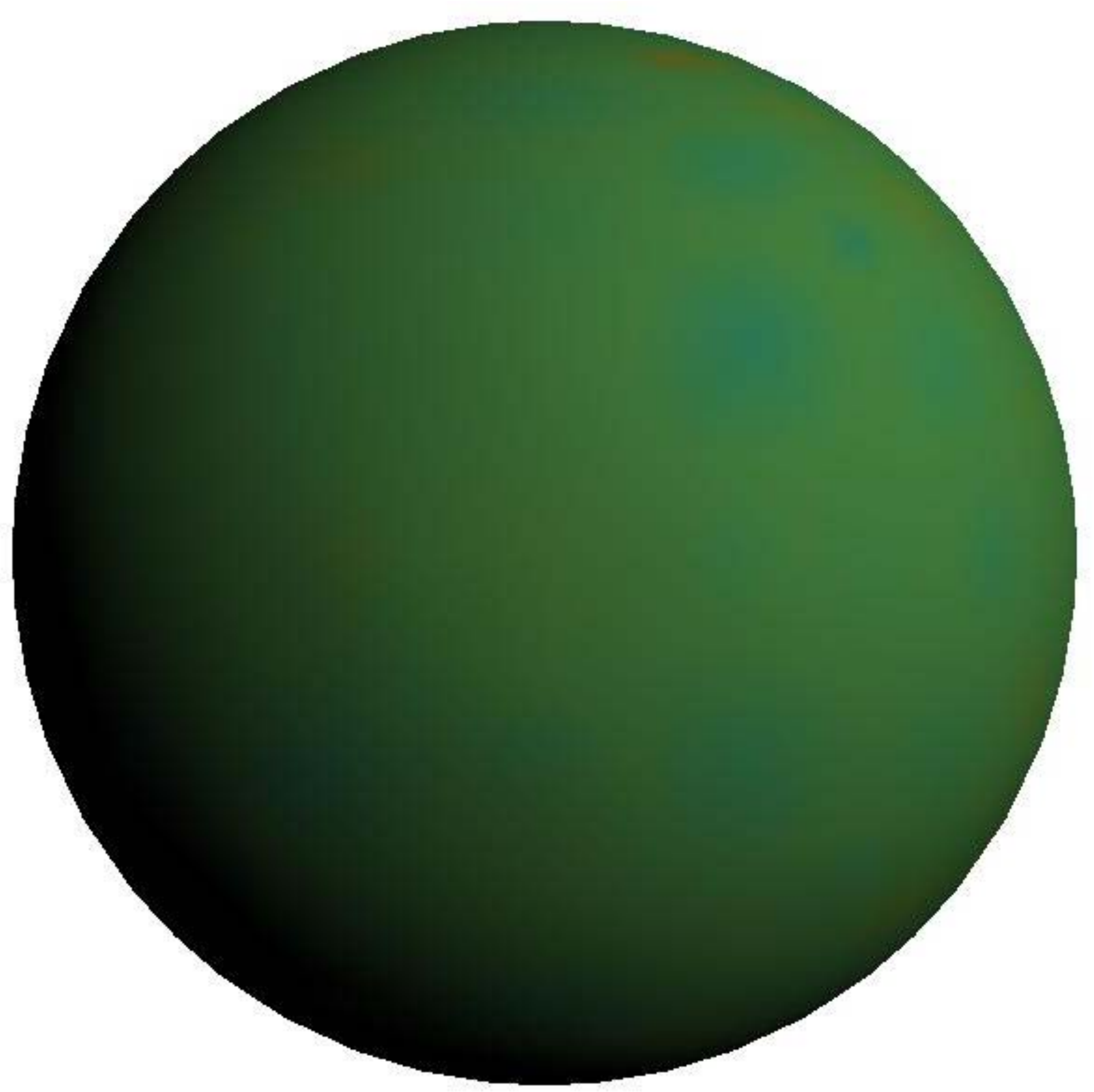}  & \ \ \ \ \
\includegraphics[scale=0.2]{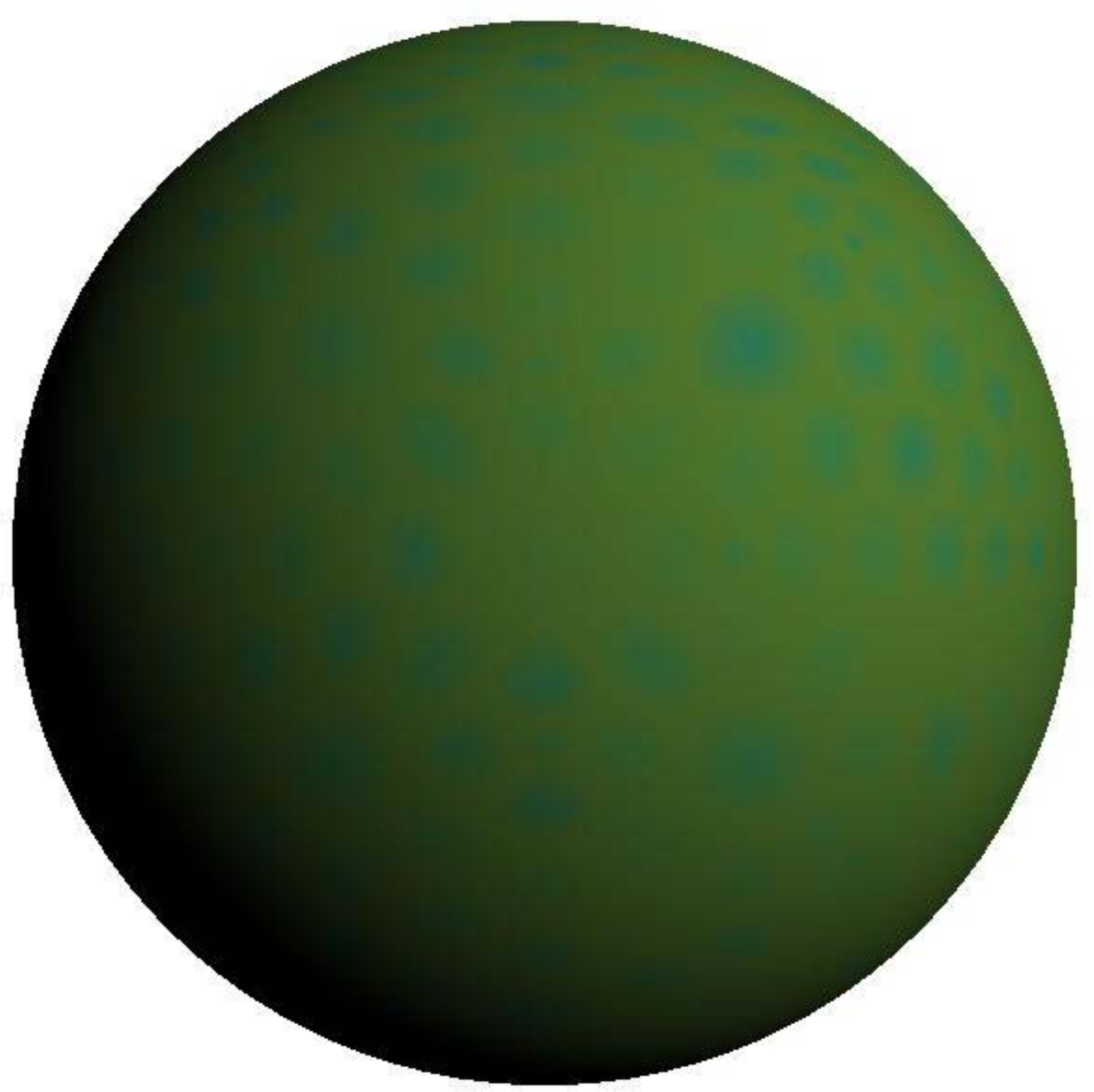}  & \ \ \ \ \
\includegraphics[scale=0.2]{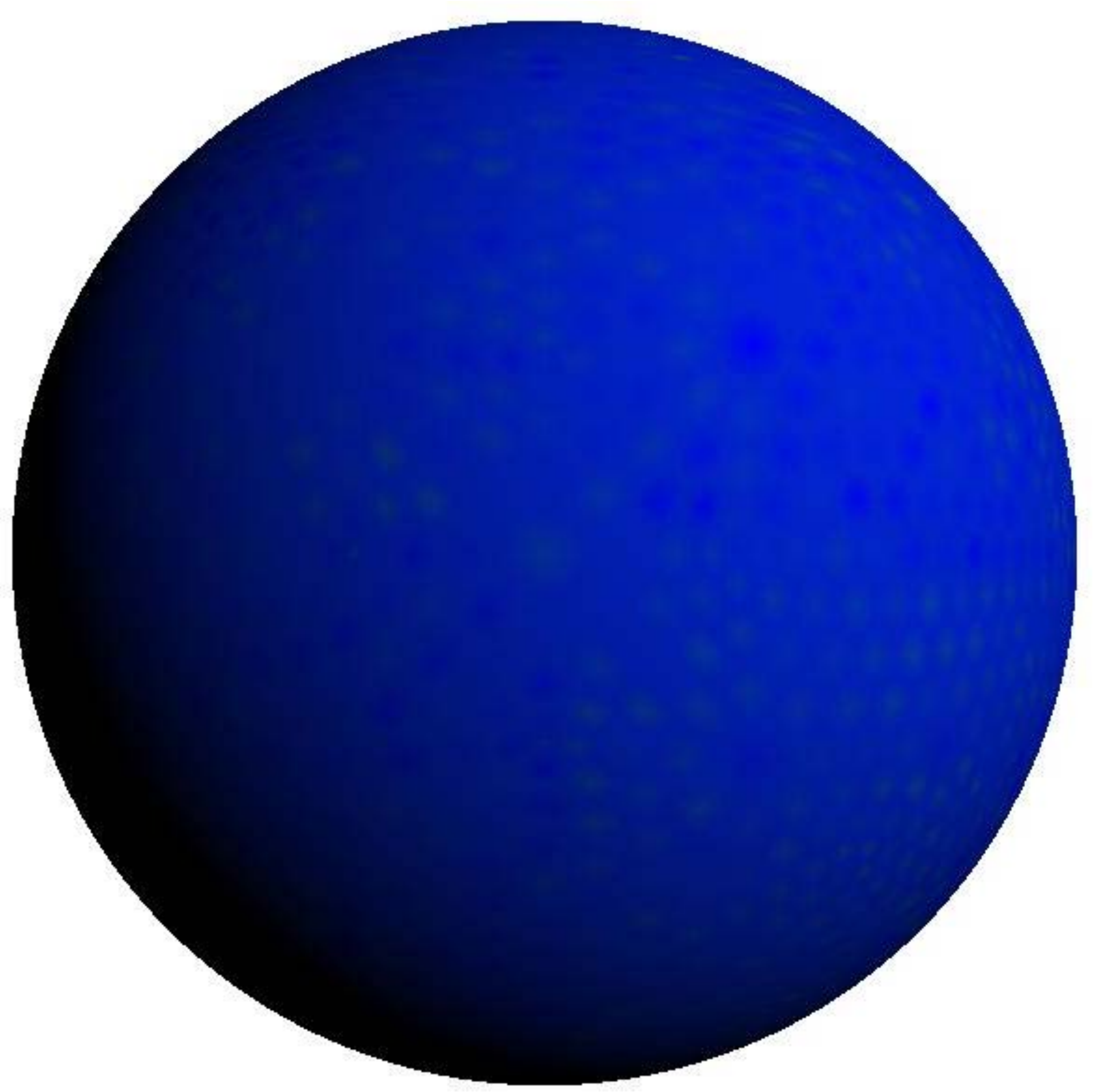}  \\
($a''$) &  ($b''$) & ($c''$)  \\
\end{tabular}
}  \caption{{\small Surface ${\cal S}_4$. ($a$), ($b$)
and ($c$) are three control meshes where one time refinement is
implemented from ($a$) to ($b$), and ($b$) to ($c$). The corresponding distribution of the error
$u-u^h$ resulting from FEM-Linear and IGA-Loop is respectively shown
in ($a'$), ($b'$), and ($c'$) of the second row, and ($a''$),
($b''$), and ($c''$) of the third row.}} \label{wholesphere}
\end{figure}

As shown in the top row of Figure \ref{wholesphere}, three different control meshes from coarse
to refined are provided. The total numbers of vertices/patches are 706/1408, 2818/5632, and 11266/22528,
respectively. They have the same limit surface ${\cal S}_4$ by the extended Loop subdivision.
We solve the Laplace-Beltrami triharmonic equation based on the proposed IGA-Loop.
The computational accuracy is shown by the error distribution $u-u^h$.
As shown in Figure \ref{wholesphere}, the accuracy of IGA-Loop is higher than that of FEM-Linear.

\begin{figure}[htb!] \centerline {
\begin{tabular}{c}
\includegraphics[scale=0.7]{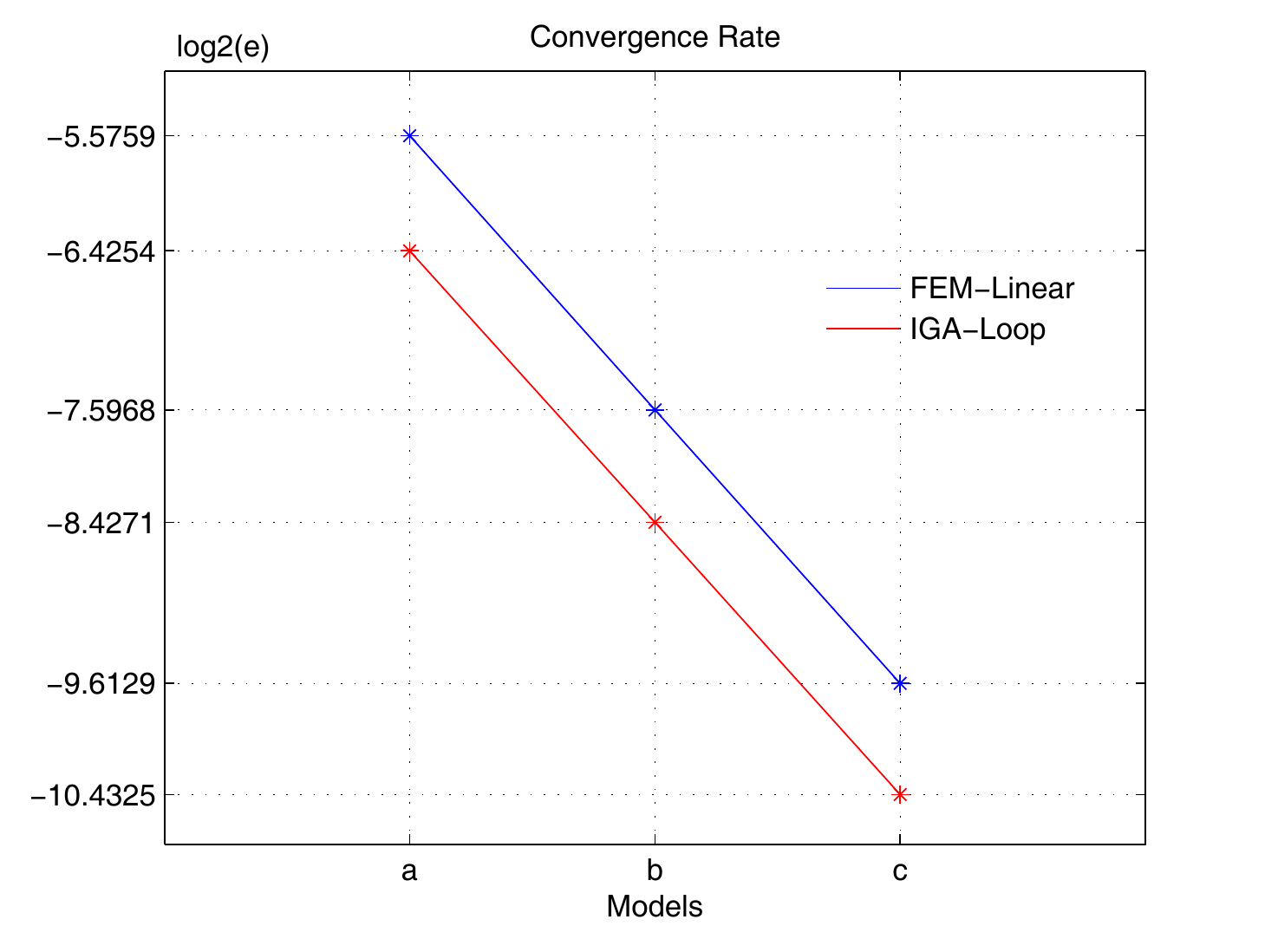}
\end{tabular}
} \caption{{\small Surface ${\cal S}_4$. Comparison of the convergence rate of
 the errors versus the refinement times between FEM-Linear and IGA-Loop. Here the numbers a, b and c on the $x$-axis correspond to the models of Figure \ref{wholesphere} ($a$), ($b$)  and ($c$) respectively, and $e$ on the $y$-axis is the $L^2$-norm error.}} \label{errsphere}
\end{figure}

Different from the above two tests, the test model for the sixth-order problem is on a closed surface. We further observe the convergence rate
of IGA-Loop against successively finer meshes of the unit sphere in Figure \ref{errsphere}. The approximate error (\ref{surffem}) in norm $L^2({\cal S})$  holds and the same convergence rate 2 is obtained. Furthermore, we can observe the approximation error obtained with the FEM-Linear discretization is approximately 1.8 times more than that of the IGA-Loop discretization. Therefore, IGA-Loop obtains the same level of accuracy more efficiently than FEM-Linear.

\subsection{Computational Cost}

Actually, the solver of our IGA-Loop can be merged into the framework of FEM-Linear, however the support of each basis function of IGA-Loop is two-ring neighbors described in Section 3. In this section, for the four test examples of Figures \ref{quacy}, \ref{quasph}, \ref{wholecylinder} and \ref{wholesphere}, we list the corresponding time costs of computing the basis functions in Tables \ref{time2th} and \ref{time46th}. To show that the extended Loop subdivision scheme does not require structured meshes and can support the same meshes with any topological structure as the standard finite elements,
the valences of the control vertices in these surfaces are in the range of 3 to 12 depicted in Figure \ref{quacy}, \ref{quasph}, \ref{wholecylinder} and \ref{wholesphere}.

Table \ref{time2th} corresponds to the examples of the Laplace-Beltrami harmonic problem where the number of control vertices and triangular patches is shown in the first and the third columns, the second and the fourth columns list the time cost (in seconds) of computing the basis functions and their derivatives because they can be pre-computed and saved in a data structure. Once refinement makes the number of triangular patches on the refined meshes increases four times and their sizes approximately decrease by half. The computational cost of IGA-Loop with the same control meshes is larger than FEM-Linear. One reason is the computation of the derivatives of the linear basis functions is unnecessary, as it happens instead with IGA-Loop. The other reason is the computation around the boundary patches is more complex for our IGA-Loop than FEM-Linear because the extended Loop subdivision scheme embodies the angle information different of every patch.

The data for the examples of the Laplace-Beltrami biharmonic problem and the Laplace-Beltrami triharmonic problem are listed in Table \ref{time46th}. Similar results to the Laplace-Beltrami harmonic problem can be observed. Here we notice the data of the Laplace-Beltrami triharmonic problem which is listed in the third and the fourth columns. The example is a closed sphere, the computation cost of IGA-Loop is almost equal to FEM-Linear. We will explain the phenomenon. Firstly the computation of the basis functions around boundary patches for IGA-Loop is removed, and the rest is only the case of computing on the interior patches. Secondly, interior patches share the same set of basis functions which depend only on the valence list of their control vertices, so we can merge the patches according to their valence list, and then use Stam's fast evaluation to treat them. Finally, as the mesh refinement proceeds, the added patches are ordinary whose valences of the three control vertices are six.

\begin{table}[htb!]
\setlength{\belowcaptionskip}{6pt}
\caption{Data of Test Suite 1} \centerline {
\begin{tabular}{|c|c|c|c|}
 \hline \  vertices/patches \  & \ \ basis func.(s)\ \  &  \  vertices/patches \  & \ \ basis func.(s) \ \ \\
\hline Figure \ref{quacy} &   {\scriptsize FEM-Linear} \ {\scriptsize IGA-Loop} & Figure \ref{quasph}
 &  {\scriptsize FEM-Linear} \ {\scriptsize IGA-Loop}  \\
\hline
221/384    &   0.02 \ \ \ \ \ \ 0.04     & 107/176     &   0.01 \ \ \ \ \ \ 0.05    \\
825/1536   &   0.06   \ \ \ \ \ \ 0.11   & 389/704     &   0.03   \ \ \ \ \ \ 0.09    \\
3185/6144  & 0.11 \ \ \ \ \ \ 0.20       &  1481/2816  & 0.05  \ \ \ \ \ \ 0.14      \\
\hline
\end{tabular}
} \label{time2th}
\end{table}

\begin{table}[ht!]
\setlength{\belowcaptionskip}{6pt}
 \caption{Data of Test Suite 2 and Test Suite 3} \centerline {
\begin{tabular}{|c|c|c|c|}
 \hline \  vertices/patches \  & \ \ basis func.(s)\ \  &  \  vertices/patches \  & \ \ basis func.(s) \ \ \\
\hline  Figure \ref{wholecylinder}   &   {\scriptsize FEM-Linear} \ {\scriptsize IGA-Loop} &  Figure \ref{wholesphere}
 &  {\scriptsize FEM-Linear} \ {\scriptsize IGA-Loop}  \\
\hline
432/768    &   0.04 \ \ \ \ \ \ 0.10    & 706/1408    &   0.02 \ \ \ \ \ \ 0.02  \\
1632/3072  &   0.11   \ \ \ \ \ \ 0.21   &   2818/5632  &   0.04   \ \ \ \ \ \ 0.05 \\
6336/12288 & 0.23 \ \ \ \ \ \ 0.42    &  11266/22528 & 0.08 \ \ \ \ \ \ 0.10    \\
\hline
\end{tabular}
} \label{time46th}
\end{table}

\section{Conclusion}
\label{conclu}

The isogeometric analysis for surface PDEs based on the extended Loop subdivision is presented.
The set of quartic box-splines corresponding to each subdivided control mesh are utilized to
represent the geometry of interest, and to construct the solution space for dependent variables as well.
The finite elements induced by the extended Loop subdivision possess the ability to
represent the geometry exactly which fully agrees with the concept of isogeometric analysis.

We have evaluated the performance of the proposed IGA-Loop in dealing with the various surface
PDEs, including the Laplace-Beltrami harmonic/biharmonic/triharmonic equations on different limit surfaces.
As shown in the visual and quantitative results, the proposed method can yield desirable
performance in terms of the accuracy, convergence and computational cost for solving the above
classical surface PDEs. As demonstrated with both the theoretical and numerical results,
the proposed IGA-Loop approach is proved to be second-order accuracy in the sense of $L^2$  norm of the surface.
Through various comparisons, the performance of IGA-Loop is also outperformed over
the standard linear finite elements. The solution of these problems will help to establish a computational framework for the
geometric PDEs solved by IGA-Loop in our future work.

\section*{Acknowledgments}

Qing Pan is supported by National Natural Science Foundation of China (NSFC) (No.11671130),
Hunan Provincial Natural Science Foundation of China (No.2018JJ2248, 2019JJ50395), and the Open Project
Program of the National Laboratory of Pattern Recognition of China (NLPR). Chong Chen is supported
by the Key Research Project of Beijing National Natural Science Foundation (No.Z180002).
Gang Xu was supported by NSFC (No.61761136010, 61772163) and Zhejiang Provincial Natural
Science Foundation (No.LR16F020003).

\end{document}